\documentclass[10pt]{article}
\usepackage{songmei}
\usepackage[toc,page]{appendix}
\usepackage{subfigure}
\usepackage{csquotes}

\numberwithin{equation}{section}
\numberwithin{theorem}{section}

\def\vphi{\varphi}

\def\op{{\rm op}}

\def\bX{{\boldsymbol X}}
\def\bK{{\boldsymbol K}}

\def\bi{{\boldsymbol i}}

\def\cC{{\mathcal C}}

\def\sF{{\mathsf{F}}}

\def\bp{{\boldsymbol p}}

\def\de{{\rm d}}
\def\bx{{\boldsymbol x}}

\def\bW{{\boldsymbol W}}

\def\cF{{\mathcal F}}
\def\Unif{{\sf Unif}}

\def\bU{{\boldsymbol U}}
\def\bv{{\boldsymbol V}}
\def\bM{{\boldsymbol M}}
\def\bZ{{\boldsymbol Z}}

\def\bI{{\boldsymbol I}}
\def\bz{{\boldsymbol z}}
\def\proj{{\mathsf P}}

\def\cE{{\mathcal E}}

\def\bB{{\boldsymbol B}}

\def\bu{{\boldsymbol u}}
\def\bg{{\boldsymbol g}}

\def\bA{{\boldsymbol A}}
\def\bD{{\boldsymbol D}}

\def\bv{{\boldsymbol v}}

\def\reals{{\mathbb R}}

\def\de{{\rm d}}
\def\bx{{\boldsymbol x}}

\def\bW{{\boldsymbol W}}
\def\bG{{\boldsymbol G}}

\def\cF{{\mathcal F}}
\def\Unif{{\rm Unif}}
\def\bU{{\boldsymbol U}}
\def\bv{{\boldsymbol v}}
\def\bz{{\boldsymbol z}}
\def\proj{{\mathsf P}}

\def\cE{{\mathcal E}}
\def\bt{{\boldsymbol t}}

\def\cD{{\mathcal D}}
\def\arctanh{\operatorname{arctanh}}

\def\supp{\operatorname{supp}}

\def\bu{{\boldsymbol u}}
\def\bg{{\boldsymbol g}}

\def\bh{{\boldsymbol h}}

\def\bA{{\boldsymbol A}}
\def\bH{{\boldsymbol H}}
\def\cS{{\mathcal S}}
\def\cM{{\mathcal M}}

\def\cX{{\mathcal X}}

\def\bY{{\boldsymbol Y}}

\def\GOE{{\rm GOE}}
\def\TAP{{\rm TAP}}
\def\bbm{{\boldsymbol m}}
\def\sh{{\mathsf h}}

\def\cB{{\mathcal B}}

\newcommand{\indic}[1]{\mathbf{1}_{#1}}
\def\Vol{\operatorname{Vol}}
\def\Bayes{\mathrm{Bayes}}
\def\Tr{\operatorname{Tr}}

\def\TAP{\mathrm{TAP}}
\def\AMP{\mathsf{AMP}}
\def\VB{\mathrm{VB}}
\def\Sym{\mathrm{Sym}}

\def\sD{\mathsf{D}}
\def\bfA{\mathbf{A}}

\usepackage{hyperref}
\hypersetup{
    colorlinks,
    linkcolor={blue!80!black},
    citecolor={green!50!black},
}
\colorlet{linkequation}{blue}

\newif\iflong
\longfalse 

\begin{document}

\title{Local convexity of the TAP free energy and AMP convergence for $\Z_2$-synchronization}

\author{Michael Celentano\thanks{Department of Statistics, University of California,
Berkeley. E-mail: \texttt{mcelentano@berkeley.edu}} \and
Zhou Fan\thanks{Department of Statistics and Data Science, Yale
University. E-mail: \texttt{zhou.fan@yale.edu}} \and
Song Mei\thanks{Department of Statistics, University of California,
Berkeley. E-mail: \texttt{songmei@berkeley.edu}}}
\date{}

\maketitle

\begin{abstract}
  We study mean-field variational Bayesian inference using the TAP approach, for
    $\Z_2$-synchronization as a prototypical example of a high-dimensional Bayesian
    model. We show that for any signal strength $\lambda>1$ (the weak-recovery
    threshold), there exists a unique local minimizer of the TAP free energy functional 
    near the mean of the Bayes posterior law. Furthermore, the TAP free energy
in a local neighborhood of this minimizer is strongly convex. Consequently, a
natural-gradient/mirror-descent algorithm achieves linear convergence to this
minimizer from a local initialization, which may be obtained by a constant number
of iterations of Approximate Message Passing (AMP). This provides a rigorous foundation for variational inference in high dimensions via minimization of the TAP free energy.

    We also analyze the finite-sample convergence of AMP, showing that AMP
    is asymptotically stable at the TAP minimizer for any $\lambda>1$, and is
    linearly convergent to this minimizer from a spectral initialization
    for sufficiently large $\lambda$. Such a guarantee is stronger than
    results obtainable by state evolution analyses, which only describe a fixed
    number of AMP iterations in the infinite-sample limit.

    Our proofs combine the Kac-Rice formula and Sudakov-Fernique Gaussian comparison inequality to analyze the complexity of critical points that satisfy strong convexity and stability conditions within their local neighborhoods.
\end{abstract}

\tableofcontents

\section{Introduction}

Variational inference is an increasingly popular method for performing
approximate Bayesian inference, and is widely used in applications
ranging from document classification to population
genetics \cite{blei2003latent,liang2007infinite,carbonetto2012scalable,raj2014faststructure}.
For large-scale problems, variational methods provide an appealing alternative
to Markov Chain Monte Carlo procedures, particularly in settings where MCMC
may be computationally prohibitive to apply. We refer readers to the classical expositions
\cite{jordan1999introduction, wainwright2008graphical} and the recent review
\cite{blei2017variational} for an introduction.

In ``mean-field'' models where the posterior
distribution $p(\bx|\bY)$ of parameters $\bx$ given data $\bY$ may be
close to being a product measure, a common approach to variational
inference is to approximate $p(\bx|\bY)$ by a product law. The most widely used
such approximation minimizes the KL-divergence to
$p(\bx|\bY)$ over the class $\mathcal{Q}$ of product measures,
\begin{equation}\label{eq:naivemeanfield}
\hat{q}(\bx)=\argmin_{q \in \mathcal{Q}}
\sD_{\mathrm{KL}}(q(\bx)\|p(\bx|\bY)).
\end{equation}
When $\bx \in \R^n$ is high-dimensional, a problematic phenomenon may occur in
which this distribution $\hat{q}(\bx)$ provides inconsistent
approximations to the posterior marginals and posterior means, even in models
where all low-dimensional marginals of $p(\bx|\bY)$ have approximately
independent coordinates. Such a phenomenon was first investigated
by Thouless, Anderson, and Palmer for the
Sherrington-Kirkpatrick (SK) model of spin glasses,
where a simple method of addressing this inaccuracy---now often called the ``TAP
correction''---was also proposed \cite{thouless1977solution}.
Manifestations of this phenomenon and
analogues of the TAP free energy for several high-dimensional statistical models have been studied in
\cite{krzakala2014variational,rangan2016inference,ghorbani2019instability,fan2021tap,qiu2022tap}, and we provide further
discussion in Section \ref{subsec:literature}.

The TAP approach to variational inference constructs a free
energy functional $\cF_\TAP$ by adding a correction term to the KL-divergence
objective (\ref{eq:naivemeanfield}). This TAP correction
accounts for dependences between pairs of coordinates of $\bx$ in their
posterior law, which are individually weak but may have a non-negligible
aggregate effect in high dimensions. Variational inference is
performed by minimizing $\cF_\TAP$, or by solving the TAP stationary equations
\begin{equation}\label{eq:TAPequations}
0=\nabla \cF_\TAP.
\end{equation}
Since the pioneering work of
\cite{kabashima2003cdma,donoho2009message,bolthausen2014iterative}, both the theory and implementation of
TAP-variational inference have been closely connected to Approximate
Message Passing (AMP) algorithms, which provide specific iterative procedures
for solving (\ref{eq:TAPequations}). TAP-variational inference has been
successfully applied via AMP to a variety of high-dimensional statistical
problems. We highlight in particular the line of work
\cite{rangan2012iterative,deshpande2014information,montanari2015non,lesieur2015phase,barbier2016mutual,lelarge2019fundamental,montanari2021estimation}
on low-rank matrix estimation, of which the $\Z_2$-synchronization
problem is a specific example.

The goal of our current paper is to address several foundational questions
regarding TAP-variational inference that, despite the above successes, remain 
poorly understood. First, the convergence of AMP is usually known only in a weak
sense, guaranteeing $\|\sqrt{n} \cdot \nabla \cF_\TAP\|_2^2<\eps$ in the limit
$n \to \infty$ for a constant number of AMP iterations $k \equiv k(\eps)$
independent of $n$. Such a guarantee is too weak to ensure, for example, even
the high-probability existence of a critical point of $\cF_\TAP$ to which the
AMP iterates converge. It does not establish whether the minimizer of
$\cF_\TAP$ is close to the true Bayes posterior mean, and indeed, these properties
remain conjectural in most models to which the AMP/TAP approach has
been applied. Second, regularity properties of the landscape of $\cF_\TAP$ are
largely unknown, making it unclear whether optimization algorithms other than
AMP can successfully implement the TAP-variational inference paradigm.

In this paper, we clarify these properties of $\cF_\TAP$ and the convergence of
AMP and other descent algorithms for the specific model of
$\Z_2$-synchronization. We build upon previous results and
techniques of \cite{fan2021tap}, which studied this model in a regime of
large signal-to-noise. Our main results will show that
for any signal strength above the weak-recovery threshold,
there exists a unique local minimizer $\bbm_\star$ of $\cF_\TAP$ near the Bayes
posterior mean, and $\cF_\TAP$ is strongly convex in a local neighborhood (of
non-trivial size) around $\bbm_\star$. Consequently,
a generic natural-gradient-descent (NGD) algorithm exhibits linear convergence
to $\bbm_\star$ from a local initialization, which may be obtained by a finite
number of iterations of AMP. We also
show that the Jacobian of the AMP map is stable at $\bbm_\star$, so that AMP
initialized in a (potentially very) small neighborhood of $\bbm_\star$ will also
converge for fixed $n$ as the number of iterations $t \to \infty$. In the large signal-to-noise regime of \cite{fan2021tap}, we show that both NGD and
AMP exhibit linear convergence to $\bbm_\star$ from a spectral initialization.

Formalizing these properties of $\cF_\TAP$ and the convergence of
generic optimization algorithms has several appeals over the existing theory
around AMP. First, it clarifies a concrete objective function for
high-dimensional variational inference, which can serve a number of practical
purposes such as assessing algorithm convergence. Second, the convergence and
state evolution of AMP are tied to probabilistic aspects of the model,
whereas NGD is always a strict descent algorithm (for small enough step size,
even in misspecified models) and may provide a more flexible and robust
approach for optimization in practice. Finally, understanding the
landscape of $\cF_\TAP$ may be useful in other contexts. For example, following
the initial posting of our work,
\cite{alaoui2022sampling,celentano2022sudakov} have used the local strong
convexity of $\cF_\TAP$ in the related SK model to argue that its
stationary point is Lipschitz in the external field. This is a central technical
ingredient in these works to show the correctness of an algorithmic stochastic
localization procedure for sampling from the SK measure.

We review relevant background on the $\Z_2$-synchronization model in
Section \ref{subsec:Z2}, and we describe our results in more detail
in Section \ref{subsec:contributions}.

\subsection{$\Z_2$-synchronization and the TAP free energy}\label{subsec:Z2}

In $\Z_2$-synchronization, we wish to estimate an unknown binary vector
$\bx \in \{-1,+1\}^n$ having the entry-wise symmetric Bernoulli prior $x_i
\overset{iid}{\sim} \Unif\{-1,+1\}$.
For a signal-to-noise parameter $\lambda>0$, we observe
\begin{equation}\label{eq:Z2-sync}
	\bY = \frac{\lambda}{n} \bx \bx^\sT + \bW,\;\; \text{where}\;\; \bW \sim \GOE(n).
\end{equation}
Thus $\bW$ is symmetric Gaussian noise, having entries $(w_{ii}:i=1,\ldots,n)
\overset{iid}{\sim} \cN(0,2/n)$ independent of $(w_{ij}:1 \leq i<j \leq n)
\overset{iid}{\sim} \cN(0,1/n)$.
Equivalently, $\bW=(\bZ+\bZ^\top)/\sqrt{2n}$ where $(z_{ij}:i,j=1,\ldots,n)
\overset{iid}{\sim} \cN(0,1)$.

The parameter $\bx$ is identifiable only up to $\pm$ sign, and the posterior law
$p(\bx|\bY)$ has the corresponding sign symmetry $p(\bx|\bY)=p(-\bx|\bY)$. Thus
we will
consider estimation of the sign-invariant rank-one matrix $\bX=\bx\bx^\top$.
The Bayes posterior-mean estimate of this matrix is
\begin{equation}\label{eq:XBayes}
\widehat{\bX}_{\Bayes}=\E[\bx\bx^\sT \mid \bY].
\end{equation}
The asymptotic squared-error Bayes risk of this estimator
was characterized by Deshpande, Abbe, and Montanari
in \cite{deshpande2016asymptotic}:
\begin{equation}\label{eq:Bayesrisk}
\lim_{n \to \infty} \frac{1}{n^2}
\E[\|\widehat{\bX}_\Bayes-\bx\bx^\top\|_\sF^2]
=\begin{cases} 1-q_*(\lambda)^2 & \text{ if } \lambda>1 \\
1 & \text{ if } \lambda \leq 1,\end{cases}
\end{equation}
where $q_\star(\lambda)>0$ is the solution to a fixed-point equation
(\ref{eq:fixedpointq}). Thus for $\lambda<1$, no non-trivial estimation is
possible in the large-$n$ limit, as the optimal Bayes risk coincides with that 
of the trivial estimator $\widehat \bX = \bzero$. In contrast, for
$\lambda>1$, the Bayes estimator achieves positive entry-wise
correlation with $\bx\bx^\top$.

\cite{deshpande2016asymptotic} studied also an AMP algorithm for approximately
computing $\widehat{\bX}_\Bayes$. Starting from initializations
$\bh^0,\bbm^{-1} \in \R^n$, this algorithm takes the form
\begin{align}
\bbm^k&=\tanh(\bh^k)\nonumber\\
\bh^{k+1}&=\lambda \bY\bbm^k-\lambda^2[1-Q(\bbm^k)]\bbm^{k-1}\label{eq:AMP}\tag{{\sf AMP}}
\end{align}
where $Q(\bbm)=\|\bbm\|_2^2/n$. The analyses of \cite{deshpande2016asymptotic}
imply that for any $\lambda>1$ and $\eps>0$, starting from an informative
initialization $\bh^0$, there exists an iterate $k \equiv k(\lambda,\eps)$
of AMP for
which $\|\bbm^k(\bbm^k)^\top-\widehat{\bX}_\Bayes\|_{\sF}^2/n^2<\eps$,
with high probability for all large $n$. More recent results of
\cite{montanari2021estimation} imply that such a guarantee holds also for AMP
with a spectral initialization.

The TAP free energy in this $\Z_2$-synchronization model is defined 
for $\bbm \in (-1,1)^n$ by
\begin{equation}\label{eq:TAP}\tag{{\sf TAP}}
\cF_{\TAP}(\bbm) ={-}\frac{\lambda}{2 n} \< \bbm, \bY \bbm\> - \frac{1}{n}
\sum_{i=1}^n \sh( m_i)  - \frac{\lambda^2}{4}[1 - Q(\bbm)]^2
\end{equation}
where $Q(\bbm)=\|\bbm\|_2^2/n$ as above, and $\sh(m)$ is the binary
entropy function
\begin{equation}\label{eq:entropy}
\sh(m)=-\frac{1+m}{2}\log \frac{1+m}{2}-\frac{1-m}{2}\log\frac{1-m}{2}.
\end{equation}
This function $\cF_\TAP$ has the sign symmetry $\cF_\TAP(\bbm)=\cF_\TAP(-\bbm)$,
corresponding to the above sign symmetry of the posterior law.
The first two terms of (\ref{eq:TAP}) coincide\footnote{Up to an additive
constant, and a replacement of $\E_{\bx \sim q}[\langle \bx,\bY\bx \rangle]$
by $\langle \bbm,\bY\bbm \rangle$ which incurs negligible error} with
the KL-divergence $\sD_{\mathrm{KL}}(q(\bx)\|p(\bx|\bY))$ for a product measure 
$q(\bx)$ on $\{-1,+1\}^n$, upon parameterizing $q$ by its mean $\bbm=\E_{\bx \sim
q}[\bx] \in (-1,1)^n$. The third term of (\ref{eq:TAP}) is the TAP correction. Applying $\sh'(m)={-}\arctanh(m)$, the stationary
condition $0=\nabla \cF_\TAP(\bbm)$ may be rearranged as the TAP mean-field
equations
\[\bbm=\tanh\Big(\lambda \bY\bbm-\lambda^2[1-Q(\bbm)]\bbm\Big),\]
and the AMP algorithm (\ref{eq:AMP}) is an iterative scheme
for computing a fixed point of these equations.

In \cite{fan2021tap}, an upper bound for the expected number of critical points
of $\cF_\TAP$ in sub-regions of the domain $(-1,1)^n$ was derived for
any $\lambda>0$. Using this result, for $\lambda>\lambda_0$ a large enough
absolute constant, it was shown that the global minimizer
$\bbm_\star$ of $\cF_\TAP$ satisfies
$\E[\|\bbm_\star\bbm_\star^\top-\widehat{\bX}_\Bayes\|_\sF^2]/n^2 \to 0$,
and that this holds more generally for any critical point $\bbm$ of $\cF_\TAP$
in the domain
\[\cS=\{\bbm \in (-1,1)^n:\cF_\TAP(\bbm)<-\lambda^2/3\}.\]
As a consequence, it was also shown that
$\E[\|\bbm_\star\bbm_\star^\top-\widehat{\bX}_\Bayes\|_\sF^2]/n^2$ must be
bounded away from 0 for the minimizer $\bbm_\star$ of the naive mean-field
objective (\ref{eq:naivemeanfield}) parametrized similarly by $\bbm$.
We note that the landscape guarantees in \cite{fan2021tap} do not
extend to the entire weak-recovery regime $\lambda>1$. The analyses for large
$\lambda>\lambda_0$ also fall short of showing uniqueness (up to sign) of the
TAP critical point $\bbm_\star$ in $\cS$, and of establishing
polynomial-time convergence of AMP or other optimization algorithms for computing $\bbm_\star$.

\subsection{Contributions}\label{subsec:contributions}

Our current work establishes the following properties of $\cF_\TAP(\bbm)$ and of
descent algorithms for minimizing this objective function.

\begin{enumerate}

\item \textbf{Existence of Bayes-optimal TAP local minimizer.}
For any $\lambda > 1$, we show there exists a local minimizer
$\bbm_\star$ of $\cF_\TAP$ such that $\| \bbm_\star \bbm_\star^\top - \widehat
\bX_{\Bayes}\|_{\sF}^2/n^2 \rightarrow 0$ in probability. This
strengthens the guarantee of \cite{fan2021tap} that was shown for large
$\lambda>\lambda_0$. Subject to the validity of a numerical conjecture about a
deterministic low-dimensional variational problem
(see Remark \ref{remark:globalmin}), our results imply that this is also
the global minimizer of $\cF_\TAP$ for any $\lambda>1$.
\item \textbf{Local strong convexity of the TAP free energy.}
For any $\lambda>1$, we show that $\cF_\TAP$ is strongly convex in a
$\sqrt{\eps n}$-neighborhood of this local minimizer $\bbm_\star$. Hence this
local minimizer is the unique critical point satisfying
$\| \bbm_\star \bbm_\star^\top - \widehat
\bX_{\Bayes}\|_{\sF}^2/n^2<\iota(\eps)$, for some constant $\iota(\eps)>0$.
\item \textbf{Local convergence of natural gradient descent.}
We introduce a natural gradient descent (NGD) algorithm for minimizing
$\cF_\TAP$, which is equivalently a mirror descent procedure that adapts to the
curvature of $\cF_\TAP$ near the boundaries of $(-1,1)^n$. For any $\lambda>1$,
we prove that NGD achieves linear convergence to $\bbm_\star$ from an 
initialization within this $\sqrt{\eps n}$-neighborhood. This initialization
may be obtained by first performing a fixed number of iterations of AMP,
thus yielding a polynomial-time algorithm for computing $\bbm_\star$.
\item \textbf{Stability of AMP.}
For any $\lambda > 1$, we show that the AMP map is stable at $\bbm_\star$,
in the sense of having a Jacobian with spectral radius strictly less than 1. 
Thus, AMP initialized in a sufficiently small neighborhood of $\bbm_\star$ will
also linearly converge to $\bbm_\star$.
\item \textbf{Finite-$n$ convergence of AMP and NGD.} Finally, for
$\lambda>\lambda_0$ a large enough absolute constant, our results combine with
those of \cite{fan2021tap} to show that $\bbm_\star$ is the global
minimizer and unique critical point (up to sign) of $\cF_\TAP$ in the domain
$\{\bbm:\cF_\TAP(\bbm)<-\lambda^2/3\}$.
In this signal-to-noise regime, we prove that both AMP and NGD alone exhibit
linear convergence to $\bbm_\star$ from a spectral initialization.

We emphasize that this convergence of AMP is established in the sense
$\lim_{k \to \infty} \bbm^k=\bbm_\star$ for fixed dimension $n$, which is
stronger than the guarantee $\limsup_{n \to \infty}
\|\bbm^k-\bbm_\star\|_2^2/n<\eps$ for fixed $k \equiv k(\lambda,\eps)$ that
is obtainable by standard analyses of the AMP state evolution.
\end{enumerate}

The main challenge in understanding the landscape of $\cF_\TAP$ locally
near $\bbm_\star$ is that---for any constant signal strength $\lambda$---this
point $\bbm_\star$ does not converge to the true signal vector 
$\bx \in\{-1,+1\}^n$ as $n \to \infty$, but rather remains random in
$(-1,1)^n$. Thus it is not enough to study the landscape of $\cF_\TAP$ in a
vanishing neighborhood of $\bx$ using, for example, the uniform convergence arguments \cite{sun2018geometric, mei2018landscape}. The above results instead pertain to the 
geometry of $\cF_\TAP$ in a random region of the cube $(-1,1)^n$.

We will prove these results using a combination of the Kac-Rice formula and
Gaussian comparison inequalities. We provide a detailed overview of this
proof in Section \ref{sec:local}. The Kac-Rice formula has been successfully
applied to study the complexity of critical points for various non-convex
function landscapes. However, to our knowledge, our argument for using Kac-Rice
to study also the local geometry around a particular critical point is novel.
We believe that this technique may be of independent interest for some recent
analyses of related disordered
systems \cite{bolthausen2018morita,ding2019capacity,fan2021replica}, where
conditioning on a sequence of AMP iterates was used as a surrogate for
conditioning on an actual TAP critical point.

\subsection{Further related literature}\label{subsec:literature}

\subsubsection{Variational inference}

The terminology ``variational inference'' encompasses a large
family of methods for approximate Bayesian inference
\cite{blei2012probabilistic,pearl1982reverend,minka2001family,yedidia2003understanding}, based upon approximating
a variational representation to the evidence or marginal log-likelihood of the
observed data. Variational inference has been incorporated into many software packages including Pyro
\cite{bingham2019pyro}, Infer.NET \cite{minka2014infer}, and Edward
\cite{tran2016edward}.

There has been renewed interest in theoretical analyses of variational inference
in recent years, focusing on a number of common desiderata:
\cite{hall2011theory,hall2011asymptotic,bickel2013asymptotic,wang2019frequentist,gaucher2021optimality}
study properties of consistency and asymptotic normality for estimates of
low-dimensional parameters in latent variable models (i.e.\ of the prior
``hyperparameters'' in Bayesian contexts), using variational approximations
for the marginal log-likelihood. In particular, 
\cite{bickel2013asymptotic,gaucher2021optimality} establish such guarantees
for the mean-field variational approximation in stochastic block models (SBMs),
which are closely related to the $\Z_2$-synchronization model of our work.
\cite{mukherjee2018mean,plummer2020dynamics,zhang2020theoretical} study 
the optimization landscape and convergence properties of iterative coordinate
ascent (CAVI) and block coordinate ascent (BCAVI) algorithms, with
\cite{zhang2020theoretical} showing that BCAVI achieves an optimal
exponentially-vanishing rate of estimation error for the latent community
membership vector in SBMs with asymptotically diverging signal strength.
\cite{zhang2020convergence,alquier2020concentration,cherief2019consistency,yang2020alpha,ray2021variational}
study rates of posterior contraction for both variational Bayes and
$\alpha$-fractional variational Bayes methods, establishing conditions under
which the variational posteriors may enjoy the same optimal rates of contraction
in a frequentist Bernstein-von-Mises sense as the true Bayes posteriors.
In particular, \cite{alquier2020concentration,cherief2019consistency,yang2020alpha} discuss
applications of these results to low-rank matrix estimation problems,
including matrix completion, probabilistic PCA, and topic models.

In our work, we study the $\Z_2$-synchronization model with bounded signal
strength, which is in a different asymptotic regime from the above posterior
contraction results for SBMs and low-rank matrix estimation.
Fixing the true parameter $\bx$ as the all-1's vector, the Bayes estimate
for $\bx$ in our setting has a marginal
distribution of coordinates that converges to a non-degenerate limit law, and
an asymptotically non-vanishing per-coordinate Bayes risk.

Our focus on such a setting is motivated in part by our belief that in
many applications, Bayesian approaches to inference may be favored because the
data is in a regime of limited signal-to-noise that is
far from theoretical regimes of posterior contraction. Instead, information
in the hypothesized prior is important in informing inference, and the
desideratum is then to obtain an accurate estimate of the posterior distribution under this prior. Our results
are oriented towards this goal, showing (in a simple but illustrative model)
that minimizing the TAP free energy yields a variational approximation which
consistently estimates the posterior marginals, even when the posterior
distribution itself does not concentrate strongly around the true parameter.

\subsubsection{TAP free energy and the naive mean-field approximation}

Thouless, Anderson, and Palmer introduced in \cite{thouless1977solution} the TAP
equations (and the associated TAP free energy) as a system of asymptotically
exact mean-field equations in the SK model. For spin
glasses, the validity of the TAP equations and their relation to the Gibbs
measure have been extensively studied---see for example
\cite{plefka1982convergence,de1983weighted,bray1984weighted,cavagna2003formal}
in the physics literature, and \cite{talagrand2010mean, chatterjee2010spin,
bolthausen2014iterative, auffinger2019thouless, chen2018tap,
chen2018generalized, belius2019tap, subag2021free} for rigorous
mathematical results. Direct optimization of an analogous TAP
free energy (a.k.a.\ approximate Bethe free energy) was proposed for
Bayesian linear and generalized linear models in
\cite{krzakala2014variational,rangan2016inference}, which recognized
that its critical points are in exact correspondence with fixed points of AMP.
$\Z_2$-synchronization corresponds to the SK model with an added ferromagnetic
bias, and the form of the TAP free energy that we study is identical to the
(high-temperature) TAP free energy in the SK model with this added
ferromagnetic component.

We emphasize that both the TAP approach and the ``naive'' mean-field
approach of (\ref{eq:naivemeanfield}) have received significant attention in
the theoretical literature. A line of work
\cite{chatterjee2016nonlinear,basak2017universality,eldan2018gaussian,jain2018mean,yan2020nonlinear,augeri2020nonlinear} on the theory
of non-linear large deviations establishes that the naive mean-field
approximation to the free energy (i.e.\ the marginal log-likelihood in Bayesian
models) is asymptotically accurate to leading order, without the need for a
TAP correction, under a condition that the log-density has a ``low-complexity
gradient''. In Ising models with couplings matrix $\bY \in \R^{n \times n}$
having $O(1)$ operator
norm, such a condition holds when $\bY$ is nearly low-rank in the sense
$\|\bY\|_F^2=o(n)$ \cite{basak2017universality}. It does not hold
for $\Z_2$-synchronization with any fixed signal strength $\lambda$, where
\cite{ghorbani2019instability,fan2021tap} contrasted variational inference
based on the TAP and naive mean-field approximations. In particular,
\cite{ghorbani2019instability} showed that for $\lambda \in (1/2,1)$, naive
mean-field variational Bayes may yield a ``falsely informative'' variational
posterior, and \cite{fan2021tap} showed that critical points of the naive
mean-field free energy cannot correspond to consistent approximations of the
posterior mean for any sufficiently large but fixed value of $\lambda$.

\subsubsection{Spiked matrix models and $\Z_2$-synchronization}

Spiked matrix models have been a mainstay in the statistical literature since
their introduction by \cite{johnstone2001distribution}. $\Z_2$-synchronization
is a specific example of the spiked model with Bernoulli prior, and also of
more general synchronization problems over compact groups
\cite{singer2011angular,bandeira2020non}. The Bayes risks in
$\Z_2$-synchronization and other spiked matrix models were studied in
\cite{deshpande2016asymptotic, barbier2016mutual, krzakala2016mutual,
lelarge2019fundamental}. For $\Z_2$-synchronization, non-trivial signal 
estimation above the weak-recovery threshold $\lambda=1$ can also be achieved
by spectral methods \cite{baik2005phase, peche2006largest}
and semi-definite programming
\cite{montanari2016semidefinite,javanmard2016phase}, although such methods do
not achieve the asymptotically optimal Bayes risk (\ref{eq:Bayesrisk}).

$\Z_2$-synchronization has been studied in part as a simpler
analogue of the symmetric two-component SBM that replaces
the noise $\bfA-\E[\bfA]$ of the adjacency matrix $\bfA$ by Gaussian noise,
and it is possible to make formal connections between estimation in these models
via universality arguments
\cite{deshpande2016asymptotic,montanari2016semidefinite}. We believe that
certain aspects of our analyses and results may also be extendable to the SBM
via universality arguments developed for AMP 
in \cite{bayati2015universality, chen2021universality,
wang2022universality, dudeja2022spectral} and for minimizers of optimization
objective functions with random data in
\cite{montanari2017universality,hu2020universality,montanari2022universality,han2022universality},
and this would be interesting to explore in future work.

\subsubsection{AMP algorithms}

AMP algorithms were proposed and studied in
\cite{kabashima2003cdma,donoho2009message} for Bayesian linear regression and
compressed sensing. They may be derived by approximating belief propagation on
dense graphical models,
see e.g.\ \cite{donoho2010message, montanari2012graphical}.
Various generalizations of AMP have been developed, including the Generalized AMP
algorithm of \cite{rangan2011generalized} and the Vector AMP algorithm of
\cite{rangan2019vector}, and we refer to \cite{feng2022unifying} for a recent
review. The state evolution formalism of AMP was introduced in
\cite{donoho2009message} and rigorously established in
\cite{bolthausen2014iterative, bayati2011dynamics}. This has since been
generalized in \cite{javanmard2013state, berthier2020state,
montanari2021estimation}. A finite-$n$ analysis of AMP was performed in
\cite{rush2018finite}, which extended the validity of the state evolution to
$o(\log n / \log \log n)$ iterations. Following the initial posting
of our work, \cite{weiLi2022} established a different finite-$n$ guarantee for
AMP via a novel decomposition of the AMP iterates, which applies for
$o(n/(\log^7 n))$ iterations in the $\Z_2$-synchronization problem
with signal strength $\lambda \in (1,1.2)$.

\subsubsection{Gaussian comparison inequalities}

The proofs of our main results rely heavily on Slepian's comparison inequality
\cite{slepian1962one} and its later development by Sudakov-Fernique
\cite{sudakov1971gaussian,sudakov1979geometric, fernique1975regularite}, to
reduce the study of $\cF_\TAP$ to a simpler Gaussian process. This approach is
related to a recent line of work that generalizes Gordon's inequality 
\cite{gordon1985some, kahane1986inegalite} to a Convex Gaussian Minimax
Theorem (CGMT) \cite{stojnic2013framework, oymak2013squared, thrampoulidis2015regularized, miolane2021distribution, celentano2020lasso}.

\subsubsection{Kac-Rice formula and complexity analysis}

Physics calculations of the complexity of critical points in spin glass models
using the Kac-Rice formalism can be found in \cite{bray1980metastable,
cavagna2003formal, crisanti2003complexity, fyodorov2004complexity,
crisanti2005complexity}. This method was made rigorous for spherical spin
glasses in \cite{auffinger2013random, auffinger2013complexity,subag2017complexity},
and a more recent line of work \cite{arous2019landscape,
maillard2020landscape, fan2021tap, baskerville2021loss, baskerville2022spin,
arous2021landscape} has used this approach to analyze non-convex
function landscapes in other high-dimensional probabilistic and statistical models.

\section{Main results}

\subsection{Local analysis of the TAP free energy}

Our first result shows the existence and uniqueness of a local minimizer of the
TAP free energy $\cF_\TAP$ near the Bayes estimator (c.f.\ Eq.\ (\ref{eq:XBayes})), for any signal
strength $\lambda>1$. We also establish strong convexity of $\cF_\TAP$ in a
$\sqrt{\eps n}$-neighborhood around this minimizer, as well as the stability
of the AMP map
\begin{equation}\label{eq:TAMP}
T_{\AMP}(\bbm,\bbm_-)=\Big(\tanh\Big(\lambda \bY \bbm-\lambda^2[1-Q(\bbm)]
\bbm_-\Big),\;\bbm\Big)
\end{equation}
at this local minimizer. This is the map for which the AMP
iterations (\ref{eq:AMP}) may be expressed as
$(\bbm^{k+1},\bbm^k)=T_{\AMP}(\bbm^k,\bbm^{k-1})$.

\begin{theorem}[Local convexity and AMP stability]\label{thm:local}
Fix any $\lambda>1$. There exist $\lambda$-dependent constants $\eps, t > 0$ and
$r \in (0, 1)$ such that for any fixed $\iota>0$, with probability
approaching 1 as $n \to \infty$, the following all occur.
\begin{enumerate}
\item[(a)] (Bayes-optimal TAP local minimizer)
Let $\widehat{\bX}_{\Bayes}=\E[\bx\bx^\top \mid \bY]$. There exists a critical point and local minimizer $\bbm_\star$ of $\cF_\TAP(\bbm)$ such that
\begin{equation}\label{eq:bayesoptimal}
\frac{1}{n^2}\|\bbm_\star\bbm_\star^\sT-\widehat{\bX}_{\Bayes}\|_\sF^2<\iota.
\end{equation}
For sufficiently small $\iota>0$ (which is $\lambda$-dependent and $n$-independent), this is the unique critical point
satisfying (\ref{eq:bayesoptimal}) up to $\pm$ sign.
\item[(b)] (Local strong convexity of TAP free energy)
Let $\lambda_{\min}(\cdot)$ denote the smallest eigenvalue. For this local minimizer $\bbm_\star$, we have
\[\lambda_{\min}\Big(n \cdot \nabla^2 \cF_{\TAP}(\bbm)\Big)>t>0
\text{ for all } \bbm \in (-1,1)^n \cap \ball_{\sqrt{\eps n}}(\bbm_\star).\]
In particular, $\cF_{\TAP}$ is strongly convex over
$(-1,1)^n \cap \ball_{\sqrt{\eps n}}(\bbm_\star)$.
\item[(c)] (Local stability of AMP) Let $\de T_{\AMP} \in \reals^{2n \times 2n}$
be the Jacobian of the AMP map (\ref{eq:TAMP}), and let $\rho(\cdot)$ denote 
the spectral radius. For this local minimizer $\bbm_\star$, we have
\[\rho\Big(\de T_{\AMP}(\bbm_\star,\bbm_\star)\Big)<r<1.\]
\end{enumerate}
\end{theorem}

Combining with the global landscape analysis of \cite{fan2021tap}, this implies
the following immediate corollary for large enough signal strength $\lambda$.

\begin{corollary}[Global landscape for large
$\lambda$]\label{cor:globallandscape}
For an absolute constant $\lambda_0>0$, suppose $\lambda>\lambda_0$. Then with
probability approaching 1 as $n \to \infty$,
the local minimizers $\pm \bbm_\star$ guaranteed by Theorem
\ref{thm:local} are the global minimizers of $\cF_\TAP$. Furthermore, they are
the only critical points of $\cF_\TAP$ in the domain
\[\cS=\Big\{\bbm \in (-1,1)^n:\cF_\TAP(\bbm)<-\lambda^2/3\Big\}.\]
\end{corollary}

A proof sketch of Theorem \ref{thm:local} can be found in Section \ref{sec:local}, and its detailed proof can be found in Appendix \ref{appendix:local}. The proof of Corollary \ref{cor:globallandscape} can be found in Appendix \ref{appendix:TAP_global_landscape}.

\subsection{Convergence of algorithms}

We study convergence of the AMP algorithm (\ref{eq:AMP}),
with the spectral initialization
\begin{equation}\label{eq:initialization}\tag{{\sf SI}}
\bh^0=\text{principal eigenvector of } \bY \text{ with }
\|\bh^0\|_2=\sqrt{n\lambda^2(\lambda^2-1)}, ~~
\bbm^{-1}=\lambda \bh^0.
\end{equation}
We choose this scaling for $\bh^0$ as in
\cite[Section 2.4]{montanari2021estimation} to simplify the AMP state evolution.

We introduce also the following more ``generic'' first-order 
natural gradient descent (NGD) algorithm, with a step size parameter $\eta>0$:
\begin{align}
\bbm^k&=\tanh(\bh^k)\nonumber\\
\bh^{k+1}&=\bh^k-\eta n \cdot \nabla \cF_\TAP(\bbm^k)\nonumber\\
&=(1-\eta)\bh^k+\eta\Big(\lambda \bY
\bbm^k-\lambda^2[1-Q(\bbm^k)]\bbm^k\Big).\label{eq:NGD}\tag{{\sf NGD}}
\end{align}
We call this algorithm ``natural gradient descent'' because we may apply
$(\de/\de h) \tanh(h)=1-\tanh(h)^2$ to write the $\bbm$-gradient
$\nabla \cF_\TAP(\bbm^k)$ equivalently as a preconditioned $\bh$-gradient,
\[\nabla \cF_\TAP(\bbm^k)
=\bI(\bbm^k)^{-1} \cdot \nabla_\bh \cF_\TAP(\tanh(\bh^k)),
\qquad \bI(\bbm)=\diag\left(\frac{1}{1-\bbm^2}\right),\]
where $\bI(\bbm)$ is proportional to the Fisher information matrix in a model of
$n$ independent Bernoulli $\{-1,+1\}$ variables with mean $\bbm \in \R^n$.
This identifies (\ref{eq:NGD}) as a natural gradient method
\cite{amari1998natural}. We note that setting the step size $\eta=1$
yields an algorithm similar to (\ref{eq:AMP}), but with $\bbm^{k-1}$
replaced by $\bbm^k$. For simplicity, we will consider the same spectral
initialization $\bh^0$ for this algorithm as for AMP
in (\ref{eq:initialization}),
although here this specific choice of initialization is less important.

Alternatively, the iterations (\ref{eq:NGD}) may be understood as a
mirror-descent/Bregman-gradient method in the $\bbm$-parameterization
\cite{nemirovskij1983problem,beck2003mirror}. Recalling the 
binary entropy function $\sh$ from (\ref{eq:entropy}), we define
\begin{equation}\label{eq:Bregman}
\begin{aligned}
L=\frac{1}{\eta}, ~~& H(\bbm)=\frac{1}{n}\sum_{i=1}^n \sh(m_i),\\
D_{-H}(\bbm,\bbm')&=-H(\bbm)+H(\bbm')+\langle \nabla H(\bbm'),
\bbm-\bbm' \rangle
\end{aligned}
\end{equation}
where $L$ is the inverse step size, $-H(\bbm)$ is a separable convex prox
function, and $D_{-H}(\bbm,\bbm')$ is its associated Bregman divergence. Then
it may be checked that (\ref{eq:NGD}) takes the equivalent mirror-descent form
\begin{equation}\label{eq:mirrordescent}
\bbm^{k+1}=\argmin_{\bbm \in (-1,1)^n} \cF_\TAP(\bbm^k)
+\langle \nabla \cF_\TAP(\bbm^k), \bbm-\bbm^k \rangle
+L \cdot D_{-H}(\bbm,\bbm^k).
\end{equation}
One motivation for studying this algorithm, rather than ordinary gradient
descent in the $\bbm$-parameterization, is that the Hessian
$\nabla^2 \cF_\TAP(\bbm)$ is not uniformly bounded over $(-1,1)^n$, and instead
diverges as $\bbm$ approaches
the boundaries of the cube. The form (\ref{eq:mirrordescent}) naturally adapts
to this non-uniform curvature of $\cF_\TAP$, allowing for a convergence
analysis using techniques of \cite{bauschke2017descent,lu2018relatively} for
minimizing functions that are not strongly smooth in the Euclidean metric.

Combining the local strong convexity of Theorem \ref{thm:local}, the state
evolution of spectrally-initialized AMP, and this type of convergence analysis
for NGD, we deduce the following result, whose proof can be found in Section \ref{sec:proof_local_convergence_NGD} and Appendix \ref{appendix:global}.

\begin{theorem}[Computation of Bayes-optimal TAP minimizer]\label{thm:localconvergence}
Fix any $\lambda>1$. There exist $\lambda$-dependent constants
$C,\mu,\eta_0>0$ and $T \geq 1$
such that with probability approaching 1 as $n \to \infty$, the following occurs.

Fix any step size $\eta \in (0,\eta_0)$,
let $\bbm^{T} \in (-1,1)^n$ be the $T^\text{th}$ iteration of (\ref{eq:AMP})
from the spectral initialization (\ref{eq:initialization}), and let $\bbm^{T+k}
\in (-1, 1)^n$ be obtained by $k$ iterations of (\ref{eq:NGD}) with
step size $\eta$ from the initialization $\bbm^T$. Let $\bbm_\star$ be the Bayes-optimal local minimizer of
$\cF_\TAP$ in Theorem \ref{thm:local}. Then for some choice of sign $\pm$ and
every $k \geq 1$, 
\begin{align*}
\cF_\TAP(\bbm^{T+k})-\cF_\TAP(\pm \bbm_\star)&<C(1-\mu\eta)^k,\\
\|\bbm^{T+k}-(\pm \bbm_\star)\|_2&<C(1-\mu\eta)^k\sqrt{n}.
\end{align*}
In particular, $\lim_{k \to \infty} \bbm^{T+k} \in \{+\bbm_\star,-\bbm_\star\}$.
\end{theorem}

This theorem implies that for any fixed value of $\lambda>1$, the
Bayes-optimal local minimizer $\bbm_\star$ of $\cF_\TAP$ guaranteed by
Theorem \ref{thm:local} may be computed in time that is polynomial
in the problem size $n$ (in the usual sense of linear convergence).
Let us remark that the convergence analysis of NGD in this result
is purely geometric, relying only on the smoothness and local convexity
properties of $\cF_\TAP$. We hence expect that a similar convergence analysis
may be performed for momentum-accelerated or stochastic variants of NGD,
such as those developed recently
in \cite{hanzely2021accelerated,gutman2022perturbed,dragomir2021fast}.

For sufficiently large signal strength $\lambda$, where the more
global landscape of $\cF_\TAP$ is clarified by Corollary
\ref{cor:globallandscape}, our next result Theorem \ref{thm:globalconvergence} verifies that the hybrid AMP/NGD
approach in Theorem \ref{thm:localconvergence} is not needed, and that either
algorithm alone can achieve linear convergence to the global TAP minimizer
$\bbm_\star$ from a spectral initialization. The proof of Theorem \ref{thm:globalconvergence} can be found in Section \ref{sec:proof_NGD_global} and \ref{sec:proof_AMP_global}, and Appendix \ref{appendix:global}.

\begin{theorem}[Convergence of AMP and NGD for large $\lambda$]
\label{thm:globalconvergence}
For an absolute constant $\lambda_0>0$, suppose $\lambda>\lambda_0$ and let $\bbm_\star$ be the global minimizer of $\cF_\TAP$ in Corollary \ref{cor:globallandscape}. 
Then there exist $\lambda$-dependent constants $C,\mu,\eta_0>0$ and $\alpha \in (0,1)$
such that with probability approaching 1 as $n \to \infty$, the following all
occur.
\begin{enumerate}
\item[(a)] (Convergence of AMP)
Let $\bbm^k$ be the $k^\text{th}$ iterate of \ref{eq:AMP} from the spectral
initialization (\ref{eq:initialization}).
For some choice of sign $\pm$ and every $k \geq 1$,
\[\cF_\TAP(\bbm^k)-\cF_\TAP(\pm \bbm_\star)<C\alpha^k,
\qquad \|\bbm^k-(\pm \bbm_\star)\|_2<C\alpha^k\sqrt{n}.\]
\item[(b)] (Convergence of NGD)
Fix any step size $\eta \in (0,\eta_0)$, and
let $\bbm^k$ be the $k^\text{th}$ iterate of \ref{eq:NGD} from the spectral
initialization (\ref{eq:initialization}) with step size $\eta$.
For some choice of sign $\pm$ and every $k \geq 1$,
\[\cF_\TAP(\bbm^k)-\cF_\TAP(\pm \bbm_\star)<C(1-\mu\eta)^k,
\qquad \|\bbm^k-(\pm \bbm_\star)\|_2<C(1-\mu\eta)^k\sqrt{n}.\]
\end{enumerate}
In particular, for both algorithms, $\lim_{k \to \infty} \bbm^k \in
\{+\bbm_\star,-\bbm_\star\}$.
\end{theorem}

\begin{remark}
We believe that the requirement $\lambda>\lambda_0$ sufficiently large in
Theorem \ref{thm:globalconvergence} is artificial, and that this result also holds for
all $\lambda>1$. This is supported by numerical simulations in
Section \ref{sec:simulations} below. Let us clarify that such a guarantee for AMP
does not follow from its state evolution combined with its local stability
shown in Theorem \ref{thm:local}(c): The state evolution ensures convergence
to a $\sqrt{\eps n}$-neighborhood of $\bbm_\star$, for any fixed $\eps>0$, in a
finite number of AMP iterations. However, the local stability in
Theorem \ref{thm:local}(c) does not quantify the size of the neighborhood of 
$\bbm_\star$ in which AMP is then guaranteed to converge to $\bbm_\star$.
\end{remark}

\begin{remark}
Part of our analysis of Theorem \ref{thm:globalconvergence}(a) still
uses the state evolution for AMP with spectral initialization developed
in \cite{montanari2021estimation}. This result would hold equally if AMP is
initialized with a vector $\bbm_1$ that is independent of the noise matrix $\bW$
and has non-vanishing correlation with $\bbm_\star$, by the validity
of the AMP state evolution also in this setting. For a random initialization
that is uncorrelated with $\bbm_\star$, we note that an analysis of AMP seems
challenging even in this setting of large but fixed $\lambda>\lambda_0$, as the
algorithm would still require $O(\log(n))$ iterations to achieve a
non-negligible correlation with $\bbm_\star$, and existing finite-$n$ analyses
of AMP \cite{rush2016,weiLi2022} do not seem to immediately apply to describe
this early phase of optimization. In Theorem \ref{thm:globalconvergence}(b),
the spectral initialization is used to ensure that NGD is initialized in a
basin of attraction of $\bbm_\star$, and analyses of the global
landscape of $\cF_\TAP$ in \cite{fan2021tap} are also insufficient to show that 
this basin of attraction includes random initializations.
\end{remark}

\section{Numerical simulations}\label{sec:simulations}

\subsection{Convergence of algorithms}

\begin{figure}[t]
\centering
\includegraphics[width = 0.4\linewidth]{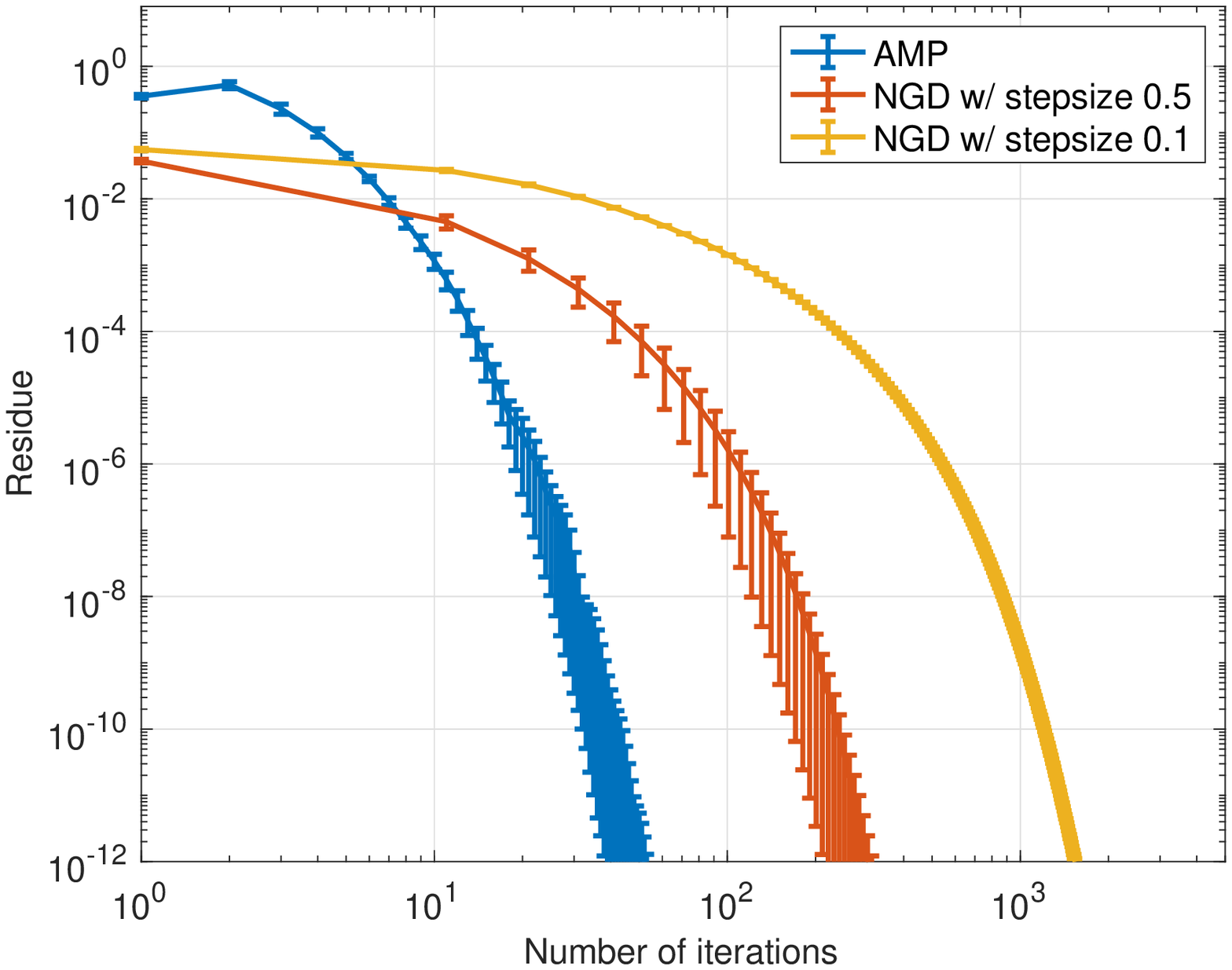}~~~~
\includegraphics[width = 0.4\linewidth]{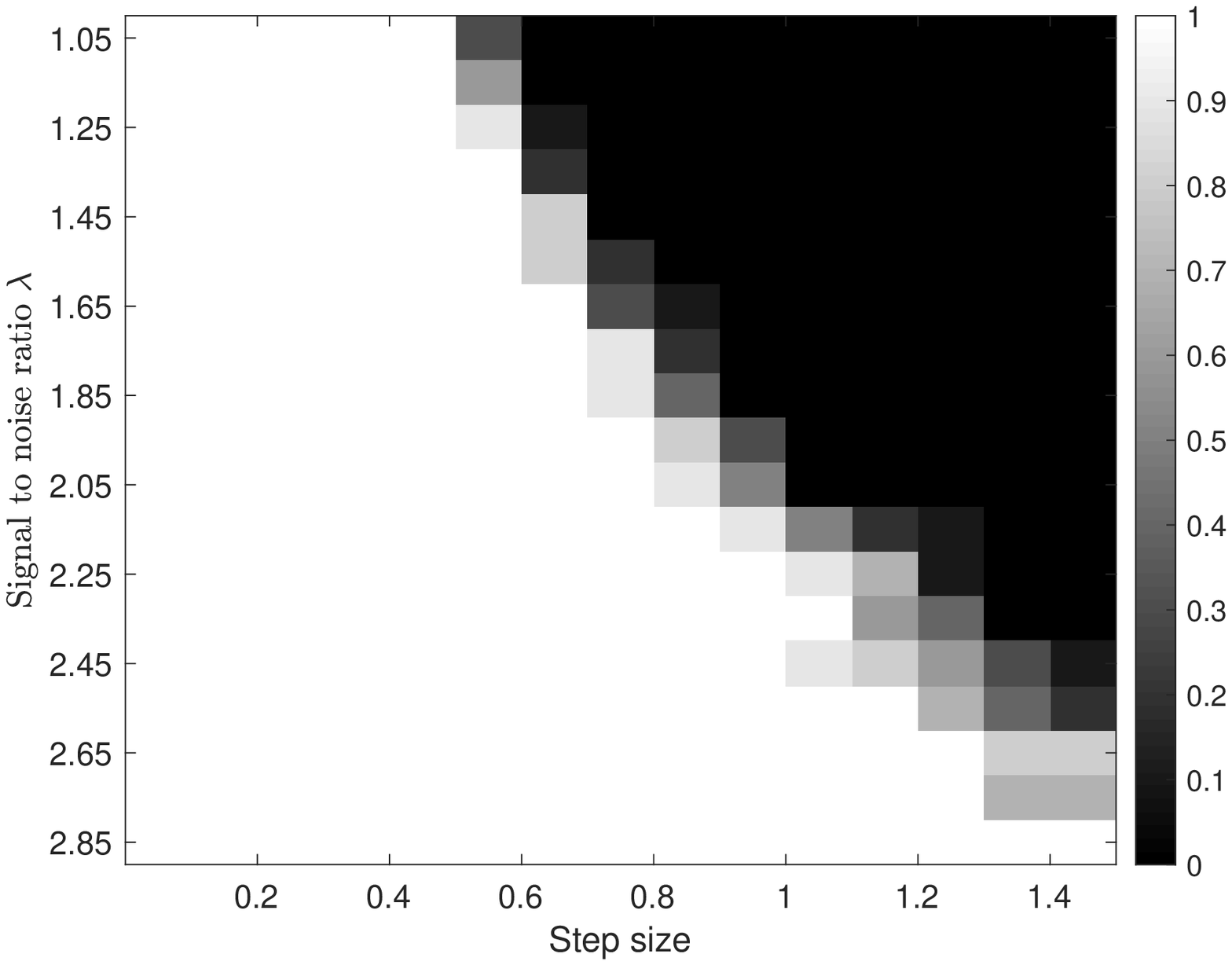}
\caption{Convergence of AMP and NGD from a spectral initialization. Left: Residual squared
error $\min\{ \| \bbm^{k} - \bbm_\star \|_2^2 / n, \| \bbm^{k} + \bbm_\star
\|_2^2 / n\}$ versus number of iterations $k$ (both on a log-scale), for
signal-to-noise ratio $\lambda=1.5$. The mean curve is averaged over
$10$ independent instances, and the error bars report $1/\sqrt{10}$ times the
standard deviation across instances. Right: Success probability of NGD for
convergence to $\bbm_\star$, for varying signal-to-noise ratios $\lambda$ and
step sizes $\eta$. In both panels, $n = 500$.
}\label{fig:AMP_NGD}
\end{figure}

We perform numerical simulations to confirm the global convergence of AMP and
NGD for all $\lambda>1$, and to compare their convergence rates. 
We initialize both AMP and NGD using
the spectral initialization (\ref{eq:initialization}).

In Figure \ref{fig:AMP_NGD}(a), we plot the residual squared error $\min\{ \|
\bbm^{k} - \bbm_\star \|_2^2 / n, \| \bbm^{k} + \bbm_\star \|_2^2 / n\}$, where
$\bbm^k$ is the $k^\text{th}$ iterate of AMP or NGD with different step sizes,
and $\bbm_\star = \argmin_{\bbm} \cF_\TAP(\bbm)$. (We first compute $\bbm_\star$
up to high numerical accuracy using AMP.)
For each algorithm, we simulated 10 random instances of
$\bY \in \R^{n \times n}$ according to the $\Z_2$-synchronization model
(\ref{eq:Z2-sync}),
with $n=500$ and $\lambda=1.5$. Figure \ref{fig:AMP_NGD}(a) shows that AMP and
NGD with step sizes 0.1 and 0.5 all consistently achieve convergence to
$\bbm_\star$, where AMP has the fastest rate of convergence.

In Figure \ref{fig:AMP_NGD}(b), we report the success probability of NGD for
achieving convergence to $\bbm_\star$, for various step sizes $\eta$ (horizontal
axis) and signal-to-noise ratios $\lambda>1$ (vertical axis). The success
probability is defined as the fraction of the 10 random instances of $\bY$
for which NGD
achieved residual squared error $10^{-4}$ within $k=12000$ iterations.
Figure \ref{fig:AMP_NGD}(b) suggests that NGD with step size $\eta<0.4$
converges for any $\lambda>1$, and illustrates that as $\lambda$ increases, NGD
allows for a larger step size in achieving this convergence.

\subsection{Universality with respect to the noise distribution}

\begin{figure}[t]
\centering
\includegraphics[width = 0.4\linewidth]{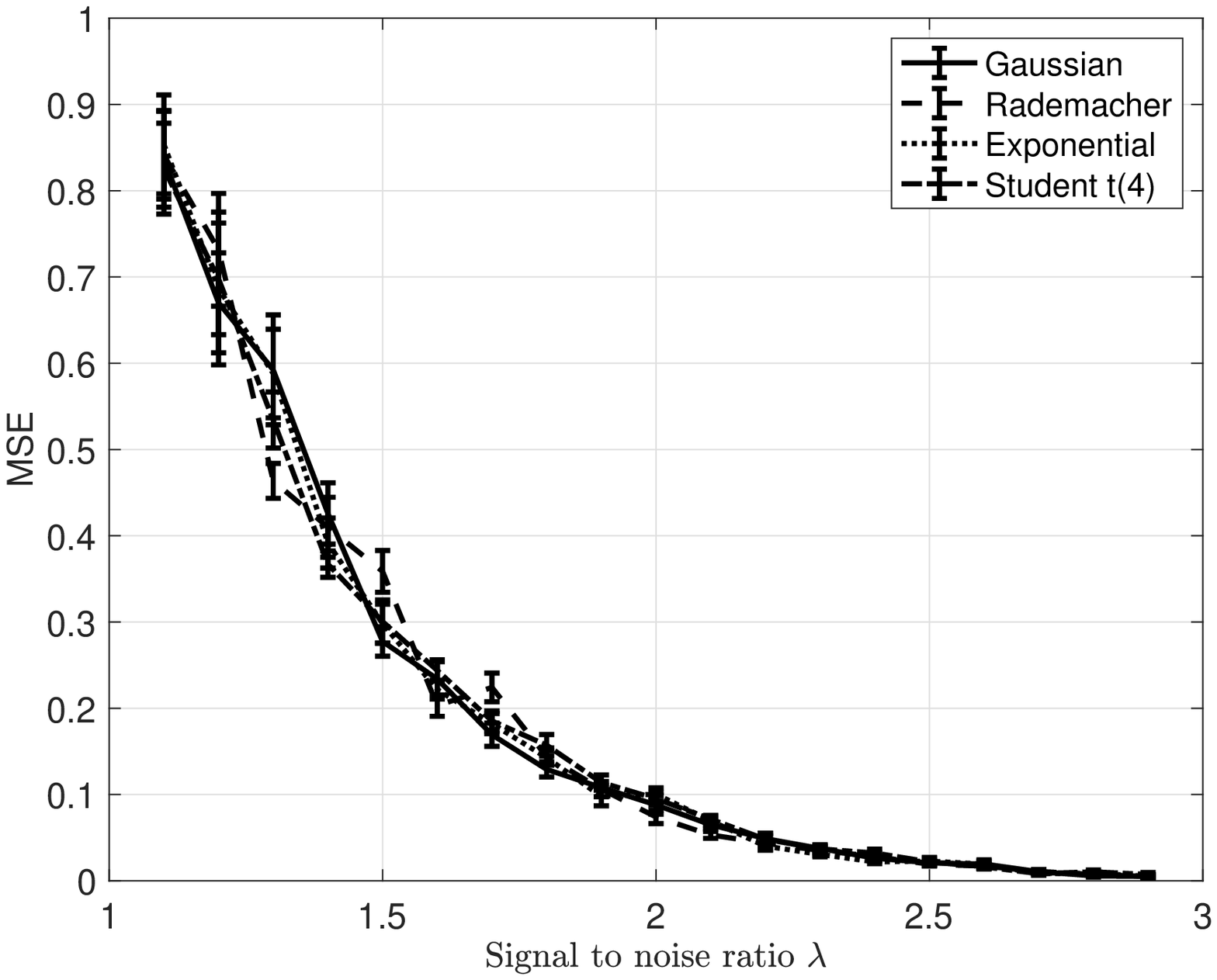}~~~~
\includegraphics[width = 0.4\linewidth]{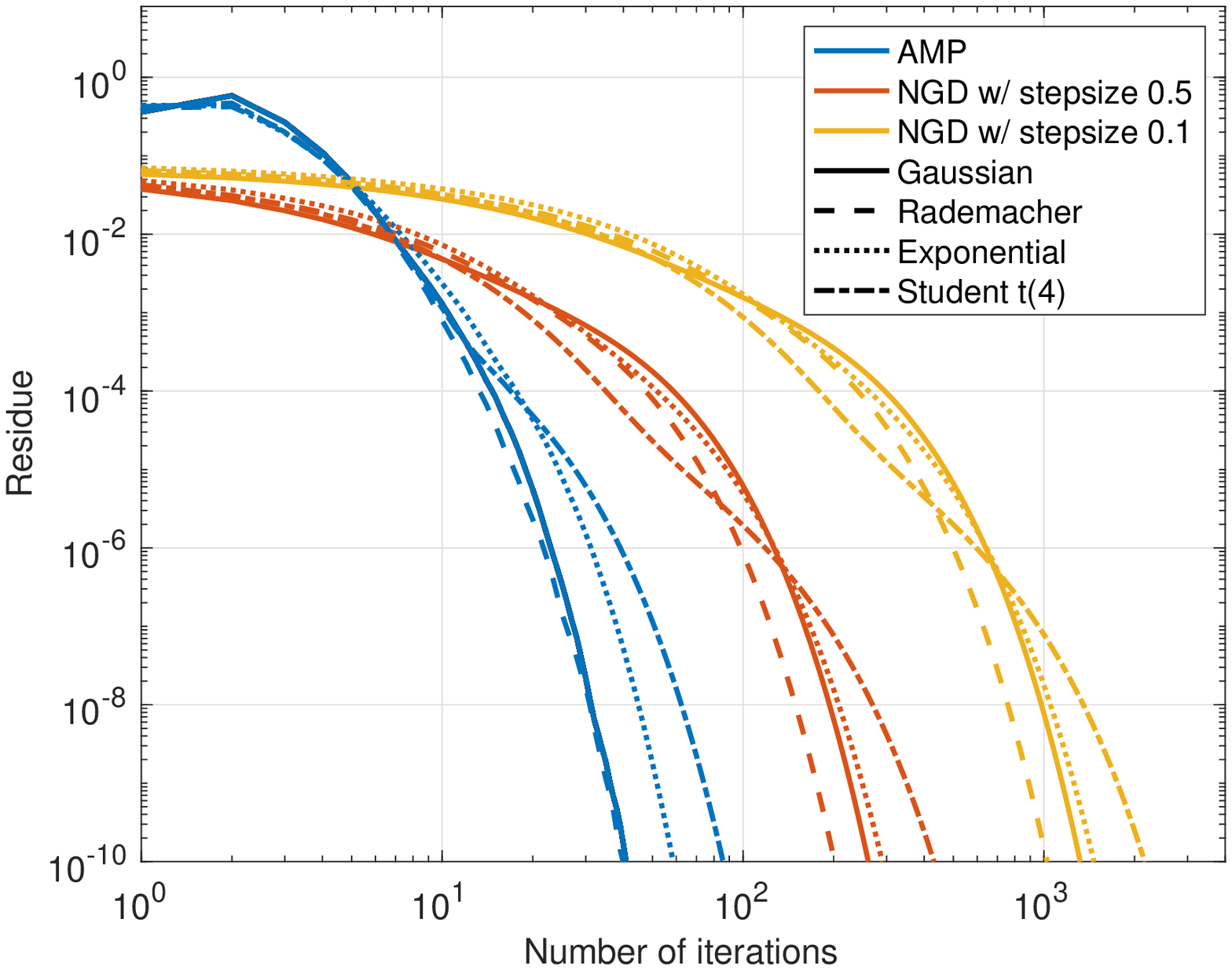}
\caption{Universality with respect to the noise distribution. Left: Estimation mean squared error $\min\{\| \bbm_\star - \bx \|_2^2 /
n, \| \bbm_\star - \bx \|_2^2 / n \}$ versus the signal-to-noise ratio
$\lambda$, for different noise ensembles. The mean curve is averaged
over 10 independent instances, and the error bars report $1/\sqrt{10}$ times the
standard deviation across instances. Right: Residual squared error $\min\{ \|
\bbm^{k} - \bbm_\star \|_2^2 / n, \| \bbm^{k} + \bbm_\star \|_2^2 / n\}$ versus
the number of iterations $k$, for different noise ensembles and
signal-to-noise ratio $\lambda=1.5$. In both panels, $n=500$. }\label{fig:universality}
\end{figure}

Although we analyze AMP and NGD for Gaussian noise, we expect the 
properties of these estimators and of the TAP free energy landscape to be
robust under sufficiently light-tailed distributions of noise entries. Here, we verify this numerically for  three examples of symmetric
non-Gaussian noise matrices $\bW$:
\begin{itemize}
\item Rademacher: $(w_{ij}:1 \leq i \leq j \leq n) \overset{iid}{\sim}
\Unif\{-1/\sqrt{n}, 1/\sqrt{n}\}$.
\item Double-exponential (Laplace): $\bW=(\bG+\bG^\sT)/\sqrt{2n}$, where $(G_{ij}: 1 \le i, j \le n) \overset{iid}{\sim} (1 / \sqrt{2}) \exp\{ -  \sqrt{2} \cdot \vert x \vert \}$. 
\item Student's t: $\bW= (\bG+\bG^\sT)/\sqrt{2n}$, where $(G_{ij}: 1
\le i, j \le n) \overset{iid}{\sim} t(\nu) / \sqrt{\nu/(\nu-2)}$ and the
degrees-of-freedom is $\nu = 4$.
\end{itemize}
In all three examples, all entries $w_{ij}$ have mean 0, and
all off-diagonal entries $w_{ij}$ have variance $1/n$.

In Figure \ref{fig:universality}(a), we report the estimation mean squared error
(MSE) $\min\{\| \bbm_\star - \bx \|_2^2 / n, \| \bbm_\star + \bx \|_2^2 / n\}$
versus $\lambda$, where $\bbm_\star = \argmin_{\bbm} \cF_{\TAP}(\bbm)$
is computed from AMP up to high accuracy as before, and the
noise matrix $\bW$ is generated from either the assumed Gaussian (GOE) model or 
from the above three non-Gaussian ensembles. In Figure \ref{fig:universality}(b), we
report the residual squared error $\min\{ \| \bbm^{k} - \bbm_\star \|_2^2 / n,
\| \bbm^{k} + \bbm_\star \|_2^2 / n\}$ versus the number of algorithm iterations
$k$, for the same four noise ensembles. These figures show that properties of the
TAP minimizers and of the AMP and NGD iterates are indeed robust to
these distributions
of the noise entries, even for some heavy-tailed distributions. 

We also tested Student's t-distribution with
degrees-of-freedom $\nu=3$, and observed that when $\lambda \in (1,2)$ and
$n=500$, AMP oscillates between two points rather than converging to a fixed
point. Instead, the NGD algorithm with a sufficiently small step size continues to converge to the global minimizer.

\subsection{Comparing TAP and mean-field variational Bayes}

\begin{figure}[t]
\centering
\includegraphics[width = 0.42\linewidth]{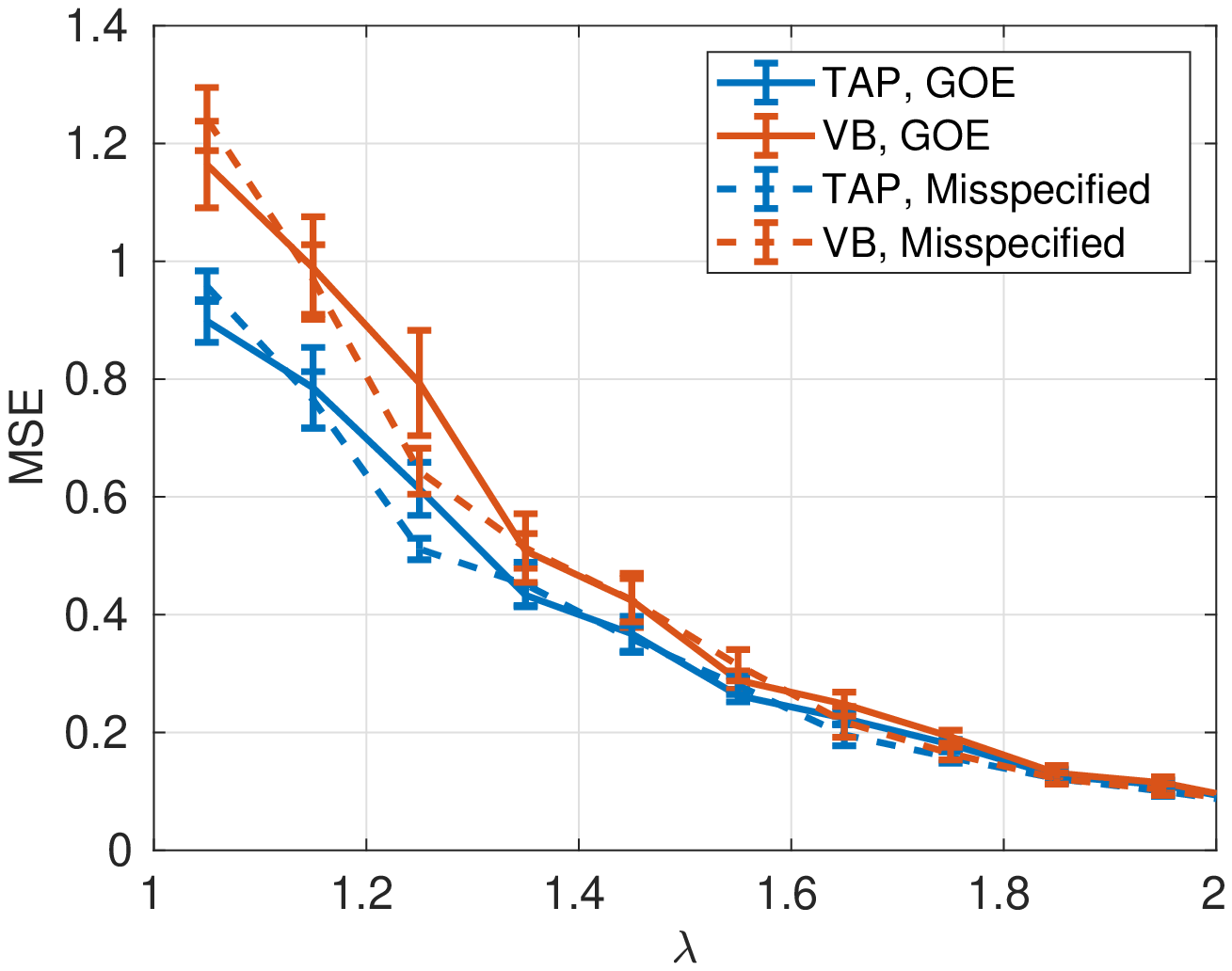}
\caption{Comparison of TAP with mean-field VB. The plot shows mean squared errors 
of the TAP and VB minimizers in both a correctly specified and a misspecified
model, for signal-to-noise ratio $\lambda \in [1, 2]$ and $n=500$. The mean curve is
averaged over 10 independent instances, and the error bars report $1/\sqrt{10}$
times the standard deviation across instances.}\label{fig:compare_TAP_VB}
\end{figure}

We compare the TAP approach to naive mean-field variational Bayes 
(mean-field VB), under both a correctly specified noise model and a misspecified
model that lies outside of the preceding universality class. 

For $\Z_2$-synchronization, parametrizing (\ref{eq:naivemeanfield}) by the mean
vector $\bbm=\E_{\bx \sim q}[\bx]$ gives the mean-field VB free energy
\[
\cF_{\VB}(\bbm)=- \frac{1}{n} \sum_{i = 1}^n \sh(m_i) - \frac{\lambda}{2n} \< \bbm, \bY \bbm\>. 
\]
This coincides with (\ref{eq:TAP}) upon removing the TAP correction term.

In Figure \ref{fig:compare_TAP_VB}, we compare the mean squared errors
$\min \{ \| \bbm_\star - \bx \|_2^2 / n, \| \bbm_\star + \bx \|_2^2 / n\}$ for
the minimizers $\bbm_\star$ of $\cF_\TAP$ and of $\cF_\VB$, when $\bY$ is
generated according to the following two models: 
\begin{itemize}
\item The correctly specified $\Z_2$-synchronization model (\ref{eq:Z2-sync}).
\item A misspecified model $\bY = (\lambda / n) \bx \bx^\sT + \bW$,
where $\bx \sim \Unif(\{-1,+1\}^n)$ has the assumed discrete uniform prior, but
$\bW=\bU\bD\bU^\sT$ does not have independent entries. We choose
$\bU \in \R^{n \times n}$ as a uniformly sampled orthogonal matrix, and $\bD = \diag(d_1,\ldots,d_n)$ where $(d_i:1 \leq i \leq
n) \overset{iid}{\sim} \Unif([-\sqrt{3}, \sqrt{3}])$. By this scaling of $\bW$,
we have $\| \bW \|_{\sF}^2 / n \approx 1$ which matches the scaling of $\bW \sim
\GOE(n)$.
\end{itemize}

For both free energies, we compute their (possibly local) minimizers using the NGD iterations $\bh^{k + 1} = \bh^{k} - \eta n \nabla_\bbm \cF(\bbm^k)$, with step
size $\eta=0.1$ and a spectral initialization. We observe that NGD
typically converged within $k=8000$ iterations (in the sense of achieving a
small gradient), despite the lack of a theoretical convergence guarantee 
in certain settings. Under this model
misspecification, the minimizer $\bbm_\star$ of $\cF_\TAP$ defined according to
(\ref{eq:TAP}) is no longer expected to be asymptotically exact
for the Bayes posterior mean in the true generating model. Nonetheless, we
observe that $\bbm_\star$ which minimizes $\cF_\TAP$ achieves lower
mean squared error than that which minimizes $\cF_\VB$, in both the
well-specified and misspecified examples. For larger values of $\lambda$, the
difference in mean squared error between these approaches becomes harder to
discern, although the theory implies (in the well-specified setting) that this
difference is asymptotically non-vanishing for any fixed $\lambda>1$.

\section{Local analysis of the TAP free energy}\label{sec:local}

In this section, we describe the main ideas and steps in the proof of
Theorem \ref{thm:local}.

We will prove that each statement of the theorem holds with probability
approaching 1 conditional on the signal vector $\bx \in \{-1,+1\}^n$. By
symmetry, this
conditional probability is the same for any given vector $\bx \in \{-1,+1\}^n$,
so we may assume without loss of generality
\[\bx=\ones=(1,1,\ldots,1).\]
Conditional on $\bx$, the only remaining randomness is in the noise matrix $\bW
\sim \GOE(n)$, and $\cF_{\TAP}(\bbm)$ is a Gaussian process
indexed by $\bbm \in (-1,1)^n$.

The proof combines information derived from the Kac-Rice formula for the
expected number of critical points of Gaussian processes, the Sudakov-Fernique
Gaussian comparison inequality, and the AMP state evolution. It is helpful to
summarize the type of information each of these tools will provide:

\begin{description}
	\item[\textbf{Kac-Rice formula.}] We use the Kac-Rice formula to upper bound the
expected number of critical points of the TAP free energy in certain regions of
the domain $(-1,1)^n$, or for which the TAP Hessian or AMP Jacobian violate the
stated properties of Theorem \ref{thm:local}. In particular, by establishing
upper bounds that are vanishing as $n \to \infty$, we prove the non-existence
of such critical points with high probability.\\

	\item[\textbf{Sudakov-Fernique inequality.}] We use the Sudakov-Fernique
inequality to lower bound the infima of Gaussian processes defined by
$\bW \sim \GOE(n)$ with the infima of Gaussian processes defined by a standard Gaussian vector $\bg \in \R^n$. We then
analyze the latter to obtain variational lower bounds for large $n$. There are three Gaussian processes to which we apply this technique:
	\begin{itemize}

		\item The TAP free energy itself, to obtain lower bounds on its
minimum value over regions of $(-1,1)^n$.

		\item A Gaussian process whose infimum gives the minimum
eigenvalue of the TAP Hessian over subsets of $(-1,1)^n$, to show
local strong convexity of the TAP free energy.

		\item A Gaussian process whose infimum is related to the spectral
radius of the AMP Jacobian, to show local stability of the AMP map.\\

	\end{itemize}

	\item[\textbf{AMP state evolution.}] We use the AMP state evolution to evaluate
the TAP free energy at the iterates of AMP, giving upper bounds for the
TAP free energy value near the Bayes estimator.

\end{description}

\noindent The information provided by each of these three tools is distinct,
and the proof of Theorem \ref{thm:local} combines the information we can extract from each.

We outline the four main steps of the proof in Section
\ref{sec:localproofoutline}. 
These steps
are discussed in Sections \ref{sec:TAP-lb} through \ref{sec:AMP-stab}, and the
technical arguments that execute each step are deferred to
Appendix \ref{appendix:local}.

\subsection{Proof outline}\label{sec:localproofoutline}

For small parameters $\delta,\eta>0$, we define two deterministic subsets
$\cB_\delta,\cD_\eta \subset (-1,1)^n$ based on the empirical
distribution of coordinates of $\bbm \in (-1,1)^n$. These subsets will
contain the desired TAP local minimizer $\bbm_\star$ with high probability
(conditional on $\bx=\ones$).

For $\lambda>1$, let $q_\star=q_\star(\lambda)$ be the unique solution in
$(0,1)$ (cf.\ Proposition \ref{prop:qstar}) to the fixed-point equation
\begin{align}
q_\star&=\E_{G \sim \cN(0,1)}\big[\tanh(\lambda^2 q_\star + \lambda\sqrt{q_\star}
G)^2\big].\label{eq:fixedpointq}
\end{align}
Define
\begin{align}
h_\star&=\E_{G \sim \cN(0,1)}[\log 2 \cosh(\lambda^2 q_\star+
\lambda \sqrt{q_\star}G)]-\lambda^2 q_\star,\label{eq:hstar}\\
e_\star&=-\frac{\lambda^2}{4}(1-2q_\star-q_\star^2)
-\E_{G \sim \cN(0,1)}[\log 2\cosh(\lambda^2 q_\star+\lambda\sqrt{q_\star}G)].
\label{eq:estar}
\end{align}
For any point $\bbm \in (-1,1)^n$, denote
\[Q(\bbm)=\frac{1}{n}\|\bbm\|_2^2, \qquad M(\bbm)=\frac{1}{n}\bbm^\sT \ones,
\qquad H(\bbm)=\frac{1}{n}\sum_{i=1}^n \sh(m_i)\]
where $\sh(\cdot)$ is the binary entropy function from (\ref{eq:entropy}).
We define the first subset $\cB_\delta$ as
\begin{align}\label{eqn:set_B_delta}\cB_\delta=\Big\{\bbm \in (-1,1)^n:\;|Q(\bbm)-q_\star|,
|M(\bbm)-q_\star|,|H(\bbm)-h_\star|<\delta\Big\}.\end{align}
Let
\begin{equation}\label{eq:mustar}
\mu_\star=\text{distribution of } \tanh(\lambda^2 q_\star+\lambda
\sqrt{q_\star}G) \text{ when } G \sim \cN(0,1),
\end{equation}
which will be the limiting empirical distribution of coordinates of
$\bbm_\star$. 
For $\bbm \in (-1,1)^n$, let $\hat{\mu}_{\bbm}$ be the empirical distribution of coordinates of $\bbm$, i.e., 
\begin{equation}
\hat{\mu}_{\bbm} = \frac{1}{n} \sum_{i = 1}^n \delta_{m_i}.
\end{equation}
Denote by $W(\mu,\mu')$ the Wasserstein-2
distance between $\arctanh \mu$ and $\arctanh \mu'$, where $\arctanh \mu$ is 
shorthand for the law of $\arctanh m$ when $m \sim \mu$. That is, we have 
\begin{align}
W(\mu, \mu') &\equiv W_2(\arctanh \mu, \arctanh \mu')\nonumber\\
&= \left(\inf_{\text{couplings } \nu \text{ of } (\mu,\mu')}
\int (\arctanh m- \arctanh m')^2 \de \nu(m,m')\right)^{1/2}. 
\label{eqn:W_distance}
\end{align}
We review properties of this distance in Appendix \ref{appendix:wass2}. We define the second subset
$\cD_\eta$ as
\begin{align}\label{eqn:set_D_eta}\cD_\eta=\Big\{\bbm \in (-1,1)^n:\;W(\hat{\mu}_{\bbm},\mu_\star)<\eta\Big\}.\end{align}

The proof of Theorem \ref{thm:local} then consists of four steps (all
conditional on $\bx=\ones$):
\begin{enumerate}
\item For sufficiently small $\delta>0$, we use the Sudakov-Fernique inequality
to lower bound the value of $\cF_\TAP$ on
$\cB_{\delta} \setminus \cB_{\delta/2}$. Comparing with the value of $\cF_\TAP$
achieved by an iterate $\bbm^k \in \cB_{\delta/2}$ of AMP,
we show that $\cF_\TAP$ must have a local minimizer $\bbm_\star$ in
$\cB_\delta$, and $\cF_\TAP(\bbm_\star) \approx e_\star$.\\
\item For any fixed $\eta>0$,
we use a Kac-Rice upper bound to show that with high probability, any such
local minimizer $\bbm_\star$ cannot belong to $\cB_\delta \setminus \cD_\eta$.
Thus it must belong to $\cB_\delta \cap \cD_\eta$.\\
\item For $\eps,t>0$ sufficiently small, we apply a second Kac-Rice
upper bound to show that for all critical points $\bbm_\star \in
\cD_\eta$, $\lambda_{\min}(n \cdot \nabla^2 \cF_\TAP)\geq t$ everywhere in a
$\sqrt{\eps n}$-ball around $\bbm_\star$.
We analyze the Kac-Rice bound by
representing $\lambda_{\min}(n \cdot \nabla^2 \cF_\TAP)$ over this ball as the infimum
of a Gaussian process, and lower bounding its value by a second application of
the Sudakov-Fernique inequality.

This implies that $\cF_\TAP$ is strongly convex near any local minimizer
$\bbm_\star \in \cB_\delta \cap \cD_\eta$ of Steps 1 and 2. This convexity then
ensures that there exists a unique such local minimizer satisfying
(\ref{eq:bayesoptimal}), establishing Theorem \ref{thm:local}(a--b).\\
\item To show Theorem \ref{thm:local}(c), we relate each (possibly complex)
eigenvalue $\mu$ of $\de T_{\AMP}(\bbm_\star,\bbm_\star)$
to a zero eigenvalue of a corresponding
``Bethe Hessian'' of $\cF_{\TAP}$ \cite{saade2014spectral}.
We extend the Kac-Rice/Sudakov-Fernique
argument of Step 3 from $\nabla^2 \cF_\TAP$ to this Bethe Hessian, and show
that it is positive definite whenever $|\mu|$ exceeds some constant
$r(\lambda) \in (0,1)$. Thus all
eigenvalues of $\de T_{\AMP}$ satisfy $|\mu| \leq r(\lambda)$.
\end{enumerate}

The next four sections describe these steps in greater detail.

\subsection{Sudakov-Fernique lower bound for the TAP free energy}\label{sec:TAP-lb}

We record here the following application of the
Slepian/Sudakov-Fernique comparison inequality for Gaussian processes.

\begin{lemma}\label{lem:slepian}
Let $\cX$ be a separable metric space, and let $f:\cX \to \R$ and $\bv:\cX \to
\R^n$ be bounded measurable functions on $\cX$. Let $\bW \sim \GOE(n)$ and
$\bg \sim \cN(\bzero,\id_n)$. Then
\[\E\left[\sup_{x \in \cX} \bv(x)^\top \bW \bv(x)+f(x)\right]
\leq \E\left[\sup_{x \in \cX} \frac{2}{\sqrt{n}}\|\bv(x)\|_2
\langle \bg, \bv(x) \rangle+f(x)\right].\]
\end{lemma}

Note that (conditional on $\bx=\ones$)
$-\cF_{\TAP}(\bbm)$ is a Gaussian process of this form, where
$\cX=(-1,1)^n$ and $\bv(\bbm)=\sqrt{\lambda/2n} \cdot \bbm$. Then applying this
comparison lemma and an analysis of the comparison process,
we obtain the following lower bound for $\cF_{\TAP}(\bbm)$ in terms of a
low-dimensional, deterministic variational formula.
\begin{lemma}\label{lem:global_min}
	Fix any $\lambda>1$, and suppose $\bx=\ones$. Fix any $\eps>0$ and
	two compact sets $K \subseteq [0,1]^2 \times [0,\log 2]$
	and $K' \subset \reals^3$.
	Then for some $(\lambda,K',\eps)$-dependent constant $c>0$ and all large $n$,
	with probability at least $1-e^{-cn}$,
	\begin{equation}\label{eq:valuelowerbound}
	\inf_{\bbm \in (-1,1)^n:\,(Q(\bbm),M(\bbm),H(\bbm)) \in K}\;\cF_{\TAP}(\bbm)
	>\inf_{(q,\vphi,h) \in K} \sup_{(\gamma, \tau, \nu) \in K'}
	E_{\lambda}(q, \vphi, h; \gamma, \tau, \nu)-\eps
	\end{equation}
	where
	\begin{align}
	E_{\lambda}(q, \vphi, h;\gamma, \tau, \nu)&=-\frac{\lambda^2}{2}
	\vphi^2-\frac{\lambda^2}{4}(1 - q)^2 - h + \frac{q \gamma}{2} + \vphi \tau 
	+\nu h\nonumber\\
	&\hspace{0.2in}-\E_{G \sim \cN(0, 1)} \Big\{ \sup_{m \in (-1,1)}
	\Big[ \lambda \sqrt{q} \cdot Gm + \frac{\gamma m^2}{2}
	+ \tau m + \nu\sh(m) \Big] \Big\}.\label{eq:Elambda}
	\end{align}
\end{lemma}
Lemma \ref{lem:global_min} makes precise the statement that
\begin{equation}\label{eq:barE}
\bar{E}_\lambda(q,\varphi,h)=\sup_{(\gamma,\tau,\nu) \in K'}
E_\lambda(q,\varphi,h;\gamma,\tau,\nu)
\end{equation}
is a lower bound for $\cF_\TAP(\bbm)$ when
$Q(\bbm) \approx q$, $M(\bbm)\approx \varphi$, and $H(\bbm) \approx h$. We may
show that $\bar{E}_\lambda(q,\varphi,h)$ 
has a local minimizer at $(q,\varphi,h)=(q_\star,q_\star,h_\star)$ and is
strongly convex around this minimizer, and hence give a more explicit lower
bound for $\cF_\TAP(\bbm)$ when $\bbm \in \cB_\delta$ for sufficiently small
$\delta>0$. 

\begin{lemma}\label{lem:global_min_analysis}
	Fix any $\lambda>1$, and let $E_{\lambda}(q,\varphi,h;\gamma,\tau,\nu)$ be
	as defined in Lemma \ref{lem:global_min}. Then
	\begin{equation}\label{eq:global_min_point}
	\sup_{(\gamma,\tau,\nu) \in \reals^3}
	E_{\lambda}(q_\star,q_\star,h_\star;\gamma,\tau,\nu)=E_{\lambda}(q_\star,q_\star,h_\star;0,\lambda^2
	q_\star,1)=e_\star.
	\end{equation}
	Fix any subset $K' \subseteq \reals^3$ containing
	$(0,\lambda^2q_\star,1)$ in its interior, and define $\bar{E}_\lambda$
by (\ref{eq:barE}). Then for some
	$\lambda,K'$-dependent constants $\delta,c>0$ and all $(q,\varphi,h)$
	satisfying $|q-q_\star|,|\vphi-q_\star|,|h-h_\star| \leq \delta$,
	\begin{equation}\label{eq:global_min_bound}
	\bar{E}_{\lambda}(q,\vphi,h) \geq e_\star
	+c(q-q_\star)^2+c(\vphi-q_\star)^2+c(h-h_\star)^2.
	\end{equation}
\end{lemma}

Lemmas \ref{lem:global_min} and \ref{lem:global_min_analysis} together imply
that the energy value $\cF_\TAP(\bbm)$ is bounded away from $e_\star$ on the
domain $\bbm \in \cB_\delta \setminus \cB_{\delta/2}$. The AMP state evolution
may be applied to show that AMP iterates eventually enter $\cB_{\delta/2}$, and
achieve a TAP free energy value arbitrarily close to $e_\star$ (c.f. Lemma \ref{lem:AMP}). Combined, these
yield the following corollary.

\begin{corollary}\label{cor:critexists}
Fix any $\lambda>1$ and $\delta>0$, and suppose $\bx=\ones$.
Then with probability approaching 1 as $n \to \infty$,
there exists a critical point and local minimizer $\bbm_\star$ of $\cF_{\TAP}$
belonging to $\cB_\delta$ and satisfying
$|\cF_{\TAP}(\bbm_\star)-e_\star|<\delta$.
\end{corollary}

The detailed proofs of this section are contained in Appendix \ref{appendix:TAP-lb-proofs}.

\begin{figure}\label{fig:global_min}
\centering
\includegraphics[width = 0.32\linewidth]{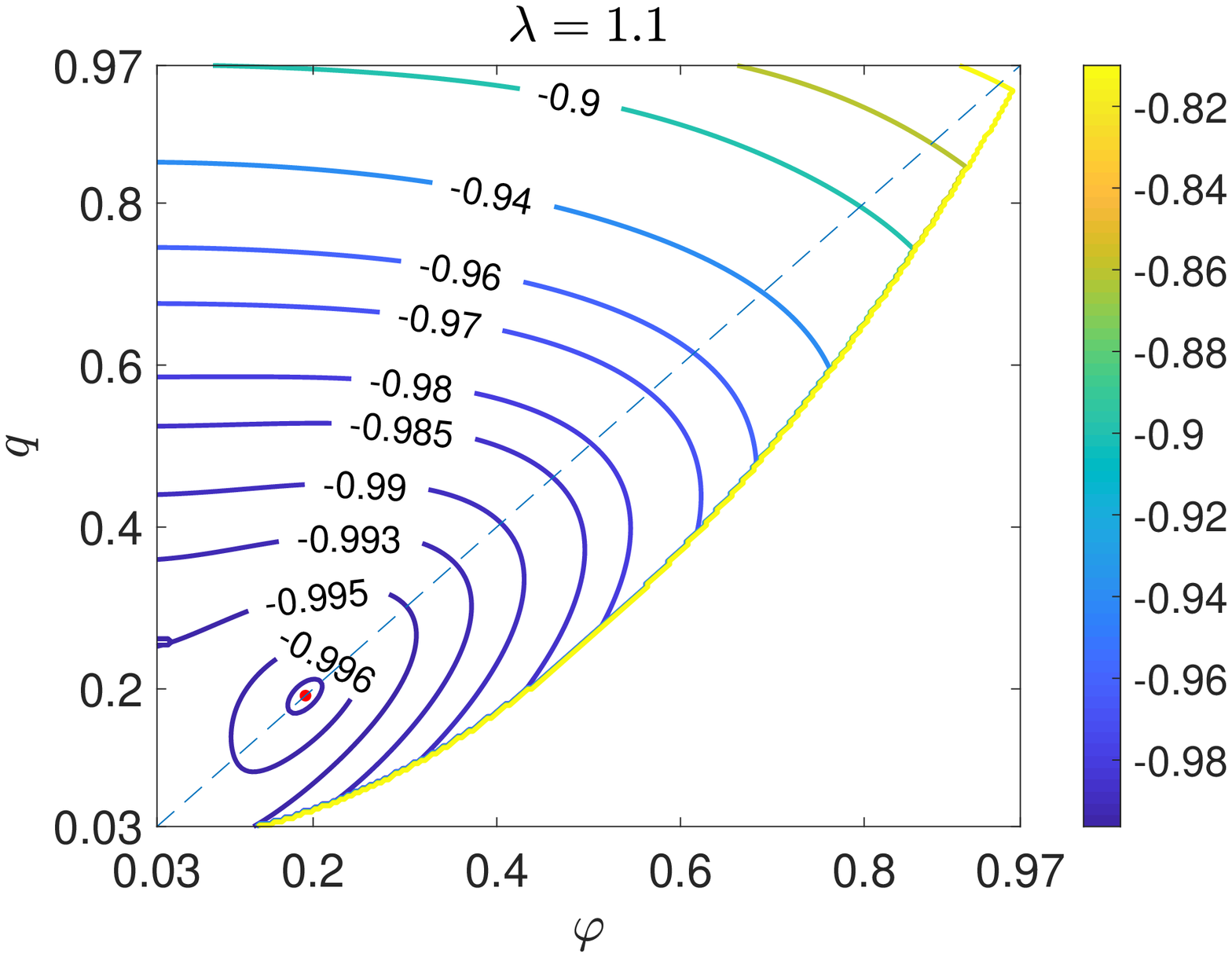}
\includegraphics[width = 0.32\linewidth]{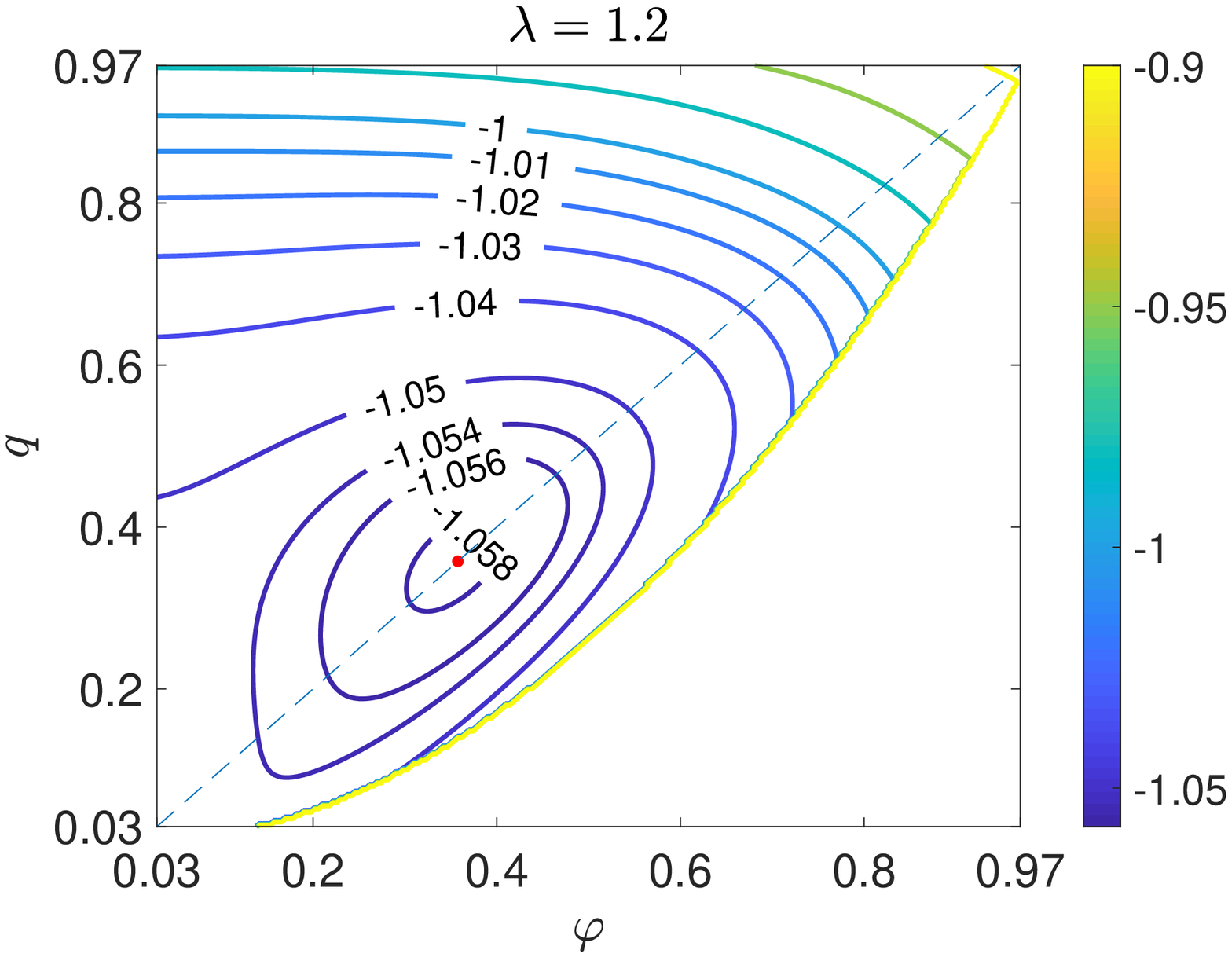}
\includegraphics[width = 0.32\linewidth]{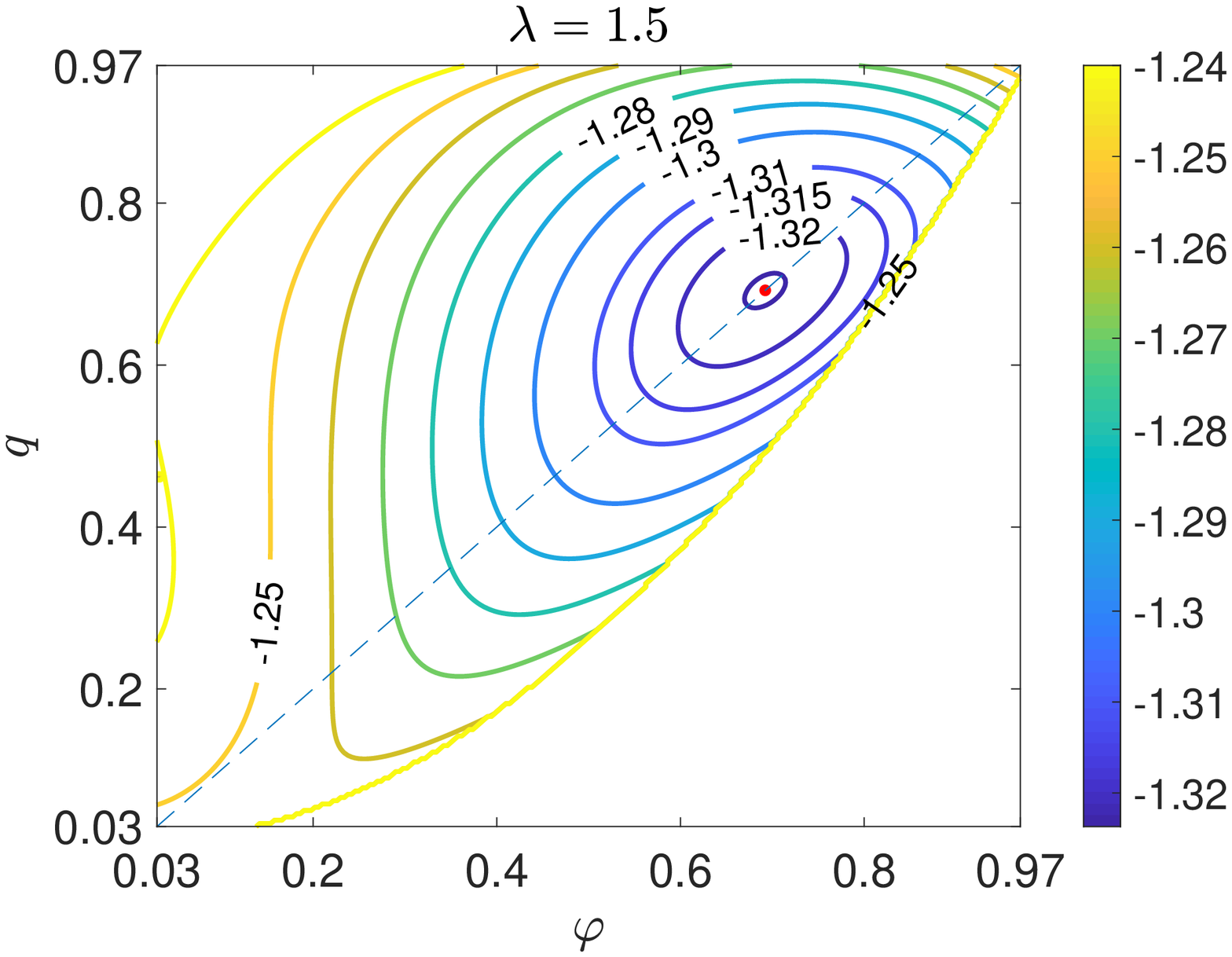}
\caption{The contour plot of the function $\bar{E}_\lambda(q,\varphi)$ as
defined in Eq. (\ref{eqn:bar_E_q_phi}). Here we take $\lambda = 1.1, 1.2, 1.5$. The global
minimum is at $(q, \varphi) = (q_\star(\lambda), q_\star(\lambda))$ where
$q_\star(1.1) \approx 0.1917$, $q_\star(1.2) \approx 0.3577$, $q_\star(1.5) \approx 0.6923$. The dashed line is $q = \vphi$. }
\end{figure}

\begin{remark}\label{remark:globalmin}
We conjecture, based on numerical evidence,
that $(q_\star,q_\star,h_\star)$ is in fact the global
minimizer of $\bar{E}_\lambda(q,\varphi,h)$ for all $\lambda>1$:
We may first restrict $E_\lambda$ to $\nu=1$ and
$\tau=\lambda^2 \varphi$, to obtain the further lower bound
\begin{equation}\label{eqn:bar_E_q_phi}
\bar{E}_\lambda(q,\varphi,h) \geq \bar{E}_\lambda(q,\varphi)
=\sup_\gamma E_\lambda(q,\varphi;\gamma)\end{equation}
where
\[E_\lambda(q,\varphi;\gamma)=\frac{\lambda^2}{2}\varphi^2
-\frac{\lambda^2}{4}(1-q)^2+\frac{q\gamma}{2}
-\E_{G \sim \cN(0,1)}\Big[\sup_{m \in (-1,1)} \lambda \sqrt{q} \cdot Gm
+\frac{\gamma m^2}{2}+\lambda^2 \varphi m+\sh(m)\Big].\]
Numerical evaluations of this function
$\bar{E}_\lambda(q,\varphi)$ over the relevant domain
$q \in (0,1)$ and $|\varphi|<\sqrt{q}$ are presented in Figure
\ref{fig:global_min}. For all tested values of $\lambda>1$, these evaluations support the claim
that $\bar{E}_\lambda(q,\varphi)$ has the unique \emph{global} minimizer
$(q,\varphi)=(q_\star,q_\star)$. This claim then implies
that $(q_\star,q_\star,h_\star)$ is also the
unique global minimizer of $\bar{E}_\lambda(q,\varphi,h)$, by the global 
convexity of $h \mapsto \bar{E}_\lambda(q_\star,q_\star,h)$ and its strong
convexity near its minimizer $h_\star$.

Subject to the validity of this numerical conjecture, Lemma \ref{lem:global_min}
may be used to show that $\cF_\TAP(\bbm)$ is also bounded away from $e_\star$
for all $\bbm \in (-1,1)^n \setminus \cB_\delta$. Our subsequent arguments will
then imply that for any $\lambda>1$, with probability approaching 1, the
(unique) local minimizer $\bbm_\star$ described by Corollary
\ref{cor:critexists} and Theorem \ref{thm:local} is in fact the global
minimizer of $\cF_\TAP$. (All theoretical results stated in this work will
be established using only that $\bbm_\star$ is a local minimizer of $\cF_\TAP$, and
they will not require the validity of this conjecture.)
\end{remark}

\subsection{Kac-Rice localization of critical points}\label{sec:KacRice}

We now use a Kac-Rice upper bound to show that the critical point(s)
$\bbm_\star$ described by Corollary \ref{cor:critexists} must belong to the
more restrictive set $\cB_\delta \cap \cD_\eta$ (c.f. Eq. (\ref{eqn:set_B_delta}) and (\ref{eqn:set_D_eta})).

Define functions $\bg$ and $\bH$, which are the gradient and Hessian of the renormalized TAP free energy
\begin{align}
	\bg(\bbm)&=n \cdot \nabla \cF_{\TAP}(\bbm)
	=-\lambda \bY \bbm+\arctanh(\bbm)+\lambda^2[1 -
Q(\bbm)]\bbm,\label{eq:grad}\\
	\bH(\bbm)&=n \cdot \nabla^2 \cF_{\TAP}(\bbm)=-\lambda\bY+\diag\Big(\frac{1}{1 - \bbm^2} \Big)+\lambda^2 [1 - Q(\bbm)] \id
	-\frac{2\lambda^2}{n}\bbm \bbm^\sT.\label{eq:hess}
\end{align}
\noindent We apply the following Kac-Rice upper bound from \cite{fan2021tap}.

\begin{lemma}\label{lem:KacRicestandard}
Fix any $\lambda>0$, suppose $\bx=\ones$, and let $T \subseteq (-1,1)^n
\setminus \{\bzero\}$ be any (deterministic) Borel-measurable set. Then
\[\E\Big[\big|\{\bbm \in T:\bg(\bbm)=\bzero\}\big|\Big] \leq
\int_T \E\Big[|\det \bH(\bbm)|\;\Big|\; \bg(\bbm)=\bzero\Big]\,
p_{\bg(\bbm)}(\bzero)\de \bbm\]
where $p_{\bg(\bbm)}(\bzero)$ is the Lebesgue-density of the distribution of
$\bg(\bbm)$ at $\bg(\bbm)=\bzero$.
\end{lemma}

Applying this bound, we eliminate the possibility that the critical point(s)
described by Corollary \ref{cor:critexists} belong to
$\cB_\delta \setminus \cD_\eta$, as stated in the following
lemma. Thus they belong to $\cB_\delta \cap \cD_\eta$ as desired.

\begin{lemma}\label{lem:localizetoDeta}
Fix any $\lambda>1$ and $\eta>0$, and suppose $\bx=\ones$. Then for some
$(\lambda,\eta)$-dependent constants $c,\delta>0$ and all large $n$,
\[\P\Big[\text{ there exists } \bbm \in
\cB_\delta:\;\bg(\bbm)=\bzero,\;|\cF_{\TAP}(\bbm)-e_\star|<\delta,\;
\bbm \notin \cD_\eta\;\Big]<e^{-cn}.\]
\end{lemma}

Let us make two high-level clarifications regarding the proof:
First, to show Lemma \ref{lem:localizetoDeta}, we wish to apply Lemma
\ref{lem:KacRicestandard} with $T$ being the set
\[\Big\{\bbm \in \cB_\delta \setminus \cD_\eta:
|\cF_{\TAP}(\bbm)-e_\star|<\delta\Big\}.\]
We cannot do so directly, because $\cF_{\TAP}(\bbm)$ is random,
and hence this is not a deterministic subset of $(-1,1)^n$.
However, restricted to points $\bbm$
where $\bg(\bbm)=\bzero$, the identity $0=\bbm^\top \bg(\bbm)$ allows us to
re-express $\bbm^\top \bY\bbm$ and $\cF_{\TAP}(\bbm)$ as deterministic functions
of $\bbm$. Lemma \ref{lem:localizetoDeta} is then obtained by replacing
$|\cF_{\TAP}(\bbm)-e_\star|<\delta$ with an equivalent deterministic condition
to define $T$.

Second, we remark that the Sudakov-Fernique argument of the preceding section
cannot be used here to similarly localize $\bbm_\star$ to $\cD_\eta$, by
bounding the TAP free energy value outside $\cB_\delta \setminus \cD_\eta$.
This is because there exists
$\bbm \in (-1,1)^n$ with one coordinate very close to $\pm 1$, so that
$W(\hat{\mu}_{\bbm},\mu_\star)$ is arbitrarily large (c.f. Eq. (\ref{eqn:W_distance})) and $\bbm \notin
\cD_\eta$, but $\cF_\TAP(\bbm)$ is arbitrarily close to $e_\star$ in
value. Thus, we use this separate Kac-Rice argument, and the
condition $\bbm \in \cB_\delta$ as an input to the Kac-Rice analysis,
to establish the localization to $\cD_\eta$ in Lemma \ref{lem:localizetoDeta}.

The detailed proofs of this section are contained in Appendix \ref{appendix:KacRiceproofs}.

\subsection{Sudakov-Fernique lower bound for local strong
convexity}\label{sec:strongconvex}

We now show that the TAP free energy is strongly convex in a local neighborhood
of any critical point $\bbm_\star \in \cD_\eta$.
For a parameter $\eps>0$, define
\begin{equation}\label{eq:elleps}
\ell_\eps^+(\bbm,\bW)=\inf\Big\{\lambda_{\min}\big(\bH(\bu)\big):
\bu \in (-1,1)^n \cap \ball_{\sqrt{\eps n}}(\bbm)\Big\}.
\end{equation}
The dependence of $\ell_\eps^+$ on $\bW$ is via $\bH(\bu)$. We will make this
dependence implicit in what follows, and write simply
$\ell_\eps^+(\bbm)=\ell_\eps^+(\bbm,\bW)$. If $\ell_\eps^+(\bbm_\star) \geq
t>0$, then the TAP free energy is strongly convex on a $\sqrt{\eps n}$-ball
around $\bbm_\star$, as desired. We use a Kac-Rice upper bound to show,
with high probability, no critical points of $\cF_\TAP$ belong to
the set
\[\Big\{\bbm \in \cD_\eta:\ell_\eps^+(\bbm)<t\Big\}\]
for some sufficiently small constant $t>0$.

The condition $\ell_\eps^+(\bbm)<t$ is again random, so this is not a
deterministic subset of $(-1,1)^n$. We address this using the following
extension of the Kac-Rice upper bound in Lemma \ref{lem:KacRicestandard}.

\begin{lemma}\label{lem:KacRicebound}
Fix any $\lambda>0$, and suppose $\bx=\ones$.
Let $\Sym_n$ be the space of real symmetric $n \times n$
matrices, $T \subseteq (-1,1)^n \setminus \{\bzero\}$ any (deterministic)
Borel-measurable set,
and $\ell:T \times \Sym_n \to \reals$ any Borel-measurable function. Let $c>0$
and $t \in \reals$ be any (possibly $n$-dependent) values, and let
$U \sim \Unif([-c,c])$ be a uniform random variable independent of $\bW$.
Define
\[\cC=\Big\{\bbm \in T:\bg(\bbm)=\bzero \text{ and } \ell(\bbm,\bW)+U<t\Big\}.\]
Then
\begin{equation}\label{eq:generalKacRice}
\E[|\cC|] \leq \int_T \E\Big[|\det \bH(\bbm)| \cdot \ones\{\ell(\bbm,\bW)+U<t\}
\;\Big|\;\bg(\bbm)=\bzero\Big]\;p_{\bg(\bbm)}(\bzero)\,\de \bbm,
\end{equation}
where $p_{\bg(\bbm)}(\bzero)$ is the Lebesgue-density of the distribution of
$\bg(\bbm)$ at $\bg(\bbm)=\bzero$, and the expectations are over both $U$ and
$\bW$.
\end{lemma}

\noindent (Introducing this auxiliary variable $U$ alleviates the
need to check a technical condition that $\ell(\bbm,\bW)=t$ and
$\bg(\bbm)=\bzero$ do not simultaneously occur at any $\bbm \in T$, when
applying the Kac-Rice lemma.)

In \cite{fan2021tap}, an upper bound on the determinant $|\det \bH(\bbm)|$ was
established via a spectral analysis of $\bH(\bbm)$, which shows
$\E[|\det \bH(\bbm)|^2 \mid \bg(\bbm)=\bzero] \leq e^{c(\eta)n}$ for
$\bbm \in \cD_\eta$ and a constant $c(\eta) \to 0$ as $\eta \to 0$.
Thus, to show that (\ref{eq:generalKacRice}) is vanishing,
we complement this by showing an exponentially small upper bound for the
probability $\P[\ell_\eps^+(\bbm)+U<t \mid \bg(\bbm)=\bzero]$. We do
this again using the Sudakov-Fernique inequality of Lemma \ref{lem:slepian}, to
obtain the variational lower bound on the conditional mean
$\E[\ell_\eps^+(\bbm) \mid \bg(\bbm)=\bzero]$ stated in part (a) of the
following lemma. This bound is shown to be positive in part (b).

\begin{lemma}\label{lem:gordon-variational}
Suppose $\lambda>1$ and $\bx=\ones$. Define
	\begin{align}
	H_\lambda^+&(p,u;\alpha,\kappa,\gamma) = -\Big[2\lambda^2p^2+\lambda^2
	u^2-2\lambda^2(1 - q_\star) p^2 / q_\star- \alpha u  - \kappa p \Big]
	+\lambda^2(1-q_\star) + \gamma \nonumber\\
		&- \E_{m \sim \mu_\star}\Big[ \Big(4 \lambda^2  (1 - p^2 / q_\star) +
	(2 z(m)  p / q_\star + \alpha   + \kappa m)^2 \Big) \Big/ \Big( \frac{4}{1 -
	m^2} - 4\gamma \Big) \Big]
	\label{eq:Hlambda}
	\end{align}
	where $z(m) = \arctanh m -  \lambda^2 q_\star + \lambda^2 (1 - q_\star)m$.

\begin{enumerate}
	\item[(a)]
Fix any $t>0$ and compact domain $K' \subset \reals^2 \times (-\infty,1)$.
For some $(\lambda,K',t)$-dependent constants $\eps,\eta>0$, and all large $n$,
	\[
		\inf_{\bbm \in \cD_\eta}
		\E\Big[\;\ell_\eps^+(\bbm) \;\Big|\; \bg(\bbm)=\bzero\;\Big] \geq
		\mathop{\inf_{u \in [-1,1]}}_{p  \in [-\sqrt{q_\star},
		\sqrt{q_\star}]} \sup_{(\alpha,\kappa,\gamma) \in K'}
	H_\lambda^+(p,u;\alpha,\kappa,\gamma)-t
	\]
	\item[(b)]
	Suppose $K'$ contains $(0,0,0)$ in its interior. Then there is a
	$(\lambda,K')$-dependent constant $t_0>0$ for which
	\[\mathop{\inf_{u \in [-1,1]}}_{p  \in [-\sqrt{q_\star},
		\sqrt{q_\star}]} \sup_{(\alpha,\kappa,\gamma) \in K'}
	H_\lambda^+(p,u;\alpha,\kappa,\gamma)>t_0>0.\]

\end{enumerate}
\end{lemma}

The desired upper bound for $\P[\ell_\eps^+(\bbm)+U<t \mid \bg(\bbm)=\bzero]$
then follows by concentration of $\ell_\eps^+(\bbm)$ around its mean.
Applying this to (\ref{eq:generalKacRice})
yields the following corollary on local strong convexity.
\begin{corollary}\label{cor:strongconvexity}
Fix any $\lambda>1$, and suppose $\bx=\ones$. Then there exist
$\lambda$-dependent constants $\eps,\eta,t,c>0$ such that, for all large $n$,
\begin{align}
\P\Big[\text{ there exist } \bbm \in \cD_\eta \text{ and }
\bu \in (-1,1)^n &\cap \ball_{\sqrt{\eps n}}(\bbm): \bg(\bbm)=\bzero \text{ and } \lambda_{\min}(\bH(\bu))<t\;\Big]<e^{-cn}.
\label{eq:stronglyconcave}
\end{align}
\end{corollary}

Finally, this convexity implies Theorem \ref{thm:local}(a--b)
by the following argument: Letting $\bbm^k$ be a sufficiently large iterate of
AMP, we may pick a local minimizer $\bbm_\star$
in Corollary \ref{cor:critexists} such that there is a strict descent path from
$\bbm^k$ to $\bbm_\star$. Strong convexity of $\cF_\TAP$ around $\bbm_\star$
and an upper bound on $\cF_\TAP(\bbm^k)-\cF_\TAP(\bbm_\star)$
then imply an upper bound on the Euclidean distance $\|\bbm_\star-\bbm^k\|_2$.
Then this point $\bbm_\star$ must satisfy (\ref{eq:bayesoptimal}) by the
Bayes-optimality of the AMP iterate $\bbm^k$.
Furthermore, the local convexity of $\cF_\TAP$ implies that such a 
point $\bbm_\star$ is unique. We provide the details of this argument
in Appendix \ref{appendix:strongconvexproofs}.

\subsection{Local stability of AMP}\label{sec:AMP-stab}

We now describe the proof of Theorem \ref{thm:local}(c).
Let us write the input and output of $T_{\AMP}$ in (\ref{eq:TAMP}) as
\[(\bbm_+,\bbm)=T_\AMP(\bbm,\bbm_-).\]
Differentiating by the chain rule, the Jacobian of $T_\AMP$ may be expressed as
\begin{equation}\label{eq:dTAMP}
\de T_\AMP(\bbm,\bbm_-)=\begin{pmatrix}
\diag(1-\bbm_+^2) \cdot [\lambda \bY+2\lambda^2 \bbm_-\bbm^\sT/n] &
-\diag(1-\bbm_+^2) \cdot \lambda^2[1-Q(\bbm)] \\
\id & \bzero \end{pmatrix}. 
\end{equation}
At any point $\bbm_\star \in (-1,1)^n$ where $\bg(\bbm_\star)=\bzero$, we have
$T_{\AMP}(\bbm_\star,\bbm_\star)=(\bbm_\star,\bbm_\star)$. Thus
$\de T_{\AMP}(\bbm_\star,\bbm_\star)=\bB(\bbm_\star)$ for the matrix
\[\bB(\bbm)=\begin{pmatrix}
\diag(1-\bbm^2) \cdot [\lambda \bY+2\lambda^2 \bbm\bbm^\sT/n] &
-\diag(1-\bbm^2) \cdot \lambda^2[1-Q(\bbm)] \\
\id & \bzero \end{pmatrix}.\]

In Appendix \ref{appendix:AMPstabproofs}, we first verify the simple algebraic
identity that the eigenvalues $\mu \in \C$ of this matrix $\bB(\bbm)$, for any
$\bbm \in (-1,1)^n$, are exactly those values $\mu \in \C$ for which the ``Bethe
Hessian'' matrix
\begin{equation}\label{eq:betheHessiancomplex}
\mu\Big({-}\lambda \bY-\frac{2\lambda^2}{n}\bbm\bbm^\top\Big)
+\lambda^2[1-Q(\bbm)]\bI+\mu^2\diag\Big(\frac{1}{1-\bbm^2}\Big)
\end{equation}
is singular. Applying this relation, we then show the following
deterministic lemma relating the spectral radius of $\bB(\bbm)$ to the smallest
eigenvalue of the above matrix for real arguments $\mu=\pm r$.

\begin{lemma}\label{lem:stabilitycondition}
Fix any $\lambda>1$. There exist $\lambda$-dependent constants $\delta>0$ and
$r_0 \in (0,1)$ such that for any $r \in (r_0,1)$ and $\bbm \in (-1,1)^n$ with
$|Q(\bbm)-q_\star|<\delta$, if we have
\begin{equation}\label{eq:betheHessian}
\lambda_{\min}\bigg[\pm r\Big({-}\lambda \bY-\frac{2\lambda^2}{n}\bbm\bbm^\sT\Big)
+\lambda^2[1-Q(\bbm)]\id+r^2\diag\Big(\frac{1}{1-\bbm^2}\Big)\bigg]>0
\end{equation}
for both choices of sign $\pm$, then $\rho(\bB(\bbm))<r<1$.
\end{lemma}

To prove Theorem \ref{thm:local}(c), by a simple continuity argument,
it will suffice to consider exactly $r=1$ in (\ref{eq:betheHessian}) and
to show that (\ref{eq:betheHessian}) holds with high probability
at $\bbm=\bbm_\star$ for both choices of sign $\pm$. For $r=1$ and sign $+$,
the matrix in (\ref{eq:betheHessian}) is precisely the Hessian 
$\bH(\bbm)$, whose smallest eigenvalue at $\bbm=\bbm_\star$ was bounded in the
preceding section. The case of sign $-$ is a minor extension of these
arguments: Define
\begin{align*}
\bH^-(\bbm)&=\left(\lambda \bY+\frac{2\lambda^2}{n}\bbm\bbm^\sT\right)
+\diag\Big(\frac{1}{1-\bbm^2}\Big)+\lambda^2[1-Q(\bbm)]\id,\\
\ell_\eps^-(\bbm)&=\inf\Big\{\lambda_{\min}\big(\bH^-(\bu)\big):
\bu \in (-1,1)^n \cap \ball_{\sqrt{\eps n}}(\bbm)\Big\}.
\end{align*}
We show the following lemma using the Sudakov-Fernique inequality, analogously
to Lemma \ref{lem:gordon-variational}.

\begin{lemma}\label{lem:gordon-variational-gen}
Suppose $\lambda>1$ and $\bx=\ones$. Define
\[\begin{aligned}
H_\lambda^-&(p,u;\alpha,\kappa,\gamma) = \Big[2\lambda^2p^2+\lambda^2
u^2-2\lambda^2(1 - q_\star) p^2 / q_\star- \alpha u  - \kappa p \Big]
+\lambda^2(1-q_\star) + \gamma \\
&~- \E_{m \sim \mu_\star}\Big[ \Big(4 \lambda^2  (1 - p^2 / q_\star) +
(2 z(m)  p / q_\star + \alpha   + \kappa m)^2 \Big) \Big/ \Big( \frac{4}{1 -
m^2} - 4\gamma \Big) \Big]
\end{aligned}\]
where $z(m)=\arctanh m -  \lambda^2 q_\star + \lambda^2 (1 - q_\star)m$.
Then the statements of Lemma \ref{lem:gordon-variational} hold also
with $\ell_\eps^+(\bbm)$ and $H_\lambda^+$ replaced by
$\ell_\eps^-(\bbm)$ and $H_\lambda^-$.
\end{lemma}

Now applying this result in the Kac-Rice upper bound of
Lemma \ref{lem:KacRicebound} for $\ell(\bbm,\bW)=\ell_\eps^-(\bbm)$, we obtain
that (\ref{eq:betheHessian}) also holds with high probability for $r=1$ and
sign $-$, implying Theorem \ref{thm:local}(c).

The detailed proofs of this section are contained in Appendix \ref{appendix:AMPstabproofs}.

\section{Convergence of optimization algorithms}\label{sec:proof_algorithm}

In this section, we describe the main ideas in the proofs of
Theorems \ref{thm:localconvergence} and \ref{thm:globalconvergence}.
It again suffices to show that the results hold with high probability
conditional on $\bx=\ones$. The detailed proofs of this section are contained in Appendix \ref{appendix:global}.

\subsection{Convergence of NGD with local initialization} \label{sec:proof_local_convergence_NGD}

Theorem \ref{thm:localconvergence} is a consequence of the local strong
convexity of $\cF_\TAP$ established in Theorem \ref{thm:local}(b) and
the following local convergence result for the natural gradient
algorithm (\ref{eq:NGD}).  

\begin{lemma}\label{lem:NGD}
Fix any $\lambda>1$, $t>0$, and $\eps \in (0,1)$.
Consider the event where $\bbm_\star$ in
Theorem \ref{thm:local}(a) exists and is unique up to sign, and
$\|\bW\|_\op<3$ and $\lambda_{\min}(n \cdot \nabla^2 \cF_\TAP(\bbm))>t$ 
for every $\bbm \in (-1,1)^n \cap \ball_{\sqrt{\eps n}}(\bbm_\star)$.
Consider any initialization $\bbm^0=\tanh(\bh^0)$ such that
\begin{equation}\label{eq:basecaseNGD}
\cF_\TAP(\bbm^0)<\cF_\TAP(\bbm_\star)+t\eps/8, \qquad
\|\bbm^0-\bbm_\star\|_2<\sqrt{\eps n}.
\end{equation}
There exist $(\lambda,t,\eps)$-dependent constants $C,\mu,\eta_0>0$ such
that if (\ref{eq:NGD}) with any step size $\eta \in (0,\eta_0)$ is initialized at $\bbm^0$,
then on this event, for every $k \geq 1$ we have
\begin{align}
\cF_\TAP(\bbm^k)&<\cF_\TAP(\bbm_\star)
+C\left(1+\frac{\|\arctanh(\bbm^0)\|_2}{\sqrt{n}}\right)(1-\mu \eta)^k,
\label{eq:valueconvergence}\\
\|\bbm^k-\bbm_\star\|_2&<C\sqrt{n}
\left(1+\frac{\|\arctanh(\bbm^0)\|_2}{\sqrt{n}}\right)(1-\mu \eta)^k
\label{eq:mconvergence}
\end{align}
\end{lemma}

The proof of this lemma applies the mirror-descent form of NGD given in
(\ref{eq:mirrordescent}), together with an observation that on the above
event, $\cF_\TAP$ is strongly smooth and strongly convex over $(-1,1)^n \cap
\ball_{\sqrt{\eps n}}(\bbm_\star)$ relative to the prox function $-H(\bbm)$,
in the sense of \cite{bauschke2017descent,lu2018relatively}
\[\mu \cdot \nabla^2 ({-}H(\bbm)) \preceq
\nabla^2 \cF_\TAP(\bbm) \preceq L \cdot \nabla^2 ({-}H(\bbm))\]
for some constants $L,\mu>0$. We may then
adapt a convergence analysis of \cite{lu2018relatively} to show that, for the
above initialization, NGD with sufficiently small step size $\eta>0$ must
remain in this strongly convex neighborhood and exhibit the above linear
convergence to $\bbm_\star$.
For any $\lambda>1$, the event in Lemma \ref{lem:NGD} holds with high
probability by Theorem \ref{thm:local}. The required initial
condition (\ref{eq:basecaseNGD}) is also with high probability achieved by a
sufficiently large iteration of AMP, as may be deduced from the AMP state
evolution. Combined, this yields Theorem \ref{thm:localconvergence}. The detailed proofs of Lemma \ref{lem:NGD} and Theorem \ref{thm:localconvergence} are contained in Appendix \ref{appendix:analysis_NGD}.

\subsection{Convergence of NGD from spectral initialization} \label{sec:proof_NGD_global}

For large $\lambda$, to show the result of Theorem \ref{thm:globalconvergence}(b)
that NGD alone converges to $\pm \bbm_\star$ from a
spectral initialization, recall the domain
\[\cS=\Big\{\bbm \in (-1,1)^n:\cF_\TAP(\bbm)<-\lambda^2/3\Big\}\]
as defined in Corollary \ref{cor:globallandscape}. For a parameter
$q \in (0,1)$,
define the deterministic subset
\[\cM_q=\Big\{\bbm \in (-1,1)^n:M(\bbm)>q\Big\},\]
where recall $M(\bbm)=\bbm^\top \ones/n$.
We first establish the following more quantitative characterization of 
the landscape of $\cF_\TAP$.

\begin{lemma}\label{lem:quantitativelandscape}
Fix any integer $a \geq 5$, and set $q=1-\lambda^{-a}$. Suppose $\bx=\ones$.
For a constant $\lambda_0(a)>0$, if $\lambda>\lambda_0(a)$, then there
are $(a,\lambda)$-dependent constants $C,c,t>0$ such that
with probability at least $1-Ce^{-cn}$,
\begin{enumerate}
\item[(a)] Every point $\bbm \in \cS \setminus \cM_q$ satisfies
\[\big\|\sqrt{n} \cdot \nabla \cF_\TAP(\bbm)\big\|_2^2>t.\]
\item[(b)] Every point $\bbm \in \cS \cap \cM_q$ satisfies
\[n \cdot \nabla^2 \cF_\TAP(\bbm) \succ
\frac{1}{2}\,\diag\left(\frac{1}{1-\bbm^2}\right) \succeq \frac{1}{2}\,\id.\]
\end{enumerate}
\end{lemma}

Part (b) of this lemma is sufficient to imply
Theorem \ref{thm:globalconvergence}(b) on the convergence of NGD: The
initialization $\bbm^0=\tanh(\bh^0)$ defined by (\ref{eq:initialization})
will belong to the region $\cS \cap \cM_q$ with high probability, so
Lemma \ref{lem:NGD} may again be used to show linear convergence to
$\pm \bbm_\star$. The detailed proof of Lemma \ref{lem:quantitativelandscape} is contained in Appendix \ref{appendix:TAP_global_landscape} and the detailed proof of Theorem \ref{thm:globalconvergence}(b) is contained in Appendix \ref{appendix:analysis_NGD}.

\subsection{Convergence of AMP from spectral initialization}\label{sec:proof_AMP_global}

For Theorem \ref{thm:globalconvergence}(a) on the convergence of AMP, we
directly prove contractivity of the map $T_\AMP$ defined in (\ref{eq:TAMP})
locally near $\bbm_\star$, in a parameterization by coordinates $\bp$
that lie ``between'' $\bh$ and $\bbm=\tanh(\bh)$: Define two strictly increasing
functions $\Gamma,\Lambda: \R \to \R$ as
\[ \Gamma(h) = \int_0^h \sqrt{1 - \tanh(s)^2}\, \de s, 
\qquad \Lambda(p)=\tanh(\Gamma^{-1}(p)),\]
and consider $\bp=\Gamma(\bh)$. Then $\bbm=\tanh(\bh)=\Lambda(\bp)$.
We write as shorthand
\[\frac{\de \bbm}{\de \bh}=\de_h \tanh(\bh), \qquad
\frac{\de \bp}{\de \bh}=\de_h \Gamma(\bh), \qquad
\frac{\de \bbm}{\de \bp}=\de_p \Lambda(\bp)\]
where these are vectors in $\R^n$, and the derivatives are applied entry-wise.
These definitions of $\Gamma$ and $\Lambda$ are designed so as to
factor the identity $1-\bbm^2=\de \bbm/\de \bh$ into the pair of identities
\[\sqrt{1-\bbm^2}=\frac{\de \bbm}{\de \bp}=\frac{\de \bp}{\de \bh}.\]
 (This reparameterization by $\bp$ may seem mysterious, and is 
carefully chosen to precondition the Jacobian of the AMP map and enable an
operator norm bound for this Jacobian. We provide a heuristic motivation for
this reparametrization in Remark \ref{rmk:motivation} in Appendix
\ref{appendix:AMPglobal}.)

The range of $\bp=\Gamma(\bh)$ is the cube $\Omega^{(p)}=(-\pi/2,\pi/2)^n$.
We denote the AMP map (\ref{eq:TAMP}) in the $\bp$-parameterization as
$T_\AMP^{(p)}:\Omega^{(p)} \times \Omega^{(p)} \to
\Omega^{(p)} \times \Omega^{(p)}$, defined by
\[T_\AMP^{(p)}(\bp,\bp_-)=(\Lambda \otimes \Lambda)^{-1} \circ
T_\AMP\big((\Lambda \otimes \Lambda)(\bp,\bp_-)\big).\]
Thus, reparameterizing by $\bp^k=\Gamma(\bh^k)$, the AMP iterations
(\ref{eq:AMP}) take the form $(\bp^{k+1},\bp^k)=T_\AMP^{(p)}(\bp^k,\bp^{k-1})$.

\begin{lemma}\label{lem:AMPcontractive}
Consider the metric $\|(\bp,\bp')\|_\lambda=\|\bp\|_2+\lambda^{-1/5}\|\bp'\|_2$.
Fix $q=1-\lambda^{-5}$ and $\bx=\ones$. For an absolute constant $\lambda_0>0$,
suppose $\lambda>\lambda_0$. Then with probability at least $1-Ce^{-cn}$
for $\lambda$-dependent constants $C,c>0$, the following holds:
If there exists a critical point $\bbm_\star \in \cM_q$ of $\cF_\TAP$,
then for $\bp_\star=\Lambda^{-1}(\bbm_\star)$, any
$\bp,\bp_- \in \ball_{\lambda^{-7}\sqrt{n}}(\bp_\star) \cap \Omega^{(p)}$,
and $(\bp_+,\bp)=T_\AMP^{(p)}(\bp,\bp_-)$,
we have $\bp_+ \in \ball_{\lambda^{-7}\sqrt{n}}(\bp_\star) \cap \Omega^{(p)}$
and
\begin{equation}\label{eq:AMPcontractive}
\big\|(\bp_+,\bp)-(\bp_\star,\bp_\star)\big\|_\lambda
\leq 2\lambda^{-1/5}\big\|(\bp,\bp_-)-(\bp_\star,\bp_\star)\big\|_\lambda.
\end{equation}
\end{lemma}

The AMP state evolution guarantees that with probability approaching 1 as $n \to \infty$,
$\bp^{k-1},\bp^k \in \ball_{\lambda^{-7}\sqrt{n}}(\bp_\star)$
for a sufficiently large iteration $k$. Then the
contractivity guaranteed in Lemma \ref{lem:AMPcontractive} implies Theorem
\ref{thm:globalconvergence}(a). The detailed proofs of Lemma \ref{lem:AMPcontractive} and Theorem
\ref{thm:globalconvergence}(a) are contained in Appendix \ref{appendix:AMPglobal}.

\section{Discussion}

In this paper, we showed the local strong convexity of the TAP free energy
for $\Z_2$-synchronization around its Bayes-optimal local minimizer, 
and studied the finite-$n$ convergence of optimization algorithms for computing
this minimizer. 
Numerical simulations confirm that the TAP free energy can be
efficiently optimized, and that properties of its minimizer are robust to model
misspecification.
Our results provide theoretical justification for using the TAP
free energy to perform variational inference in this model.

In terms of proof techniques, our work introduced a method
of using the Kac-Rice formula to study
the local geometry of a non-convex function around its critical points. Some
intermediate results in the proof, for example the convergence of the empirical
distribution of coordinates of the TAP minimizer, are of independent interest.
We note that an analogous TAP free energy function may be defined
in broader contexts, such as for spiked matrix models with more general priors
or for linear and generalized linear models, and some of our techniques may be
useful also for analyzing the local geometries of these TAP free energy
functions around their
informative fixed points. However, the Rademacher $\{+1,-1\}$ prior in
$\Z_2$-synchronization does have several conveniences, including a fixed second
moment, an explicit form for both its entropy and its posterior mean
function, and a unique fixed-point for the equation (\ref{eq:fixedpointq}) that
defines $q_\star$. Analyses of models having priors that lack these properties
would have additional technical hurdles, and we leave the exploration of such
extensions to future work.

Finally, we proved the finite-$n$ convergence of a well-studied AMP
algorithm for this problem, which is not implied by analysis of the AMP state evolution alone. 
Our proof of this result required sufficiently large $\lambda$, but we
conjecture that the result holds for any $\lambda>1$. This conjecture is supported by
our numerical simulations and also by the stability of the AMP map around its fixed
point, which indeed holds for any $\lambda>1$. We leave this conjecture as an open question, and hope that the
techniques developed in this paper can perhaps inspire a proof. 

\section*{Acknowledgement}

M. Celentano is supported by the Miller Institute for Basic Research in Science, University of California Berkeley. Z. Fan is supported in part by NSF Grants DMS-1916198 and DMS-2142476. S. Mei is supported in part by NSF Grant DMS-2210827. 

\bibliography{references}{}
\bibliographystyle{alpha}
\newpage

\appendix

\section{Preliminaries}

\subsection{Uniform continuity of $\cF_\TAP$}

\begin{proposition}\label{prop:TAPcontinuous}
For any $\bbm,\bbm' \in (-1,1)^n$,
\[|Q(\bbm)-Q(\bbm')| \leq \sqrt{2} \cdot \frac{\|\bbm-\bbm'\|_2}{\sqrt{n}},
\qquad |H(\bbm)-H(\bbm')| \leq (\log 2+1)\left(\frac{\|\bbm-\bbm'\|_2^2}{n}\right)^{1/4},\]
\[|\cF_\TAP(\bbm)-\cF_\TAP(\bbm')|
\leq \left(\lambda \|\bY\|_\op+\frac{\lambda^2\sqrt{2}}{2}\right)
\cdot \frac{\|\bbm-\bbm'\|_2}{\sqrt{n}}
+(\log 2+1)\left(\frac{\|\bbm-\bbm'\|_2^2}{n}\right)^{1/4}.\]
\end{proposition}
\begin{proof}
First, by Cauchy-Schwarz,
\[|Q(\bbm)-Q(\bbm')|
=\frac{|\langle \bbm-\bbm',\bbm+\bbm' \rangle|}{n}
\leq \sqrt{2} \cdot \frac{\|\bbm-\bbm'\|_2}{\sqrt{n}}.\]
Next, by Markov's inequality, for any $t>0$,
\[\frac{1}{n}\sum_{i=1}^n \ones\{|m_i-m_i'| \geq t\}
\leq \frac{1}{t^2} \cdot \frac{\|\bbm-\bbm'\|_2^2}{n}.\]
Concavity of $\sh(m)$ implies
$|\sh(m)-\sh(m')| \leq \sh(|m-m'|-1)-\sh(-1)=\sh(|m-m'|-1)$.
Applying the bounds $\sh(m) \leq (m+1)^{2/3}$ and $\sh(m) \in [0,\log 2]$ for
$m \in [-1,1]$,
\[|H(\bbm)-H(\bbm')| \leq \frac{1}{t^2} \cdot \frac{\|\bbm-\bbm'\|_2^2}{n}
\cdot \log 2+\left(1-\frac{1}{t^2} \cdot \frac{\|\bbm-\bbm'\|_2^2}{n}\right)
\cdot t^{2/3}.\]
Then choosing $t=(\|\bbm-\bbm'\|_2^2/n)^{3/8}$ yields
\[|H(\bbm)-H(\bbm')| \leq (\log 2+1)
\left(\frac{\|\bbm-\bbm'\|_2^2}{n}\right)^{1/4}.\]
Finally, observe that
\begin{align*}
\frac{|\bbm^\top \bY \bbm-{\bbm'}^\top \bY\bbm'|}{n}
&\leq \frac{|\bbm^\top \bY (\bbm-\bbm')|}{n}
+\frac{|{\bbm'}^\top \bY (\bbm-\bbm')|}{n}\\
&\leq 2\|\bY\|_\op \cdot \frac{\|\bbm-\bbm'\|_2}{\sqrt{n}},\\
\Big|(1-Q(\bbm))^2-(1-Q(\bbm'))^2\Big|
&\leq 2\big|Q(\bbm)-Q(\bbm')\big|
\leq 2\sqrt{2} \cdot \frac{\|\bbm-\bbm'\|_2}{\sqrt{n}}.
\end{align*}
Combining these bounds and applying to (\ref{eq:TAP}) yields the stated bound
for $\cF_\TAP$.
\end{proof}

\subsection{Properties of $q_\star$ and $\mu_\star$}

For $\lambda>1$, recall $q_\star,h_\star,e_\star$ from
(\ref{eq:fixedpointq}--\ref{eq:estar}) and the distribution $\mu_\star$ from
(\ref{eq:mustar}), and define in addition
\begin{equation}\label{eq:bstar}
b_\star=\E_{G \sim \cN(0,1)}\big[\tanh(\lambda^2 q_\star + \lambda\sqrt{q_\star} G)^4 \big]
\end{equation}

\begin{proposition}\label{prop:qstar}
For any $\lambda>1$, there is a unique solution
$q_\star \in (0,1)$ to (\ref{eq:fixedpointq}).  This solution
$q_\star = q_\star(\lambda)$ is strictly increasing in $\lambda>1$, and
satisfies $\lim_{\lambda \to 1+} q_\star(\lambda) = 0$ and $\lim_{\lambda \to 1+} q_\star(\lambda) / (\lambda - 1) = 2$. 
Furthermore, we have
\[q_\star=\E_{m \sim \mu_\star}[m^2]=\E_{m \sim \mu_\star}[m], 
\qquad b_\star=\E_{m \sim \mu_\star}[m^4]=\E_{m \sim \mu_\star}[m^3],\]
\[\lambda^2 q_\star=\E_{m \sim \mu_\star}[m \arctanh m],
\qquad h_\star=\E_{m \sim \mu_\star}[\sh(m)].\]
Finally, we have $q_\star(\lambda) >1-1/\lambda^2$. 
\end{proposition}
\begin{proof}
\cite[Appendix B.2]{deshpande2016asymptotic} shows that the function
$f(\gamma)=\E[\tanh(\gamma+\sqrt{\gamma}G)^2]$ is  strictly
increasing and strictly concave over
$\gamma \in [0,\infty)$.  By simple calculus, we have $f(0)=0$, $f'(0)=1$, and $f''(0) = -2$. Hence, for any $\lambda>1$, there is a unique solution $\gamma_\star \in (0,\lambda^2)$ to
$f(\gamma)=\gamma/\lambda^2$, and we identify $q_\star=\gamma_\star/\lambda^2$.  Furthermore, the same argument shows that $q_\star(\lambda)$ is strictly increasing in $\lambda$, and $\lim_{\lambda \to 1+} q_\star(\lambda) = 0$. Finally, simple calculus shows that $\lim_{\lambda \to 1+} q_\star(\lambda) / (\lambda - 1) = 2$. 

The identities $q_\star=\E[m^2]$ and $b_\star=\E[m^4]$ follow by definition.
\cite[Appendix B.2]{deshpande2016asymptotic} also
shows $\E[\tanh(\gamma+\sqrt{\gamma}G)^{2k}]=\E[\tanh(\gamma+\sqrt{\gamma}G)^{2k-1}]$ for any
integer $k \geq 1$, so $q_\star=\E[m]$ and $b_\star=\E[m^3]$.
The identity $\lambda^2 q_\star=\E[m \arctanh m]$ 
follows from combining $q_\star=\E[m]$ and Gaussian integration by parts,
\[\E[\lambda \sqrt{q_\star} G
\cdot \tanh(\lambda^2q_\star+\lambda\sqrt{q_\star}G)]
=\lambda^2q_\star\E[1-\tanh(\lambda^2q_\star+\lambda\sqrt{q_\star}G)^2]
=\lambda^2q_\star(1-q_\star).\]
The identity $h_\star=\E[\sh(m)]$ follows from this,
$\sh(\tanh(x))=\log 2 \cosh(x)-x\tanh(x)$, and the definition of $h_\star$.

To show $q_\star>1-1/\lambda^2$, note that for an observation
$Z \sim \cN(\lambda^2q_\star X,\lambda^2q_\star)$ with prior
$X \sim \Unif\{-1,+1\}$, $\tanh(Z)=\E[X|Z]$ is the
posterior mean estimate of $X$. Applying
$q_\star=\E[m^2]=\E[m]$, its Bayes risk is
$\E[(\E[X|Z]-X)^2]=q_\star-2q_\star+1=1-q_\star$, which may be compared to
$\E[((1+\lambda^2q_\star)^{-1}Z-X)^2]=(1+\lambda^2q_\star)^{-1}$
for the linear estimator $(1+\lambda^2 q_\star)^{-1}Z$. Since
$q_\star>0$, this linear estimator is not almost
surely equal to the Bayes estimator,
so $1-q_\star<(1+\lambda^2q_\star)^{-1}$ strictly. Rearranging yields
$q_\star>1-1/\lambda^2$.
\end{proof}

\begin{proposition}\label{prop:qstarlargelambda}
For an absolute constant $\lambda_0>0$ and all $\lambda>\lambda_0$, we have
$q_\star>1-e^{-\lambda^2/8}$.
\end{proposition}
\begin{proof}
From the identity $q_\star=\E_{m \sim \mu_\star}[m]$, monotonicity of
$\tanh$, and the bound $\tanh(x) \geq -1$, we have
\begin{align*}
q_\star&=\E_{G \sim \cN(0,1)}[\tanh(\lambda^2q_\star+\lambda \sqrt{q_\star}G)]\\
&\geq \tanh(\lambda^2 q_\star/4) \cdot \P[G \geq -3\lambda \sqrt{q_\star}/4]
-\P[G<-3\lambda \sqrt{q_\star}/4].
\end{align*}
Since $q_\star(\lambda) \to 1$ as $\lambda \to \infty$, for sufficiently large
$\lambda$ we have $q_\star>1/2$. Then, applying $\P[G<-t]<e^{-t^2/2}$
and $\tanh(t)>1-2e^{-2t}$ for $t>0$,
\[q_\star>\tanh(\lambda^2/8)(1-e^{-9\lambda^2/64})
-e^{-9\lambda^2/64}>1-e^{-\lambda^2/8}.\]
\end{proof}

\subsection{Properties of Wasserstein-2 distance}\label{appendix:wass2}

Let $W_2(\mu,\mu')$ denote the Wasserstein-2 distance between distributions
$\mu$ and $\mu'$ on $\R$, i.e.\
\[W_2(\mu,\mu')=\left(\inf_{\text{couplings } \nu \text{ of } (\mu,\mu')}
\int (x-x')^2 \de \nu(x,x')\right)^{1/2}.\]
Note that if $\mu,\mu'$ are the empirical distributions of coordinates of
$\bx,\bx' \in \R^n$, then this implies
\[W_2(\mu,\mu')^2 \leq \frac{1}{n}\|\bx-\bx'\|_2^2.\]
This distance induces the weak convergence
$W_2(\mu_n,\mu) \to 0$ if and only if
$\E_{x \sim \mu_n}[U(x)] \to \E_{x \sim \mu}[U(x)]$ for all continuous functions
$U:\R \to \R$ satisfying $\sup_{x \in \R} |U(x)|/(1+x^2)<\infty$,
see \cite[Definition 6.7]{villani2008optimal}.

In Section \ref{sec:local}, we denoted
$W(\mu,\mu')=W_2(\arctanh \mu,\arctanh \mu')$. Then,
defining the function class
\begin{equation}\label{eq:W2class}
\mathcal{Q}=\bigg\{\text{ continuous functions } U:(-1,1) \to \reals \text{ s.t. }
\sup_{m \in (-1,1)}\frac{|U(m)|}{1+(\arctanh m)^2}<\infty\bigg\}
\end{equation}
this distance $W(\mu,\mu')$ induces the weak convergence
\[\lim_{n \to \infty} W(\mu_n,\mu)=0 \text{ if and only if }
\lim_{n \to \infty}
\E_{m \sim \mu_n}[U(m)]=\E_{m \sim \mu}[U(m)] \text{ for all } U \in
\mathcal{Q}.\]

\begin{proposition}\label{prop:wasserstein}
Let $X_1,\ldots,X_n \overset{iid}{\sim} \cN(0,1)$ and let $\hat{\mu}$ be
the empirical distribution of $X_1,\ldots,X_n$.
There is a universal constant $C>0$ such that for any $t>0$,
\[\P[W_2(\hat{\mu},\,\cN(0,1)) \geq t+Cn^{-1/2}] \leq e^{-nt^2/2}.\]
\end{proposition}

\begin{proof}
By \cite[Theorem 1]{fournier2015rate}, $\E[W_2(\hat{\mu},\cN(0,1))] \leq
Cn^{-1/2}$ for a universal constant $C>0$. By the Wasserstein-2
triangle inequality, if $\mu,\mu'$ are the empirical distributions of
$\bx,\bx' \in \reals^n$, then
\[\big|W_2(\mu,\cN(0,1))-W_2(\mu',\cN(0,1))\big|
\leq W_2(\mu,\mu') \leq \Big(n^{-1}\|\bx-\bx'\|_2^2\Big)^{1/2}.\]
So $W_2(\hat{\mu},\cN(0,1))$ is $n^{-1/2}$-Lipschitz in $(X_1,\ldots,X_n)$,
and the result follows by concentration of Gaussian measure.
\end{proof}

\begin{proposition}\label{prop:wassersteinbound}
Let $\mu,\mu'$ be two probability distributions on $\reals$,
let $X \sim \mu$ and $X' \sim \mu'$, and denote $\|\mu\|_{L_2}=(\E X^2)^{1/2}$
and $\|\mu'\|_{L_2}=(\E {X'}^2)^{1/2}$.
\begin{enumerate}
\item[(a)] For any $\alpha \in (0,1)$, suppose $q_\alpha,q_\alpha'$ satisfy
$\P[|X| \geq q_\alpha]=\P[|X'| \geq q_\alpha']=\alpha$. Then
\[\big|\E[X^2\indic{|X| \geq q_\alpha}]-\E[{X'}^2
\indic{|X'| \geq q_\alpha'}]\big| \leq W_2(\mu,\mu') \cdot
(\|\mu\|_{L_2}+\|\mu'\|_{L_2}).\]
\item[(b)] Let $f$ be any function such that
$|f(x)-f(x')| \leq L(1+|x|+|x'|)|x-x'|$ for all $x,x' \in \reals$
and some constant $L>0$. Then
\[\big|\E[f(X)]-\E[f(X')]\big| \leq L \cdot W_2(\mu,\mu')
\cdot \big(1+\|\mu\|_{L_2}+\|\mu'\|_{L_2}\big).\]
\end{enumerate}
\end{proposition}
\begin{proof}
For (a), let $W_2(-|X|,-|X'|)$ denote the Wasserstein-2 distance between the
laws of $-|X|$ and $-|X'|$. Any coupling of $(X,X')$ defines also a coupling of
$(-|X|,-|X'|)$, so $W_2(-|X|,-|X'|)^2 \leq
\E[(-|X|+|X'|)^2] \leq \E[(X-X')^2]$. Taking the infimum over all couplings
$(X,X')$ gives
\[W_2(-|X|,-|X'|) \leq W_2(\mu,\mu').\]
Now consider the quantile function $G(u)$ of $-|X|$, satisfying $G(u) \leq -x$
if and only if $u \leq \P[-|X| \leq -x]=\P[|X| \geq x]$. Let $G'(u)$ be the
quantile function of $-|X'|$, let $U \sim \Unif([0,1])$, and consider the
coupling of $(-|X|,-|X'|)$ given by $-|X|=G(U)$ and $-|X'|=G'(U)$. This is the
optimal coupling that yields
\[W_2(-|X|,-|X'|)=\E\big[(-|X|+|X'|)^2\big]^{1/2}
=\E\big[(G(U)-G'(U))^2\big]^{1/2}.\]
We have $|X| \geq q_\alpha$ if and only if $G(U) \leq -q_\alpha$ if and only if
$U \leq \P[|X| \geq q_\alpha]=\alpha$, and similarly $|X'| \geq q_\alpha'$
if and only if $U \leq \alpha$. Hence
\[\E[X^2\indic{|X| \geq q_\alpha}]-\E[{X'}^2
\indic{|X'| \geq q_\alpha'}]
=\E\big[(G(U)^2-G'(U)^2)\indic{U \leq \alpha}\big],\]
so by the Cauchy-Schwarz and Minkowski inequalities,
\begin{align*}
&~\big|\E[X^2\indic{|X| \geq q_\alpha}]-\E[{X'}^2
\indic{|X'| \geq q_\alpha'}]\big|\\
&\leq \E\big[\big|G(U)^2-G'(U)^2\big|\big]\\
&\leq \E\big[(G(U)-G'(U))^2\big]^{1/2}
\E\big[(G(U)+G'(U))^2\big]^{1/2}\\
&=W_2(-|X|,-|X'|) \cdot \E[(|X|+|X'|)^2]^{1/2}
\leq W_2(\mu,\mu') \cdot (\|\mu\|_{L_2}+\|\mu'\|_{L_2}).
\end{align*}

For (b), consider any coupling of $(X,X')$. Applying again Cauchy-Schwarz
and Minkowski,
\begin{align*}
\big|\E[f(X)]-\E[f(X')]\big|
&\leq \E[|f(X)-f(X')|]\\
&\leq L \cdot \E\big[(1+|X|+|X'|)^2\big]^{1/2} \cdot
\E\big[(X-X')^2\big]^{1/2}\\
&\leq L \cdot (1+\|\mu\|_{L_2}+\|\mu'\|_{L_2})
\cdot\E\big[(X-X')^2\big]^{1/2}.
\end{align*}
The left side does not depend on the coupling, so taking the infimum over
couplings yields (b).
\end{proof}

\subsection{Sudakov-Fernique bound}

\begin{proof}[Proof of Lemma \ref{lem:slepian}]
Denote
\[G(x)=\bv(x)^\top \bW \bv(x)+f(x),
\qquad g(x)=\frac{2}{\sqrt{n}}\|\bv(x)\|_2 \langle \bg,\bv(x) \rangle
+f(x).\]
For all $x \in \cX$, we have $\E[G(x)]=\E[g(x)]=f(x)$. Furthermore, for any
$\bv,\bv' \in \R^n$, we have
\[\E\Big[\langle \bv,\bW\bv \rangle \cdot \langle \bv',\bW\bv' \rangle\Big]
=\sum_{i=1}^n \E[W_{ii}^2 v_i^2{v_i'}^2]
+4\sum_{i<j} \E[W_{ij}^2 v_iv_jv_i'v_j']
=\frac{2}{n}\langle \bv,\bv' \rangle^2.\]
Then
\begin{align*}
\E\left[\Big(\langle \bv,\bW\bv \rangle
-\langle \bv',\bW\bv' \rangle\Big)^2\right]
&=\frac{2}{n}\Big(\|\bv\|_2^4+\|\bv'\|_2^4-2\langle \bv,\bv' \rangle^2\Big)\\
&\leq \frac{2}{n}\Big(\|\bv\|_2^4+\|\bv'\|_2^4
+2\|\bv\|_2^2\|\bv'\|_2^2
-4\|\bv\|_2\|\bv'\|_2\langle\bv,\bv'\rangle\Big)\\
&\leq \frac{2}{n}\Big(2 \|\bv\|_2^4+ 2 \|\bv'\|_2^4
-4\|\bv\|_2\|\bv'\|_2\langle \bv,\bv' \rangle\Big)\\
&=\E\left[\left(\frac{2}{\sqrt{n}} \|\bv\|_2\langle \bg,\bv \rangle
-\frac{2}{\sqrt{n}}\|\bv'\|_2\langle \bg,\bv' \rangle\right)^2\right].
\end{align*}
So $\E[(G(x)-G(x'))^2] \leq \E[(g(x)-g(x'))^2]$, and the result follows from
the Sudakov-Fernique inequality, see
e.g.\ \cite[Theorem 2.2.3]{adlertaylor}.
\end{proof}

\subsection{Kac-Rice upper bounds}

We prove the Kac-Rice upper bounds of Lemmas \ref{lem:KacRicestandard}
and \ref{lem:KacRicebound} by small extensions of arguments in
\cite{fan2021tap}.

\begin{proof}[Proof of Lemma \ref{lem:KacRicestandard}]
This follows from taking $\delta \to 0$ on both sides of
\cite[Lemma A.1]{fan2021tap}, using the monotone convergence theorem.
\end{proof}

\begin{lemma}
In the setting of Lemma \ref{lem:KacRicebound}, suppose further that $T \subset
(-1,1)^n \setminus \{\bzero\}$ is compact, and its boundary $\partial T$ has 
zero Lebesgue measure in $\reals^n$. Then with probability 1 over $\bW$ and $U$,
there are no points $\bbm \in T$ satisfying
$\bg(\bbm)=\bzero$ together with any of the following three conditions:
\begin{itemize}
\item $\det \bH(\bbm)=0$, or
\item $\bbm \in \partial T$, or
\item $\ell(\bbm,\bW)+U=t$.
\end{itemize}
On this event of probability 1,
\begin{equation}\label{eq:deterministicKacRice}
|\cC|=\lim_{r \to 0} \int_T \ones\Big\{\ell(\bbm)+U<t\Big\}
\cdot \frac{1}{\Vol(B_r(\bzero))}
\ones\{\bg(\bbm) \in B_r(\bzero)\} \cdot |\det \bH(\bbm)|\,\de \bbm.
\end{equation}
\end{lemma}
\begin{proof}
It is shown in \cite[Lemma A.3]{fan2021tap} that on an event $\cE$ of
probability 1 defined by $\bW$,
no point $\bbm \in T$ satisfies both $\bg(\bbm)=\bzero$ and
either $\det \bH(\bbm)=\bzero$ or $\bbm \in \partial T$. Letting $\mathring{T}$
be the interior of $T$, it remains to
check that with probability 1, also no point $\bbm \in \mathring{T}$
satisfies both $\bg(\bbm)=\bzero$ and $\ell(\bbm,\bW)+U=t$.
Conditioning on $\bW$
and on this event $\cE$, the TAP free energy $\cF_{\TAP}(\bbm)$ is a Morse
function over $\mathring{T}$. Thus the critical points $\bbm \in
\mathring{T}$ where $\bg(\bbm)=\bzero$ are isolated, and there are at most
countably many such points. So
\[\{t-\ell(\bbm,\bW):\,\bbm \in \mathring{T},\,\bg(\bbm)=\bzero\}\]
is a countable set of values.
Since $U$ is independent of $\bW$, this implies that $U$
does not belong to this set with probability 1 conditional on
$\bW$ and $\cE$. Hence also unconditionally with probability 1 over $\bW$ and
$U$, no point $\bbm \in \mathring{T}$ satisfies both $\bg(\bbm)=\bzero$ and
$\ell(\bbm,\bW)+U=t$, as desired.

The statement (\ref{eq:deterministicKacRice}) then follows from \cite[Theorem
11.2.3]{adlertaylor}.
\end{proof}

\begin{proof}[Proof of Lemma \ref{lem:KacRicebound}]
We write as shorthand $\ell(\bbm)$ for $\ell(\bbm,\bW)$.
Suppose first that $T$ is compact with boundary $\partial T$ having Lebesgue
measure 0. Then applying (\ref{eq:deterministicKacRice}), Fatou's lemma, and
Fubini's theorem,
\begin{align*}
&~\E[|\cC|] \\
&\leq \liminf_{r \to 0} \int_T \E\bigg[\ones\{\ell(\bbm)+U<t\}
\cdot \frac{1}{\Vol(B_r(\bzero))}\ones\{\bg(\bbm) \in B_r(\bzero)\}
\cdot |\det \bH(\bbm)|\bigg]\;\de \bbm\\
&=\liminf_{r \to 0} \int_T
\bigg(\frac{1}{\Vol(B_r(\bzero))} \int_{B_r(\bzero)}
 \E\Big[\ones\{\ell(\bbm)+U<t\} \cdot |\det \bH(\bbm)|
\;\Big|\;\bg(\bbm)=\bh\Big]p_{\bg(\bbm)}(\bh)\,\de\bh\bigg) \;\de \bbm
\end{align*}
where $p_{\bg(\bbm)}(\bh)$ is the Lesbesgue-density of $\bg(\bbm)$ at
$\bg(\bbm)=\bh$.
It may be checked from the forms (\ref{eq:grad}) and (\ref{eq:hess}) for $\bg$
and $\bH$ that for any fixed $r_0>0$,
both $p_{\bg(\bbm)}(\bh)$ and $\E[|\det \bH(\bbm)| \mid
\bg(\bbm)=\bh]$ are continuous functions of $(\bbm,\bh) \in T \times
\overline{B_{r_0}(\bzero)}$ (where $\overline{B_{r_0}(\bzero)}$ is the closure of $B_{r_0}(\bzero)$). Then applying
\[\E\Big[\ones\{\ell(\bbm)+U<t\} \cdot |\det \bH(\bbm)|
\;\Big|\;\bg(\bbm)=\bh\Big]p_{\bg(\bbm)}(\bh)
\leq \E\Big[|\det \bH(\bbm)| \;\Big|\;\bg(\bbm)=\bh\Big]p_{\bg(\bbm)}(\bh),\]
this continuity, and the compactness of $T \times \overline{B_{r_0}(\bzero)}$,
the left side is bounded over $(\bbm,\bh) \in T \times
\overline{B_{r_0}(\bzero)}$.
Thus the above integrand
\[\frac{1}{\Vol(B_r(\bzero))} \int_{B_r(\bzero)}
 \E\Big[\ones\{\ell(\bbm)+U<t\} \cdot |\det \bH(\bbm)|
\;\Big|\;\bg(\bbm)=\bh\Big]p_{\bg(\bbm)}(\bh)\,\de\bh\]
is also bounded over $\bbm \in T$, and furthermore by continuity in $\bh$,
its limit as $r \to 0$ is
\[\E\Big[\ones\{\ell(\bbm)+U<t\} \cdot |\det \bH(\bbm)|
\;\Big|\;\bg(\bbm)=\bzero\Big]p_{\bg(\bbm)}(\bzero).\]
Then applying the bounded convergence theorem, we obtain
(\ref{eq:generalKacRice}).

This establishes the result for all compact
$T \in (-1,1)^n \setminus \{\bzero\}$ whose boundary
has zero Lebesgue measure, and in particular for all hyper-rectangles $T$ in
this domain. The result for all Borel-measurable $T$ then follows from the
same argument of outer measure as in the conclusion of the proof
of \cite[Lemma A.1]{fan2021tap}.
\end{proof}

\subsection{AMP state evolution}

The following lemma collects some implications of the state evolution for AMP
starting from a
spectral initialization, as characterized in \cite{montanari2021estimation}.

\begin{lemma}\label{lem:AMP}
Suppose $\lambda>1$. Let $\{\bh^k,\bbm^k\}_{k \geq 0}$
be the iterates of the AMP algorithm (\ref{eq:AMP}) with the initializations of
(\ref{eq:initialization}), where we take the
sign $\langle \bx,\bh^0 \rangle \geq 0$.
Set $\gamma_0=\lambda^2-1$ and define recursively
$\gamma_{k+1}=\lambda^2\,\E_{G \sim
\cN(0,1)}[\tanh(\gamma_k+\sqrt{\gamma_k}G)^2]$.
\begin{enumerate}
\item[(a)] For any fixed $k \geq 0$ and any function $\psi:\R^2 \to \R$
satisfying
$|\psi(x,y)-\psi(x',y')| \leq C(1+\|(x,y)\|_2+\|(x',y')\|_2)\|(x,y)-(x',y')\|_2$
for a constant $C>0$, almost surely
\[\lim_{n \to \infty} \frac{1}{n}\sum_{i=1}^n \psi(x_i,h_i^k)
=\E[\psi(X,\gamma_k X+\sqrt{\gamma_k} G)] \text{ where } X \sim \Unif\{-1,+1\}
\perp\!\!\!\perp G \sim \cN(0,1).\]
\item[(b)] We have $\lim_{k \to \infty} \gamma_k=\lambda^2 q_\star$, and
almost surely
\[\begin{aligned}&~\lim_{k \to \infty} \lim_{n \to \infty} \left(\frac{\langle \bx,\bbm^k
\rangle}{n},Q(\bbm^k),H(\bbm^k),\frac{\langle \bbm^k,\bh^k \rangle}{n},
\frac{\|\bh^k\|_2^2}{n}\right)\\
=&~(q_\star,q_\star,h_\star,\lambda^2 q_\star,\E_{m \sim \mu_\star}[\arctanh(m)^2]).\end{aligned}\]
\item[(c)] For any fixed $\eps>0$,
\[\lim_{k \to \infty} \lim_{n \to \infty}
\P\left[\frac{1}{n^2}\left\|\widehat{\bX}_{\Bayes}-\bbm^k(\bbm^k)^\top\right\|_{\sF}^2>\eps\right]=0.\]
\item[(d)] For any fixed $\eps>0$,
\[\lim_{k \to \infty} \lim_{n \to \infty}
\P\Big[|\cF_\TAP(\bbm^k)-e_\star|>\eps\Big]=0.\]
\end{enumerate}
\end{lemma}
\begin{proof}
Part (a) follows from
\cite[Theorem 2]{montanari2021estimation}, specializing to the prior
distribution $X \sim \Unif\{-1,+1\}$ and the optimal nonlinearity
$f_k(h)=\lambda \tanh(h)$ in each iteration. (The required initialization
$\bbm^{-1}=\lambda \bh^0$ was not specified in \cite{montanari2021estimation};
this condition may be derived from the observation that the principal
eigenvector $\bh$ of $\bY$ is a fixed point of the linear AMP iterations
\[\bbm^k=\alpha \bh^k, \qquad \bh^{k+1}=\bY\bbm^k-\alpha \bbm^{k-1}
=\alpha \bY\bh^k-\alpha^2 \bh^{k-1}\]
when $1=\alpha \lambda_{\max}(\bY)-\alpha^2$. In the limit $n \to \infty$,
we have $\lambda_{\max}(\bY) \to \lambda+1/\lambda$, so this is satisfied by
$\alpha=\lambda$.)

For part (b), recall from the proof of Proposition \ref{prop:qstar} that
$f(\gamma)=\E[\tanh(\gamma+\sqrt{\gamma}G)^2]$ is increasing and concave over
$\gamma \in [0,\infty)$, with $f(0)=0$ and $f'(0)=1$. Then the iterations
$\gamma_{k+1}=\lambda^2 f(\gamma_k)$ must converge to the unique fixed point
$\gamma_\star=\lambda^2 q_\star$ from any positive initialization $\gamma_0>0$.
Applying part (a) with
\[\psi(x,y) \in \{x\tanh(y),
\tanh(y)^2,\sh(\tanh(y)),y\tanh(y),y^2\}\]
and using Proposition \ref{prop:qstar} to evaluate the
Gaussian expectations in the limit $k \to \infty$, part (b) follows.

For part (c), we have by part (b) almost surely
\begin{align*}
\lim_{k \to \infty} \lim_{n \to \infty}
\frac{1}{n^2}\left\|\bx\bx^\top-\bbm^k(\bbm^k)^\top\right\|_{\sF}^2
&=\lim_{k \to \infty} \lim_{n \to \infty}
\left(1-2\left(\frac{\langle \bx,\bbm^k \rangle}{n}\right)^2
+Q(\bbm^k)^2\right)\\
&=1-q_\star^2.
\end{align*}
This coincides with the Bayes risk $\lim_{n \to \infty} \E[
\|\widehat{\bX}_{\Bayes}-\bx\bx^\top\|_{\sF}^2]/n^2=1-q_\star^2$,
see e.g.\ \cite[Section 2.3]{lelarge2019fundamental}. 
Since $\widehat{\bX}_{\Bayes}=\E[\bx\bx^\sT \mid \bY]$ and $\bbm^k(\bbm^k)^\top$
is a function of $\bY$, we have the Pythagorean relation
\[\E\Big[\|\bx\bx^\top-\bbm^k(\bbm^k)^\sT\|_\sF^2\Big]
=\E\Big[\|\widehat{\bX}_{\Bayes}-\bbm^k(\bbm^k)^\sT\|_\sF^2\Big]
+\E\Big[\|\widehat{\bX}_{\Bayes}-\bx\bx^\sT\|_\sF^2\Big],\]
hence $\lim_{k \to \infty} \lim_{n \to \infty} 
\E[\|\widehat{\bX}_{\Bayes}-\bbm^k(\bbm^k)^\sT\|_{\sF}^2]/n^2=0$.
Part (c) then follows by Markov's inequality.

Finally, for part (d), let us use the notational shorthand $X \to c$ to mean
the convergence in probability
$\lim_{k \to \infty} \lim_{n \to \infty} \P[|X-c|>\eps]=0$ for any fixed
$\eps>0$. Observe that (c) implies
\[\frac{1}{n^2} \Big\|\bbm^k(\bbm^k)^\top-\bbm^{k+1}
(\bbm^{k+1})^\top\Big\|_\sF^2 \to 0.\]
Writing
\[\frac{1}{n^2} \Big\|\bbm^k(\bbm^k)^\top-\bbm^{k+1}
(\bbm^{k+1})^\top\Big\|_\sF^2
=Q(\bbm^k)^2+Q(\bbm^{k+1})^2-2\left(\frac{\langle \bbm^k,\bbm^{k+1} \rangle}{n}
\right)^2,\]
this and the statements $Q(\bbm^k),Q(\bbm^{k+1}) \to q_\star$ from (b) imply
$|\langle \bbm^k,\bbm^{k+1} \rangle|/n \to q_\star$.
Applying again $Q(\bbm^k),Q(\bbm^{k+1}) \to q_\star$ and Cauchy-Schwarz, this
in turn implies $\|\bbm^k \pm \bbm^{k+1}\|_2^2/n \to 0$ for some choice of sign
$\pm$. Part (b) shows $\langle \bx,\bbm^k \rangle/n \to q_\star$ and
$\langle \bx,\bbm^{k+1} \rangle/n \to q_\star$ for the same positive
constant $q_\star$, so
$\|\bbm^k-\bbm^{k+1}\|_2^2/n \to 0$ must hold with the sign $-$.
Then by part (b) and Cauchy-Schwarz, we have also
\[\langle \bbm^k,\bbm^{k+1} \rangle/n \to q_\star, 
\qquad \langle \bbm^k,\bh^{k+1} \rangle/n \to \lambda^2 q_\star.\]
Writing $\lambda \bY\bbm^k=\bh^{k+1}+\lambda^2(1-Q(\bbm^k))\bbm^{k-1}$, this
yields
\begin{align*}
\cF_\TAP(\bbm^k)&=-\frac{1}{2n}\Big[(\bbm^k)^\sT \bh^{k+1}
+\lambda^2(1-Q(\bbm^k))(\bbm^k)^\sT \bbm^{k-1}\Big]
-H(\bbm^k)-\frac{\lambda^2}{4}\Big(1-Q(\bbm)\Big)^2\\
&\to -\frac{1}{2}(\lambda^2 q_\star+\lambda^2(1-q_\star)q_\star)
-h_\star-\frac{\lambda^2}{4}(1-q_\star)^2=e_\star.
\end{align*}
\end{proof}

We remark that this proves Lemma \ref{lem:AMP} marginally over $\bx \sim
\Unif(\{-1,+1\}^n)$. The claims of parts (b--d) then hold also
conditional on any fixed $\bx \in \{-1,+1\}^n$, and in particular conditional on
$\bx=\ones$, by sign symmetry.

\section{Proofs for the local analysis of TAP and AMP}\label{appendix:local}

\subsection{Proofs for Section \ref{sec:TAP-lb}}\label{appendix:TAP-lb-proofs}

\begin{proof}[Proof of Lemma \ref{lem:global_min}]
Fixing $\bx=\ones$, we have
\[-\cF_{\TAP}(\bbm)=\frac{\lambda}{2n}\langle \bbm,\bW\bbm \rangle
+\frac{\lambda^2}{2}M(\bbm)^2
+\frac{\lambda^2}{4}\Big(1-Q(\bbm)\Big)^2+H(\bbm).\]
Let $\bg \sim \cN(\bzero, \id_n)$, and define an auxiliary Gaussian process
$G(\bbm)$ to be
\[
\begin{aligned}
G(\bbm) =&~ \frac{\lambda}{n^{3/2}} \| \bbm \|_2 \< \bg, \bbm\>  +
\frac{\lambda^2}{2}M(\bbm)^2
+ \frac{\lambda^2}{4} \Big(1 - Q(\bbm) \Big)^2 + H(\bbm).
\end{aligned}
\]
Then by Lemma \ref{lem:slepian},
for any domain $\Omega \subset (-1,1)^n$ we have
\begin{equation}\label{eq:slepianFvalue}
\E\Big[\sup_{\bbm \in \Omega} {-}\cF_{\TAP}(\bbm)\Big]
\leq \E\Big[\sup_{\bbm \in \Omega} G(\bbm)\Big].
\end{equation}
Introduce
\[\Gamma(g;p,\gamma,\tau,\nu)=\sup_{m \in (-1,1)} \lambda p \cdot gm
+\frac{\gamma m^2}{2}+\tau m+\nu \sh(m).\]
Then
\begin{equation}\label{eq:scalarGm2}
\begin{aligned}
&~\sup_{\bbm \in (-1,1)^n:\,(Q(\bbm), M(\bbm), H(\bbm)) \in K} G(\bbm) \\
=&~ \sup_{(q, \vphi, h) \in K}\; \sup_{\bbm \in (-1,1)^n:\,(Q(\bbm), M(\bbm),
H(\bbm)) = (q, \vphi, h)}
\frac{\lambda}{ n}  \< \bg, \bbm\> \sqrt{q} + \frac{\lambda^2}{2}
\vphi^2 + \frac{\lambda^2}{4}(1 - q)^2+h \\
\le&~ \sup_{(q, \vphi, h) \in K} \inf_{(\gamma, \tau, \nu) \in K'}
\sup_{\bbm \in (-1,1)^n}
\frac{\lambda}{n}  \< \bg, \bbm\> \sqrt{q} + \frac{\lambda^2}{2} \vphi^2 
+ \frac{\lambda^2}{4}(1 - q)^2 + h  \\
&\hspace{1.5in}+ \tau\Big(M(\bbm)-\varphi\Big)
+ \frac{\gamma}{2} \Big(Q(\bbm)-q\Big) +\nu\Big(H(\bbm)-h\Big)  \\
=&~ \sup_{(q, \vphi, h) \in K} \inf_{(\gamma, \tau, \nu) \in K'}
\frac{\lambda^2}{2} \vphi^2 + \frac{\lambda^2}{4}(1 - q)^2 + h
- \frac{q \gamma}{2}-\varphi \tau - h \nu
+\frac{1}{n}\sum_{i=1}^n \Gamma(g_i;\sqrt{q},\gamma,\tau,\nu) \\
\leq&~ \sup_{(q,\varphi,h) \in K} \inf_{(\gamma,\tau,\nu) \in K'}
-E_{\lambda}(q,\varphi,h;\gamma,\tau,\nu)\\
&~+\mathop{\sup_{(q,\varphi,h) \in K}}_{(\gamma,\tau,\nu) \in K'}
\left|\frac{1}{n}\sum_{i=1}^n \Gamma(g_i;\sqrt{q},\gamma,\tau,\nu)
-\E_{G \sim \cN(0,1)}[\Gamma(G;\sqrt{q},\gamma,\tau,\nu)]\right|.
\end{aligned}
\end{equation}

Let $\omega=(p,\gamma,\tau,\nu)$, let $K_\omega=[0,1] \times K'$, and let
$C,c>0$ be $(\lambda,K',\eps)$-dependent constants changing from instance
to instance. For any fixed $\omega \in K_\omega$,
the function $\Gamma(g;\omega)$ is $\lambda p$-Lipschitz
in $g$. Then the Gaussian concentration and Hoeffding
inequalities show $\|\Gamma(G;\omega)-\E[\Gamma(G;\omega)]\|_{\psi_2} \leq C$
and hence
\begin{equation}\label{eq:Gammaconcentration}
\E\left[\left\vert \frac{1}{n}\sum_{i = 1}^n
\Gamma(g_i;\omega)-\E[\Gamma(G;\omega)]\right\vert\right] \leq C/\sqrt{n}.
\end{equation}
To obtain uniform control over $K_\omega$, let $\cN$ be a $n^{-1/16}$-net
of $K_\omega$ of cardinality $|\cN| \leq Cn^{1/4}$. Fixing $g$ and applying
$\sh(m) \in [0,\log 2]$, observe that $\Gamma(g;\omega)$ is
$\frac{1}{2}$-Lipschitz in $\gamma$, 1-Lipschitz in $\tau$,
$(\log 2)$-Lipschitz in $\nu$, and $\lambda|g|$-Lipschitz in $p$.
Then $n^{-1}\sum_{i=1}^n \Gamma(g_i;\omega)$ is
$C(1+n^{-1}\sum_{i=1}^n |g_i|)$-Lipschitz in $\omega$. Applying
(\ref{eq:Gammaconcentration}) for $\omega \in \cN$, we then
obtain for all large $n$ that
\begin{align}
&\E\left[\sup_{\omega \in K_\omega}
\left\vert \frac{1}{n}\sum_{i = 1}^n \Gamma(g_i;\omega)
-\E[\Gamma(G;\omega)]\right\vert\right]\nonumber\\
&\leq \E\left[\sup_{\omega \in \cN}
\left\vert \frac{1}{n}\sum_{i = 1}^n \Gamma(g_i;\omega)
-\E[\Gamma(G;\omega)]\right\vert\right]
+Cn^{-1/16} \cdot
\E\left[1+\frac{1}{n}\sum_{i=1}^n |g_i|\right]\nonumber\\
&\leq |\cN| \cdot \frac{C}{\sqrt{n}}+Cn^{-1/16} \cdot
\E\left[1+\frac{1}{n}\sum_{i=1}^n |g_i|\right]<\eps/2.
\label{eq:uniformconcGamma}
\end{align}

Combining (\ref{eq:slepianFvalue}), (\ref{eq:scalarGm2}), and
(\ref{eq:uniformconcGamma}) and negating the sign, we arrive at
\begin{align*}
\E\left[\inf_{\bbm \in (-1,1)^n:\,(Q(\bbm),M(\bbm),H(\bbm)) \in K}
\cF_{\TAP}(\bbm)\right] \geq
\inf_{(q,\varphi,h) \in K} \sup_{(\gamma,\tau,\nu) \in K'}
E_{\lambda}(q,\varphi,h;\gamma,\tau,\nu)-\eps/2
\end{align*}
for all large $n$. Finally, writing 
$\bW=(\bZ+\bZ^\sT)/\sqrt{2n}$ where $\bZ$ has i.i.d.\ $\cN(0,1)$ entries,
observe that $\bZ \mapsto \langle \bbm,\bW\bbm \rangle$ is
$\sqrt{2/n}\|\bbm\|_2^2$-Lipschitz with respect to the Frobenius norm of $\bZ$,
and $\|\bbm\|_2^2 \leq n$. Then 
$\inf_{\bbm \in (-1,1)^n:\,(Q(\bbm),M(\bbm),H(\bbm)) \in K}
\cF_{\TAP}(\bbm)$ is $C/\sqrt{n}$-Lipschitz in $\bZ$. Applying Gaussian
concentration of measure,
\[
\begin{aligned}
\P\Big[&~\inf_{\bbm \in (-1,1)^n:\,(Q(\bbm),M(\bbm),H(\bbm)) \in K}
\cF_{\TAP}(\bbm) \\
&~~~~~~~~\leq \E\Big(\inf_{\bbm \in (-1,1)^n:\,(Q(\bbm),M(\bbm),H(\bbm)) \in K}
\cF_{\TAP}(\bbm)\Big)-\eps/2 \Big]<e^{-cn}.
\end{aligned}\]
Combining the above two displays shows (\ref{eq:valuelowerbound}).
\end{proof}

\begin{proof}[Proof of Lemma \ref{lem:global_min_analysis}]
Observe that $\sh'(m)={-}\arctanh m$ and $\sh''(m)=-1/(1-m^2) \leq -1$ for all
$m \in (-1,1)$. Thus, for $\nu>\max(\gamma,0)$
and any realization of $G$, the function
$m \mapsto \lambda \sqrt{q}\cdot Gm+\gamma m^2/2+\tau m+\nu\sh(m)$ is strictly
concave over $m \in (-1,1)$, with derivative diverging to $\mp \infty$ as
$m \to \pm 1$. Then its supremum is achieved at the unique value
$m=m(G;q,\gamma,\tau,\nu) \in (-1,1)$ satisfying the stationary
condition
\begin{equation}\label{eq:mstarstationary}
0=\lambda \sqrt{q}G+\tau+\gamma m-\nu\arctanh(m).
\end{equation}
The above strict concavity implies that the derivative in $m$
of the right side of (\ref{eq:mstarstationary}) is non-zero, so
by the implicit function theorem, $m(G;q,\gamma,\tau,\nu)$ is
analytic over $\{(G,q,\gamma,\tau,\nu) \in \R^5:q>0,\,\nu>\max(\gamma,0)\}$.

Denote
\[
\begin{aligned}
&~F(G;q,\gamma,\tau,\nu)\\
=&~\lambda \sqrt{q}\cdot Gm(G;q,\gamma,\tau,\nu)
+\frac{\gamma m(G;q,\gamma,\tau,\nu)^2}{2}+\tau m(G;q,\gamma,\tau,\nu)
+\nu \sh(m(G;q,\gamma,\tau,\nu)).\end{aligned}\]
Applying that (\ref{eq:mstarstationary}) holds at $m=m(G;q,\gamma,\tau,\nu)$ to
cancel the terms involving the derivatives of $F$ in $m$, we obtain
\begin{equation}\label{eq:firstderF}
\partial_q F=\frac{\lambda G m}{2\sqrt{q}},
\quad \partial_\gamma F=\frac{m^2}{2}, \quad \partial_\tau F=m,
\quad \partial_\nu F=\sh(m).
\end{equation}
Differentiating (\ref{eq:mstarstationary}) implicitly, we have also
\[\nabla_{q,\gamma,\tau,\nu} m=\left(\frac{\nu}{1-m^2}-\gamma\right)^{-1}
\left(\frac{\lambda G}{2\sqrt{q}},\,m,\,1,\,{-}\arctanh m\right).\]
Then, differentiating (\ref{eq:firstderF}) a second time and applying these
forms,
\begin{equation}\label{eq:secondderF}
\begin{aligned}
\nabla_{q,\gamma,\tau,\nu} \partial_q F&=-\frac{\lambda Gm}{4q^{3/2}}
+\frac{\lambda G}{2\sqrt{q}}\left(\frac{\nu}{1-m^2}-\gamma\right)^{-1}
\left(\frac{\lambda G}{2\sqrt{q}},\,m,\,1,\,{-}\arctanh m\right),\\
\nabla_{q,\gamma,\tau,\nu} \partial_\gamma F
&=m\left(\frac{\nu}{1-m^2}-\gamma\right)^{-1}
\left(\frac{\lambda G}{2\sqrt{q}},\,m,\,1,\,{-}\arctanh m\right),\\
\nabla_{q,\gamma,\tau,\nu} \partial_\tau F
&=\left(\frac{\nu}{1-m^2}-\gamma\right)^{-1}
\left(\frac{\lambda G}{2\sqrt{q}},\,m,\,1,\,{-}\arctanh m\right),\\
\nabla_{q,\gamma,\tau,\nu} \partial_\nu F&=
-\arctanh(m)\left(\frac{\nu}{1-m^2}-\gamma\right)^{-1}
\left(\frac{\lambda G}{2\sqrt{q}},\,m,\,1,\,{-}\arctanh m\right).
\end{aligned}
\end{equation}
Applying these expressions and the identity
$\arctanh m=(\lambda \sqrt{q} G+\tau+\gamma m)/\nu$ from
(\ref{eq:mstarstationary}), we may check that over a local neighborhood of any
$(q,\gamma,\tau,\nu)$ satisfying $q>0$ and $\nu>\max(\gamma,0)$, we have the
bounds $|\partial F|,|\partial^2 F| \leq C(G^2+1)$ for all first and second
partial derivatives of $F$ and some constant $C>0$ depending on $\lambda$ and
this neighborhood. Thus, the dominated convergence theorem may be applied to
differentiate $\E_{G \sim \cN(0,1)}[F(G;q,\gamma,\tau,\nu)]$ twice under the
integral. This implies also that $E_{\lambda}(q,\varphi,h;\gamma,\tau,\nu)$
as defined in (\ref{eq:Elambda}) is twice
continuously-differentiable in all arguments.

Fixing any $(q,\varphi,h)$, observe that $(\gamma,\tau,\nu) \mapsto
E_{\lambda}(q,\varphi,h;\gamma,\tau,\nu)$ is concave on all of $\reals^3$, by its
definition. At $(q,\varphi,h)=(q_\star,q_\star,h_\star)$ and
$(\gamma,\tau,\nu)=(0,\lambda^2 q_\star,1)$, the supremum over $m$ in
\[\sup_{m \in (-1,1)} \lambda \sqrt{q}\cdot Gm+\frac{\gamma m^2}{2}
+\tau m+\nu\sh(m)
=\sup_{m \in (-1,1)} \lambda \sqrt{q_\star}\cdot Gm+\lambda^2 q_\star m+\sh(m)\]
is achieved at $m=\tanh(\lambda^2 q_\star+\lambda\sqrt{q_\star}G)$, which
has the distribution $\mu_\star$. Then
\[E_{\lambda}(q_\star,q_\star,h_\star;0,\lambda^2q_\star,1)=e_\star\]
by Proposition \ref{prop:qstar}. Specializing (\ref{eq:firstderF})
to the point $(q_\star,q_\star,h_\star;0,\lambda^2q_\star,1)$, we have
\begin{align*}
\nabla_{\gamma,\tau,\nu} E_{\lambda}(q_\star,q_\star,h_\star;0,\lambda^2q_\star,1)
&=\Big(q_\star/2-\E_{m \sim \mu_\star}[m^2/2],\,q_\star
-\E_{m \sim \mu_\star}[m],\, h_\star-\E_{m \sim \mu_\star}[\sh(m)]\Big)\\
&=0.
\end{align*}
So $(\gamma,\tau,\nu)=(0,\lambda^2 q_\star,1)$ is a critical point of
$(\gamma,\tau,\nu) \mapsto E_{\lambda}(q_\star,q_\star,h_\star;\gamma,\tau,\nu)$. Then
by concavity of this function, $(\gamma,\tau,\nu)=(0,\lambda^2 q_\star,1)$ is
a global maximizer, and we obtain (\ref{eq:global_min_point}).

To show the lower bound (\ref{eq:global_min_bound}), we first check strict
concavity around $(\gamma,\tau,\nu)=(0,\lambda^2 q_\star,1)$: 
Specializing (\ref{eq:secondderF}) to the point
$(q,\varphi,h;\gamma,\tau,\nu)=(q_\star,q_\star,h_\star;0,\lambda^2 q_\star,1)$,
we have
\begin{equation}\label{eq:dualderivative}
\nabla_{\gamma,\tau,\nu}^2 E_{\lambda}(q_\star,q_\star,h_\star;0,\lambda^2 q_\star,1)
=-\E_{m \sim \mu_\star}
\left[(1-m^2)\begin{pmatrix} m \\ 1 \\ -\arctanh m
 \end{pmatrix}\begin{pmatrix} m \\ 1 \\ -\arctanh m
 \end{pmatrix}^\sT\right].
\end{equation}
The distribution of
$m \sim \mu_\star$ is supported on the full interval $(-1,1)$, and the curve
$\{(m,1,-\arctanh m):m \in (-1,1)\}$ is not
contained in any 2-dimensional subspace. Thus this Hessian is strictly
negative definite.

Applying again the implicit function theorem,
this implies that the maximizer of $(\gamma,\tau,\nu) \mapsto
(q,\varphi,h;\gamma,\tau,\nu)$ is implicitly defined as a twice
continuously-differentiable function of
$(q,\varphi,h)$, in a neighborhood of $(q_\star,q_\star,h_\star)$. Fixing any
subset $K' \subseteq \reals^3$ containing $(0,\lambda^2 q_\star,1)$ in its
interior, this maximizer must then belong to the interior of $K'$ for all
\[(q,\varphi,h) \in K=\Big\{(q,\varphi,h):|q-q_\star|,|\varphi-q_\star|,
|h-h_\star| \leq \delta\Big\}\]
and some sufficiently small $\delta$. Defining
\begin{equation}\label{eq:Ebar}
\bar{E}_\lambda(q,\varphi,h)=\sup_{(\gamma,\tau,\nu) \in K'}
E_{\lambda}(q,\varphi,h;\gamma,\tau,\nu),
\end{equation}
this function $\bar{E}_\lambda(q,\varphi,h)$ is then twice
continuously-differentiable on $K$. We proceed to show
two properties of $\bar{E}_\lambda$: For sufficiently small constants
$c,\delta>0$ and any $(q,\varphi,h) \in K$,
\begin{align}
\bar{E}_\lambda(q,\varphi,h) &\geq e_\star+c(q-q_\star)^2+c(\varphi-q_\star)^2,
\label{eq:Elower1}\\
\bar{E}_\lambda(q_\star,q_\star,h) &\geq e_\star+c(h-h_\star)^2. \label{eq:Elower2}
\end{align}

To show (\ref{eq:Elower1}), note that restricting to $\nu=1$ removes the
dependence of $E$ on $h$. We denote this restriction as
$E_{\lambda}(q,\varphi;\gamma,\tau)$. Restricting the supremum in (\ref{eq:Ebar}) to
$\nu=1$ yields the lower bound
\[\bar{E}_\lambda(q,\varphi,h) \geq 
\bar{E}_\lambda(q,\varphi)=\sup_{(\gamma,\tau):(\gamma,\tau,1) \in K'}
E_{\lambda}(q,\varphi;\gamma,\tau).\]
We compute the Hessian of $E_{\lambda}(q,\varphi;\gamma,\tau)$ at
$(q_\star,q_\star;0,\lambda^2 q_\star)$: 
Differentiating (\ref{eq:mstarstationary}) implicitly now in $q$ and $G$
at $(q,\varphi,h;\gamma,\tau,\nu)=(q_\star,q_\star,h_\star;0,\lambda^2
q_\star,1)$,
\begin{equation}\label{eq:dermGq}
\partial_G m=\lambda \sqrt{q_\star}(1-m^2),
\qquad \partial_q m=\frac{\lambda}{2\sqrt{q_\star}}G(1-m^2).
\end{equation}
Recalling that $m=\tanh(\lambda^2 q_\star+\lambda \sqrt{q_\star} G)$ at this
point, we write $\E$ for the expectation over $G \sim \cN(0,1)$.
Then, applying $b_\star=\E[m^3]=\E[m^4]$ from Proposition \ref{prop:qstar}
and $\E[Gm]=\E[\partial_G m]=\lambda\sqrt{q_\star}(1-q_\star)$
by Gaussian integration by parts,
\[\nabla^2 E_{\lambda}(q,\varphi;\gamma,\tau)
\Big|_{(q,\varphi;\gamma,\tau)
=(q_\star,q_\star;0,\lambda^2q_\star)}
=\diag(\tfrac{1}{2},1,1,1) \cdot
\begin{pmatrix} \bA_{11} & \bA_{12} \\ \bA_{21} & \bA_{22} \end{pmatrix}
\cdot \diag(\tfrac{1}{2},1,1,1)\]
where
\begin{align}
\bA_{11}&=\begin{pmatrix} -2\lambda^2
-\E\left[-\frac{\lambda G}{q_\star^{3/2}}m+\frac{2\lambda
G}{q_\star^{1/2}} \partial_q m\right] & 0 \\ 0 & -\lambda^2
\end{pmatrix}\nonumber\\
&=\begin{pmatrix} -2\lambda^2+\frac{\lambda^2(1-q_\star)}{q_\star}
-\frac{\lambda^2}{q_\star}\E[G^2(1-m^2)] & 0 \\
0 & -\lambda^2 \end{pmatrix}\nonumber\\
\bA_{12}=\bA_{21}^\sT
&=\begin{pmatrix} 1-2\E[m \partial_q m]
& -2\E[\partial_q m] & \\ 0 & 1 \end{pmatrix}
=\begin{pmatrix}
1-\frac{\lambda}{\sqrt{q_\star}}\E[Gm(1-m^2)]
 & \frac{\lambda}{\sqrt{q_\star}}\E[Gm^2] \\ 0 & 1 \end{pmatrix}
\nonumber\\
\bA_{22}
&=\begin{pmatrix} -\E[m\partial_\gamma m] &
-\E[m \partial_\tau m] \\
-\E[m \partial_\tau m] & -\E[\partial_\tau m] \end{pmatrix}
=-\begin{pmatrix} q_\star-b_\star
& q_\star-b_\star \\ q_\star-b_\star & 1-q_\star \end{pmatrix}
\label{eq:A22}
\end{align}
We may simplify the above expressions for $\bA_{11}$ and $\bA_{12}$
further using the integration by parts identities
\begin{align*}
\E[Gm^2]&=\E[2m \partial_G m]
=2\lambda \sqrt{q_\star}(q_\star-b_\star)\\
\E[Gm(1-m^2)]&=\E[(1-3m^2)\partial_G m]
=\lambda\sqrt{q_\star}(1-4q_\star+3b_\star)\\
\E[G^2(1-m^2)]&=\E[(1-m^2)-2Gm\partial_G m]\\
&=1-q_\star-2\lambda\sqrt{q_\star}\E[Gm(1-m^2)]
=1-q_\star-2\lambda^2q_\star(1-4q_\star+3b_\star),
\end{align*}
yielding
\[\bA_{11}
=\begin{pmatrix} -2\lambda^2+2\lambda^4(1-4q_\star+3b_\star)
& 0 \\ 0 & -\lambda^2 \end{pmatrix}, \quad
\bA_{12}=\begin{pmatrix} 1-\lambda^2(1-4q_\star+3b_\star) &
2\lambda^2(q_\star-b_\star) \\ 0 & 1 \end{pmatrix}.\]

Here $\bA_{22}$ is the upper-left $2 \times 2$ submatrix of
(\ref{eq:dualderivative}), which we have argued satisfies $\bA_{22} \prec 0$.
Computing explicitly its inverse,
the Hessian of $\bar{E}_\lambda(q,\varphi)$ at $(q_\star,q_\star)$ is then
given by
\[\nabla^2 \bar{E}_\lambda(q,\varphi)\Big|_{(q,\varphi)=(q_\star,q_\star)}
=\begin{pmatrix} \frac{1}{2} & 0 \\ 0 & 1 \end{pmatrix}
\Big(\bA_{11}-\bA_{12} \bA_{22}^{-1} \bA_{21}\Big)
\begin{pmatrix} \frac{1}{2} & 0 \\ 0 & 1 \end{pmatrix}
=\begin{pmatrix} \frac{1}{2} & 0 \\ 0 & 1 \end{pmatrix}
\begin{pmatrix} c_1 & -c_2 \\ -c_2 & c_2 \end{pmatrix}
\begin{pmatrix} \frac{1}{2} & 0 \\ 0 & 1 \end{pmatrix}\]
where, after some algebraic simplification, 
\[\begin{aligned}c_1=&~\frac{(1-q_\star)-2\lambda^2(1-2q_\star+b_\star)^2
+\lambda^4(1-2q_\star+b_\star)^3}
{(1-2q_\star+b_\star)(q_\star-b_\star)}
-\lambda^4(1-2q_\star+b_\star), \\
c_2=&~\frac{1}{1-2q_\star+b_\star}-\lambda^2.\end{aligned}\]
From Proposition \ref{prop:qstar}, we have the inequalities
\begin{equation}\label{eq:basicinequalities}
\begin{aligned}
q_\star-b_\star&=\E[m^2(1-m^2)]>0,\\
1-2q_\star+b_\star&=\E[(1-m^2)^2]>0,\\
1-2q_\star+b_\star&<1-q_\star<\frac{1}{\lambda^2}
\end{aligned}
\end{equation}
so $c_2>0$. We may compute the Schur-complement
\begin{align*}
c_1-(-c_2)c_2^{-1}(-c_2)&=c_1-c_2\\
&=\lambda^2\Big(1-\lambda^2(1-2q_\star+b_\star)\Big)
+\frac{1}{q_\star-b_\star}\Big(1-\lambda^2(1-2q_\star+b_\star)\Big)^2>0,
\end{align*}
where the last inequality applies $q_\star-b_\star>0$ and
$1-\lambda^2(1-2q_\star+b_\star)>0$. The statements $c_2>0$ and
$c_1-(-c_2)c_2^{-1}(-c_2)>0$ together imply
$\bA_{11}-\bA_{12} \bA_{22}^{-1} \bA_{21} \succ 0$, so 
$\bar{E}_\lambda(q,\varphi)$ is strongly convex at $(q,\varphi)=(q_\star,q_\star)$.
Then for small enough $\delta,c>0$, we obtain $\bar{E}_\lambda(q,\varphi) \geq
e_\star+c(q-q_\star)^2+c(\varphi-q_\star)^2$ for all $(q,\varphi,h) \in K$,
and hence (\ref{eq:Elower1}) holds.

To show (\ref{eq:Elower2}), let us restrict $E$ to
$(q,\varphi;\gamma,\tau)=(q_\star,q_\star;0,\lambda^2 q_\star)$ and denote this
restriction as $E_{\lambda}(h;\nu)$. Restricting the supremum in (\ref{eq:Ebar})
to $(\gamma,\tau)=(0,\lambda^2 q_\star)$ yields
\[\bar{E}_\lambda(q_\star,q_\star,h) \geq \bar{E}_\lambda(h)
=\sup_{\nu:(0,\lambda^2 q_\star,\nu) \in K'} E_{\lambda}(h;\nu).\]
The Hessian of $E_{\lambda}(h;\nu)$ is
\[\nabla^2 E_{\lambda}(h;\nu)\Big|_{(h;\nu)=(h_\star;1)}
=\begin{pmatrix} 0 & 1 \\ 1 & \E[(\arctanh m)\partial_\nu m]
\end{pmatrix}.\]
The lower-right entry is the $(3,3)$ entry of (\ref{eq:dualderivative}),
which we have argued is negative. Then
\[\nabla^2 \bar{E}_\lambda(h)\Big|_{h=h_\star}=0-1 \cdot
\E[(\arctanh m)\partial_\nu m]^{-1} \cdot 1>0,\]
so $\bar{E}_\lambda(h)$ is strongly convex at $h=h_\star$. Then
$\bar{E}_\lambda(h) \geq e_\star+c(h-h_\star)^2$ for small enough $c,\delta>0$,
implying (\ref{eq:Elower2}).

Finally, let $C_\lambda>0$ be an upper bound for $\|\nabla^2
\bar{E}_\lambda(q,\varphi,h)\|_{\op}$ over $(q,\varphi,h) \in K$. Since
$\bar{E}_\lambda(q,\varphi,h)$ is twice continuously-differentiable
on $K$ and
$\nabla \bar{E}_\lambda(q_\star,q_\star,h_\star)=0$, (\ref{eq:Elower2}) implies
\begin{align*}
\bar{E}_\lambda(q,\varphi,h) &\geq \bar{E}_\lambda(q_\star,q_\star,h)
-|\bar{E}_\lambda(q,\varphi,h)-\bar{E}_\lambda(q_\star,q_\star,h)|\\
&\geq e_\star+c(h-h_\star)^2
-\sup_{(q,\varphi,h) \in K} \|\nabla \bar{E}_\lambda(q,\varphi,h)\|
\cdot \Big[|q-q_\star|+|\varphi-q_\star|\Big]\\
&\geq e_\star+c(h-h_\star)^2
-C_\lambda\Big[|q-q_\star|+|\varphi-q_\star|+|h-h_\star|\Big]
\cdot \Big[|q-q_\star|+|\varphi-q_\star|\Big].
\end{align*}
If $|q-q_\star|+|\varphi-q_\star|<c/(4C_\lambda) \cdot |h-h_\star|$,
we apply this bound to obtain
$\bar{E}_\lambda(q,\varphi,h) \geq e_\star+(c/2)(h-h_\star)^2$. Otherwise,
we apply (\ref{eq:Elower1}), and combining these cases
yields (\ref{eq:global_min_bound}).
\end{proof}

\begin{proof}[Proof of Corollary \ref{cor:critexists}]
 Throughout the proof, for any set $K \subseteq \R^d$, we denote $\overline{K}$ as the closure of $K$. Fix any compact set $K' \subset \reals^3$ containing $(0,\lambda^2 q_\star,1)$
in its interior. Define
\[K_\delta=\Big\{(q,\varphi,h):|q-q_\star|,|\vphi-q_\star|,|h-h_\star|
<\delta\Big\}.\]
For $\delta>0$ sufficiently small, we apply Lemma \ref{lem:global_min} once
with $K=\overline{K_\delta}$ and once with
$K=\overline{K_\delta} \setminus K_{\delta/2}$.
Then Lemmas \ref{lem:global_min} and \ref{lem:global_min_analysis}
combine to show, almost surely for all large $n$ and a $\lambda$-dependent
constant $c_0>0$,
\begin{align}
\inf_{\bbm \in \overline{\cB_{\delta}}}
\cF_{\TAP}(\bbm) &> e_\star-\delta,\label{eq:FTAPlowerbound1}\\
\inf_{\bbm \in \overline{\cB_{\delta}} \setminus \cB_{\delta/2}}
\cF_{\TAP}(\bbm) &> e_\star+c_0\delta^2.
\label{eq:FTAPlowerbound2}
\end{align}

Lemma \ref{lem:AMP}(b) and (d) imply that for some sufficiently large $k$,
with probability approaching 1 as $n \to \infty$,
the AMP iterate $\bbm^k$ satisfies $\bbm^k \in \cB_{\delta/2}$ and
$\cF_{\TAP}(\bbm^k)<e_\star+(c_0/2)\delta^2$. 
Together with (\ref{eq:FTAPlowerbound2}), this
implies that there must exist a critical point and local minimizer $\bbm_\star$
of $\cF_{\TAP}$ in $\cB_{\delta/2}$, satisfying
$\cF_{\TAP}(\bbm_\star)<e_\star+(c_0/2)\delta^2<e_\star+\delta$. Applying
(\ref{eq:FTAPlowerbound1}), we get
$|\cF_{\TAP}(\bbm_\star)-e_\star|<\delta$.
\end{proof}

\subsection{Proofs for Section \ref{sec:KacRice}}\label{appendix:KacRiceproofs}

For a parameter $\eps>0$, define the truncation of the cube
\[\Sigma_\eps=\Big[-1+e^{-n^{0.6}},1-e^{-n^{0.6}}\Big]^n \;\Big\backslash\;
\ball_{\sqrt{\eps n}}(\bzero).\]
We obtain from Lemma \ref{lem:KacRicestandard} and arguments similar
to \cite{fan2021tap} the following complexity upper bound.

\begin{lemma}\label{lemma:KacRice}
Fix any $\lambda,\eps>0$ and suppose $\bx=\ones$.
Let $T \subset (-1,1)^n$ be any deterministic Borel-measurable set, and let
\[\cC=\Big\{\bbm \in T:\bg(\bbm)=\bzero\Big\}.\]
Then for a $(\lambda,\eps)$-dependent constant $C>0$ and all large $n$,
\[\E[|\cC \cap \Sigma_\eps|]
\leq \sqrt{\frac{n}{2\pi \lambda^2}} \int_{T \times \reals}
\exp\left(n \cdot J(\bbm,y)-\frac{ny^2}{2\lambda^2}+Cn^{0.9} \right)
\de \bbm\,\de y\]
where
\begin{align*}
J(\bbm,y)&=\frac{\lambda^2(1-Q(\bbm))^2}{2}
-\frac{1}{2}\log\Big(2\pi \lambda^2 Q(\bbm)\Big)\\
&\hspace{0.3in}+\frac{1}{n}\sum_{i=1}^n \Big(\log \frac{1}{1-m_i^2}
-\frac{(\arctanh m_i-\lambda^2 M(\bbm)+\lambda^2(1-Q(\bbm))m_i-ym_i)^2}
{2\lambda^2 Q(\bbm)}\Big).
\end{align*}
\end{lemma}
\begin{proof}
Lemma \ref{lem:KacRicestandard} shows
\[\E[|\cC \cap \Sigma_\eps|] \leq \int_{T \cap \Sigma_\eps}
\E\Big[|\det \bH(\bbm)|\;\Big|\;\bg(\bbm)=\bzero\Big]\,p_{\bg(\bbm)}(\bzero)
\de \bbm.\]
The result then follows from applying
\cite[Proposition 3.2 and Lemma 3.3]{fan2021tap} under the identification
$\beta=\lambda$, and noting that the integrand is non-negative so we may
upper-bound the integral over $T \cap \Sigma_\eps$ by that over $T$.
\end{proof}

\begin{proof}[Proof of Lemma \ref{lem:localizetoDeta}]
We introduce $A(\bbm)=n^{-1}\sum_{i=1}^n m_i\arctanh m_i$. If
$\bg(\bbm)=\bzero$, then
\begin{align*}
0=\frac{1}{n}\bbm^\sT \bg(\bbm)&=-\frac{\lambda}{n}\bbm^\sT \bY \bbm
+A(\bbm)+\lambda^2[1-Q(\bbm)]Q(\bbm)\\
&=2\cF_{\TAP}(\bbm)+2H(\bbm)+\frac{\lambda^2}{2}[1-Q(\bbm)]^2
+A(\bbm)+\lambda^2[1-Q(\bbm)]Q(\bbm).
\end{align*}
Consider any $\bbm \in \cB_\delta$ where
$|\cF_{\TAP}(\bbm)-e_\star|<\delta$ and $\bg(\bbm)=\bzero$.
Evaluating $2e_\star+2h_\star+(\lambda^2/2)(1-q_\star)^2
+\lambda^2(1-q_\star)q_\star=-\lambda^2 q_\star$, we obtain for any such $\bbm$
that
\[|A(\bbm)-\lambda^2 q_\star|<(3\lambda^2+4)\delta.\]
Thus, it suffices to show $\P[|\cC| \geq 1]<e^{-cn}$ where we define
\[\cC=\Big\{\bbm \in (-1,1)^n:\;\bg(\bbm)=\bzero,\;
|A(\bbm)-\lambda^2 q_\star|<(3\lambda^2+4)\delta,\;\bbm \in \cB_\delta,
\;\bbm \notin \cD_\eta\Big\}.\]
We have replaced the random condition $|\cF_{\TAP}(\bbm)-e_\star|<\delta$ with
a deterministic condition involving $A(\bbm)$.

Fix $\eps=q_\star/2$. Applying Lemma \ref{lemma:KacRice},
\begin{align*}
\P[|\cC \cap \Sigma_\eps| \geq 1] \leq \E[|\cC \cap \Sigma_\eps|]
&\leq \sqrt{\frac{n}{2\pi \lambda^2}}\int_{(-1,1)^n \times \reals}
\exp\Big(n \cdot J(\bbm,y)-\frac{ny^2}{2\lambda^2}+Cn^{0.9}\Big)\\
&\hspace{1in}\cdot \indic{|A(\bbm)-\lambda^2 q_\star|<
(3\lambda^2+4)\delta,\;\bbm \in \cB_\delta,\;\bbm \notin \cD_\eta}
\de \bbm\,\de y.
\end{align*}
We expand the square in the last term of $J(\bbm,y)$, to write
\begin{align*}
&\frac{1}{n}\sum_{i=1}^n
\frac{(\arctanh m_i-\lambda^2 M(\bbm)+\lambda^2(1-Q(\bbm))m_i-ym_i)^2}
{2\lambda^2 Q(\bbm)}\\
&=\frac{1}{n}\sum_{i=1}^n \frac{[\arctanh m_i-\lambda^2 M(\bbm)]^2}
{2\lambda^2 Q(\bbm)} \\
&\hspace{1in}-\frac{[A(\bbm)-\lambda^2 M(\bbm)^2]
\cdot [y-\lambda^2(1-Q(\bbm))]}{\lambda^2 Q(\bbm)}
+\frac{[y-\lambda^2(1-Q(\bbm))]^2}{2\lambda^2}.
\end{align*}
Then, for any $\bbm \in \cB_\delta$
satisfying $|A(\bbm)-\lambda^2 q_\star|<(3\lambda^2+4)\delta$, we have
\begin{align*}
&~J(\bbm,y)-\frac{y^2}{2\lambda^2}\\
&=\frac{\lambda^2(1-q_\star)^2}{2}
-\frac{1}{2}\log(2\pi \lambda^2 q_\star)
+(1-q_\star) \cdot [y-\lambda^2(1-q_\star)]
-\frac{[y-\lambda^2(1-q_\star)]^2}{2\lambda^2}
-\frac{y^2}{2\lambda^2}\\
&\hspace{0.2in}+\frac{1}{n}\sum_{i=1}^n \log \frac{1}{1-m_i^2}
-\frac{(\arctanh m_i-\lambda^2 q_\star)^2}{2\lambda^2 q_\star}
+r_\lambda(\delta)\Big(1+(\arctanh m_i)^2+|y|\Big),
\end{align*}
for a constant $r_\lambda(\delta)>0$ depending only on $\delta$ and $\lambda$
and satisfying $r_\lambda(\delta) \to 0$ as $\delta \to 0$.
Collecting
\[\frac{\lambda^2(1-q_\star)^2}{2}+(1-q_\star)[y-\lambda^2(1-q_\star)]
-\frac{y^2}{2\lambda^2}=-\frac{[y-\lambda^2(1-q_\star)]^2}{2\lambda^2}\]
and applying this above, we obtain
\begin{align*}
	\P[|\cC \cap \Sigma_\eps| \geq 1]
	&\leq \frac{e^{Cn^{0.9}}}{\sqrt{2}}
	\int_{(-1,1)^n \times \reals}
\exp\left(nr_\lambda(\delta)\left(1+\frac{1}{n}\sum_{i=1}^n
(\arctanh m_i)^2+|y|\right)\right) \cdot \indic{\bbm \notin \cD_\eta}\\
&\hspace{0.5in} \cdot \sqrt{\frac{n}{\pi\lambda^2}}
\exp\Big(-\frac{n(y-\lambda^2(1-q_\star))^2}{\lambda^2}\Big)\de y\\
&\hspace{0.5in}\cdot \prod_{i=1}^n \frac{1}{\sqrt{2\pi \lambda^2 q_\star}} \cdot
\frac{1}{1-m_i^2}\exp\Big(-\frac{(\arctanh m_i-\lambda^2 q_\star)^2}{2\lambda^2
q_\star}\Big)\de m_i
\end{align*}
Let us change variables to $x_i=\arctanh m_i$, and write this as an expectation
over independent random variables
$Y \sim \cN(\lambda^2(1-q_\star),\frac{\lambda^2}{2n})$
and $X_i \sim \cN(\lambda^2q_\star,\lambda^2 q_\star)$. Then the set $\cD_\eta$
is defined by the condition $W_2(\hat{\mu}_X,\cN(\lambda^2 q_\star,\lambda^2
q_\star))<\eta$, where $\hat{\mu}_X$ is the empirical distribution of
$x_1,\ldots,x_n$ and $W_2$ is the Wasserstein-2 distance. Thus, applying this
representation and Cauchy-Schwarz,
\begin{align*}
&~\P[|\cC \cap \Sigma_\eps| \geq 1]\\
&\leq \frac{e^{Cn^{0.9}}}{\sqrt{2}}
\cdot \E\Bigg[\ones\Big\{W_2(\hat{\mu}_X,\cN(\lambda^2 q_\star,\lambda^2
q_\star)) \geq \eta\Big\} \cdot
\exp\left(nr_\lambda(\delta)\left(1+\frac{1}{n}\sum_{i=1}^n
X_i^2+|Y|\right)\right)\Bigg]\\
&\leq \frac{e^{Cn^{0.9}}}{\sqrt{2}}
\cdot \P\Bigg[W_2(\hat{\mu}_X,\cN(\lambda^2 q_\star,\lambda^2
q_\star)) \geq \eta\Bigg]^{1/2}
\E\Bigg[\exp\left(2nr_\lambda(\delta)\left(1+\frac{1}{n}\sum_{i=1}^n
X_i^2+|Y|\right)\right)\Bigg]^{1/2}.
\end{align*}

By Proposition \ref{prop:wasserstein}, for a constant $c_0>0$ depending only on
$(\lambda,\eta)$ and for all large $n$,
\begin{equation}\label{eq:Uprobbound}
\P\Bigg[W_2(\hat{\mu}_X,\cN(\lambda^2 q_\star,\lambda^2
q_\star)) \geq \eta \Bigg]<e^{-c_0n}.
\end{equation}
Applying the moment generating functions of the Gaussian and non-central
chi-squared distributions,
together with $\log(1+x) \leq x$, for any $t<1/2$ we have
\begin{align*}
\E[e^{nt|Y|}]
&\leq \E[e^{ntY}]+\E[e^{-ntY}]
\leq 2e^{nt\lambda^2(1-q_\star)+nt^2 \lambda^2/4},\\
\E[e^{t\sum_{i=1}^n X_i^2}]
&\leq \exp\bigg(
\frac{nt\lambda^4q_\star^2}{1-2t}\bigg)
\cdot\big(1-2t\big)^{n/2}
\leq \exp\bigg(\frac{nt\lambda^4q_\star^2}{1-2t}-nt\bigg).
\end{align*}
Now choosing $\delta=\delta(\lambda,\eta)$ sufficiently small so that
$r_\lambda(\delta)$ is sufficiently small, we may guarantee
\[\E\Bigg[\exp\left(2nr_\lambda(\delta)\left(1+\frac{1}{n}\sum_{i=1}^n
X_i^2+|Y|\right)\right)\Bigg]^{1/2} \leq e^{c_0n/8}\]
where $c_0$ is the constant in (\ref{eq:Uprobbound}). Thus
$\P[|\cC \cap \Sigma_\eps| \geq 1]<e^{-c_0n/8}$.

Finally, recalling $\eps=q_\star/2$, the set $\cB_\delta$ does not
intersect the ball $\ball_{\sqrt{\eps n}}(\bzero)$ for any $\delta<q_\star/2$.
Thus
\[|\cC \cap \Sigma_\eps^c| \subseteq \Big\{\bbm \in (-1,1)^n:\bg(\bbm)=\bzero
\text{ and }\|\bbm\|_\infty \in (1-e^{-n^{0.6}},1)\Big\}.\]
On the event of probability $1-e^{-cn}$ where $\|\bW\|_{\op}<3$ and hence
$\|\bY\|_{\op}<\lambda+3$, we have
$\big\|\lambda \bY\bbm-\lambda^2[1-Q(\bbm)]\bbm\big\|_2
<(2\lambda^2+3\lambda)\sqrt{n}$
for all $\bbm \in (-1,1)^n$. When $\|\bbm\|_\infty \in (1-e^{-n^{0.6}},1)$,
we also have $\|\arctanh(\bbm)\|_2 \geq n^{0.6}/2$.
Applying this to (\ref{eq:grad}), we must have
$\bg(\bbm) \neq \bzero$. Thus also
$\P[|\cC \cap \Sigma_\eps^c| \geq 1]<e^{-cn}$,
so $\P[|\cC| \geq 1] \leq e^{-cn}$ for all large $n$, as desired.
\end{proof}

\subsection{Proofs for Section
\ref{sec:strongconvex}}\label{appendix:strongconvexproofs}

Let us write
\begin{equation}\label{eq:Hchi}
\bH^\chi(\bbm)=\chi\left(-\lambda \bY-\frac{2\lambda^2}{n}\bbm\bbm^\sT\right)
+\diag\Big(\frac{1}{1-\bbm^2}\Big)+\lambda^2[1-Q(\bbm)]\id, \qquad \chi \in
\{+,-\}
\end{equation}
and
\[\ell_\eps^\chi(\bbm)=\inf\Big\{\lambda_{\min}\big(\bH^\chi(\bu)\big):
\bu \in (-1,1)^n \cap \ball_{\sqrt{\eps n}}(\bbm)\Big\}.\]
These coincide with the definitions of $\ell_\eps^+$ and $\ell_\eps^-$
in Sections \ref{sec:strongconvex} and \ref{sec:AMP-stab}.

\begin{lemma}\label{lem:Hcondlaw}
Suppose $\bx=\ones$. For $\bu,\bbm \in (-1,1)^n$, let $\proj_\bbm^\perp$ be the
orthogonal projection onto the $(n-1)$-dimensional subspace orthogonal to
$\bbm$, and let
	\begin{align}
		\bz(\bbm)
			&=\arctanh(\bbm)
			-\lambda^2 M(\bbm)\ones+\lambda^2[1-Q(\bbm)]\bbm,\label{eq:zm}
		\\
		\bM^\chi(\bbm,\bu)
			&=\chi\Bigg(\frac{1}{\|\bbm\|_2^2}\Big(\bbm\bz(\bbm)^\sT+\bz(\bbm)\bbm^\sT\Big)
			-\frac{\langle \bbm,\bz(\bbm) \rangle}{\|\bbm\|_2^4} \cdot
\bbm\bbm^\sT\nonumber\\
&\hspace{1in}+\frac{\lambda^2}{n}\ones\ones^\sT
+\frac{2\lambda^2}{n}\bu\bu^\sT\Bigg)
		-\diag\Big(\frac{1}{1-\bu^2}\Big)
			-\lambda^2[1-Q(\bu)]\id.
		\label{eq:condhessshift}
	\end{align}
	Then for any $\bbm \in (-1,1)^n$ and either $\chi \in \{+,-\}$,
	conditional on the event $\bg(\bbm)=\bzero$,
	\begin{equation}\label{eq:condhess}
		\Big\{-\bH^\chi(\bu):\bu \in (-1,1)^n\Big\}\Big|_{\bg(\bbm)=\bzero}
			\overset{d}{=}
			\Big\{
				\lambda \proj_\bbm^\perp \tilde{\bW} \proj_\bbm^\perp
				+\bM^\chi(\bbm,\bu): \bu \in (-1,1)^n \Big\}
	\end{equation}
	where $\tilde{\bW} \sim \GOE(n)$ is an independent copy of $\bW$, and
	this holds as an equality in law of two Gaussian processes indexed
	by $\bu \in (-1,1)^n$.
\end{lemma}
\begin{proof}
	Writing $\bY=(\lambda/n)\ones\ones^\sT+\bW$,
$\bg(\bbm)=\bzero$ is equivalent to $\bW\bbm=\lambda^{-1}\bz(\bbm)$. Hence conditioned on
	this event, the law of $\bW$ is (see e.g.\ \cite[Lemma 4.1]{fan2021tap})
	\[\bW|_{\bg(\bbm)=\bzero}\overset{d}{=}
	\proj_\bbm^\perp \tilde{\bW} \proj_\bbm^\perp
	+\frac{1}{\lambda \|\bbm\|_2^2}\Big(\bbm\bz(\bbm)^\sT+\bz(\bbm)\bbm^\sT\Big)
	-\frac{\langle \bbm,\bz(\bbm)
	\rangle}{\lambda \|\bbm\|_2^4} \cdot \bbm\bbm^\sT.\]
	The result follows from substituting this into the expression for
$\bH^\chi(\bu)$ in (\ref{eq:Hchi}).
\end{proof}

\begin{proof}[Proof of Lemma \ref{lem:gordon-variational}(a)
and Lemma \ref{lem:gordon-variational-gen}]

Recall that by Lemma \ref{lem:Hcondlaw}, we have 
\[
\begin{aligned}
	&~ \Big\{ -\bH^\chi(\bu): \bu \in (-1, 1)^n \Big\}\Big\vert_{\bg(\bbm) =
\bzero} \stackrel{d}{=}&~ \Big\{ \lambda \proj_{\bbm}^\perp \tilde \bW
\proj_{\bbm}^\perp + \bM^\chi(\bbm, \bu): \bu \in (-1, 1)^n \Big\}.
\end{aligned}
\]
Then
\[-\ell_\eps^\chi(\bbm)\Big|_{\bg(\bbm)=\bzero} \overset{d}{=}
\sup_{\bu \in (-1,1)^n \cap \ball_{\sqrt{\eps n}}(\bbm)} 
\sup_{\bv \in \reals^n:\| \bv \|_2 = 1} \lambda \langle \proj_{\bbm}^\perp
\bv,\tilde{\bW}\proj_{\bbm}^\perp \bv \rangle+\langle \bv,\bM^\chi(\bbm,\bu)\bv
\rangle.\]
We introduce $\bg \sim \cN(\bzero,\id) \in \reals^n$ and define
the auxiliary Gaussian process
	\[G (\bbm, \bu, \bv) 
		= 
		2 \lambda \sqrt \frac{1}{n} \<\proj_{\bbm}^\perp \bv, \bg\> \| \proj_{\bbm}^\perp \bv \|_2 
		+ \langle \bv,\bM^\chi(\bbm,\bu) \bv\rangle.\]
Then by Lemma \ref{lem:slepian}, we obtain for any fixed
$\bbm \in \cD_\eta$ that
\begin{equation}\label{eq:slepianH}
\E\Big[-\ell_\eps^\chi(\bbm)\;\Big|\;\bg(\bbm)=\bzero\Big]
\leq \E\Big[ \sup_{\bu \in (-1,1)^n \cap \ball_{\sqrt{\eps n}}(\bbm)} 
\sup_{\bv \in \reals^n:\| \bv \|_2 = 1} G(\bbm, \bu, \bv) \Big]. 
\end{equation}

We now analyze $G(\bbm,\bu,\bv)$. Throughout the proof, we write
$r(n,\eps,\eta)$ for any error term depending only on $n,\eps,\eta,\lambda,K'$
and satisfying
\[\lim_{\eps,\eta \to 0} \lim_{n \to \infty} r(n,\eps,\eta)=0.\]
We write $C,c>0$ for constants depending only on $(\lambda,K')$ and changing
from instance to instance. Suprema over $(\bu,\bv)$ are
implicitly over $\bu \in (-1,1)^n \cap \ball_{\sqrt{\eps n}}(\bbm)$ and $\bv \in
\reals^n$ with $\|\bv\|_2=1$, unless otherwise stated.

Writing
$\proj_\bbm=\id-\proj_{\bbm}^\perp$ for the projection onto the span of $\bbm$,
observe that
\[\E\left[\sup_{\|\bv\|_2=1} 2\lambda \sqrt{\frac{1}{n}}
\langle \proj_\bbm \bv,\bg \rangle\|\proj_\bbm^\perp \bv\|_2\right]
\leq 2\lambda\sqrt{\frac{2}{\pi n}}.\]
Then
\[\E\Big[ \sup_{\bu,\bv} G(\bbm, \bu, \bv) \Big]
\leq \E\Big[ \sup_{\bu,\bv} 2\lambda \sqrt{\frac{1}{n}}
\langle \bv,\bg \rangle\|\proj_{\bbm}^\perp \bv\|_2
+\langle \bv,\bM^\chi(\bbm,\bu)\bv \rangle\Big]+r(n,\eps,\eta).\]
In the expression $\langle \bv,\bM^\chi(\bbm,\bu) \bv \rangle$, applying
$\|\bu-\bbm\|_2^2 \leq \eps n$, we may bound
\[\frac{2\lambda^2}{n}\Big|\langle \bu,\bv \rangle^2
-\langle \bbm,\bv \rangle^2\Big|+\lambda^2\Big|Q(\bbm)-Q(\bu)\Big|
\leq \frac{3\lambda^2}{n}\|\bbm-\bu\|_2\|\bbm+\bu\|_2
\leq 6\lambda^2\eps^{1/2}\]
to replace $Q(\bu)$ and $\langle \bu,\bv \rangle^2$ by $Q(\bbm)$ and $\langle
\bbm,\bv \rangle^2$, up to $r(n,\eps,\eta)$ error.
We may then apply, for $\bbm \in \cD_\eta$,
\begin{align*}
\frac{1}{n}\|\bbm\|_2^2&=Q(\bbm)=q_\star+r(n,\eps,\eta),\\
\frac{1}{n}\langle \bbm,\bz(\bbm) \rangle
&=2\lambda^2(1-q_\star)q_\star+r(n,\eps,\eta),\\
\frac{1}{n}\|\bz(\bbm)\|^2 &\leq C,\\
\frac{1}{\sqrt{n}}\langle \bz(\bbm),\bv \rangle
&=\frac{1}{\sqrt{n}}\langle z(\bbm),\bv \rangle+r(n,\eps,\eta),
\end{align*}
where $z(\bbm)$ in this last equality denotes the entry-wise application of
$z(m)=\arctanh m-\lambda^2 q_\star+\lambda^2(1-q_\star)m$ to $\bbm$,
and differs from $\bz(\bbm)$ defined in (\ref{eq:zm})
in the replacement of $M(\bbm)$ and $Q(\bbm)$ by $q_\star$. This yields
\[\E\Big[\sup_{\bu,\bv} G(\bbm, \bu, \bv) \Big]
\leq \E\Big[\sup_{\bu,\bv} G_\star(\bbm,\bu,\bv) \Big]+r(n,\eps,\eta)\]
where we define
\begin{align*}
		G_\star(\bbm, \bu, \bv) 
			=&~ 
			2 \lambda \sqrt \frac{1}{n} \<\bv, \bg\> \| \proj_{\bbm}^\perp \bv \|_2 
			+ \chi\Bigg(
			2\frac{\<z(\bbm), \bv\>}{\sqrt {n q_\star}} \frac{\<\bbm, \bv\> }{\| \bbm \|_2} 
			- 
			2\lambda^2 (1 - q_\star) \frac{\<\bbm, \bv\>^2}{\| \bbm \|_2^2} 
		\\
			&\hspace{0.5in}+
			\frac{\lambda^2}{n} \< \bv,\ones\>^2
			+
			2\lambda^2 q_\star \frac{\< \bbm, \bv\>^2}{\| \bbm
\|_2^2}\Bigg)
			-
			\sum_{i = 1}^n \frac{v_i^2}{1 - u_i^2}  
			- 
			\lambda^2 (1 - q_\star).
	\end{align*}

Next, we analyze $G_\star(\bbm,\bu,\bv)$. We introduce
\[p=\sqrt{q_\star} \cdot \frac{\langle \bbm,\bv \rangle}{\|\bbm\|_2} \in
[-\sqrt{q_\star},\sqrt{q_\star}],
\quad u=\frac{\langle \bv,\ones \rangle}{\sqrt{n}} \in [-1,1],
\quad K=\Big\{(p,u):|p| \leq \sqrt{q_\star},|u| \leq 1\Big\}\]
to bound the supremum over $\bv$ for fixed $\bu$ as
	\begin{align}
	&\sup_{\bv \in \reals^n:\|\bv\|_2 = 1} 
			G_\star(\bbm, \bu, \bv)\nonumber\\
		&=  
			\sup_{(p,u) \in K}
			\sup_{\bv \in \reals^n:\;\|\bv\|_2 = 1,\,
				\frac{\< \bbm, \bv\> }{ \| \bbm \|_2} = p/q_\star^{1/2},\,
				\frac{\<\bv,\ones\>}{\sqrt{n}} = u} 
				G_\star(\bbm, \bu, \bv)\nonumber \\
		&\le 
			\sup_{(p,u) \in K}
			\inf_{(\alpha,\kappa,\gamma) \in K'}
			\sup_{\bv \in \reals^n} 
				\Bigg\{ 
				2 \lambda \sqrt \frac{1}{n} \<\bv, \bg\> \Big(1 - \frac{p^2}{q_\star}\Big)^{1/2} 
				 \nonumber\\ 
				&~~~+\chi\Bigg(2\frac{\<z(\bbm), \bv\>}{\sqrt
n}\frac{p}{q_\star}
- 2 \lambda^2p^2 \frac{1 - q_\star}{q_\star} 
			+ \lambda^2 u^2 
			+ 2 \lambda^2 p^2 \Bigg) \nonumber\\
&~~~~- \sum_{i = 1}^n \frac{v_i^2}{1 - u_i^2} 
			- \lambda^2 (1 - q_\star)
			+ \gamma \Big(\| \bv \|_2^2-1\Big)
			+ 
			\chi\alpha \bigg(\frac{\< \bv, \ones\>}{\sqrt{n}}-u\bigg) 
			+ 
			\chi\kappa \bigg(
					\frac{\< \bbm, \bv\>}{\sqrt{n}} - p \frac{\| \bbm \|_2}{\sqrt{n q_\star}}
					\bigg)\Bigg\}.
		\nonumber
\end{align}
Collecting the terms above that depend on $\bv$, 
substituting $x_i$ for $v_i\sqrt{n}$, and defining
\[\Theta_\lambda(g,m,u)=\sup_{x \in \reals} \bigg\{ 2 \lambda x g\Big(1 - \frac{p^2}{q_\star}\Big)^{1/2} 
- \frac{x^2}{1 - u^2} 
					+ \gamma x^2
					+\chi x\Big[2z(m) \cdot \frac{p}{q_\star} 
					+ \alpha
					+ \kappa m\Big]  \bigg\},\]
this yields
\begin{align}
\sup_{\bv \in \reals^n:\|\bv\|_2 = 1} 
			G_\star(\bbm, \bu, \bv)
		&\leq 
		\sup_{(p,u) \in K} \inf_{(\alpha,\kappa,\gamma) \in K'}
		\bigg\{\chi\Big[
2\lambda^2p^2+\lambda^2u^2-2\lambda^2p^2\frac{1-q_\star}{q_\star}  \nonumber\\
		 	&\hspace{0.2in}
		 	-  \alpha u
-  \kappa p\frac{\| \bbm \|_2}{\sqrt{n q_\star}}
\Big]
-\lambda^2(1-q_\star)-\gamma
			+\frac{1}{n}\sum_{i=1}^n
\Theta_\lambda(g_i,m_i,u_i)\bigg\}.\label{eq:supGstar}
	\end{align}

Since $\gamma<1$ and $1/(1-u^2) \geq 1$, the maximum over $x$
in the definition of $\Theta_\lambda(g,m,u)$ is achieved at 
	\[
		x = \Big(2 \lambda  g(1 - p^2 / q_\star)^{1/2} + \chi[2 z(m)  p /
q_\star + \alpha   + \kappa m]\Big)\Big/\Big(\frac{2}{1 - u^2} - 2\gamma\Big),
	\]
yielding the explicit form
\begin{equation}\label{eq:Thetalambda}
\Theta_\lambda(g,m,u)=\Big(2 \lambda  g (1 - p^2
/ q_\star)^{1/2} + \chi[2 z(m) p / q_\star + \alpha + \kappa
m]\Big)^2\Big/\Big(\frac{4}{1 - u^2} - 4 \gamma\Big).
\end{equation}
Since $\gamma$ is less than and bounded away from 1 on the compact domain $K'$,
\[\left|\frac1{\frac1{1-u^2} - \gamma} - \frac1{\frac1{1-m^2} - \gamma}\right| 
= \left|\frac{1-u^2}{1-\gamma(1-u^2)} - \frac{1-m^2}{1-\gamma(1-m^2)}\right|
\leq C|u-m|\]
for a constant $C>0$ depending on $K'$. Then, applying this to
(\ref{eq:Thetalambda}) and recalling the definition of $z(m)$,
\[\left|\frac{1}{n}\sum_{i=1}^n \Theta_\lambda(g_i,m_i,u_i)
-\Theta_\lambda(g_i,m_i,m_i)\right| \leq
\frac{1}{n}\sum_{i=1}^n C\Big(1+g_i^2+(\arctanh m_i)^2\Big)
|u_i-m_i|\]
for a different constant $C>0$ depending also on $\lambda$.
Now applying this to (\ref{eq:supGstar}), and 
taking also the supremum over $\bu$ and the expectation over $\bg$,
we arrive at
\[\E\Big[\sup_{\bu,\bv} G_\star(\bbm,\bu,\bv)\Big] \leq
\E\Big[\sup_{(p,u) \in K} \inf_{(\alpha,\kappa,\gamma) \in K'}
-H_n^\chi(\bg,\bbm;p,u,\alpha,\kappa,\gamma)\Big]+R_n(\bbm)\]
where
\begin{align}
H_n^\chi(\bg,\bbm;p,u,\alpha,\kappa,\gamma)
&=-\chi\Big[2\lambda^2 p^2+\lambda^2u^2-2\lambda^2(1-q_\star)p^2/q_\star
-\alpha u-\kappa p\Big]\nonumber\\
&\hspace{1in}+\lambda^2(1-q_\star)+\gamma-\frac{1}{n}\sum_{i=1}^n
\Theta_\lambda(g_i,m_i,m_i),\label{eq:Hngm}
\end{align}
and
\[
\big|R_n(\bbm)\big|
\leq C\left(\left|\frac{\|\bbm\|_2}{\sqrt{nq_\star}}-1\right|
+\E\left[\sup_\bu \frac{1}{n}\sum_{i=1}^n \Big(1+g_i^2+(\arctanh m_i)^2\Big)
|u_i-m_i| \right]\right).\]

To bound this remainder $R_n(\bbm)$, note that $\bbm \in \cD_\eta$ implies
for the first term
\[\frac{\|\bbm\|_2}{\sqrt{nq_\star}}=1+r(n,\eps,\eta).\]
For the second term,
applying $n^{-1}\|\bu-\bbm\|_2^2 \leq \eps$ and Markov's inequality,
\[|u_i-m_i| \leq \eps^{1/4}+2\indic{|u_i-m_i| \geq \eps^{1/4}}, \qquad
\frac{1}{n}\sum_{i=1}^n \indic{|u_i-m_i| \geq \eps^{1/4}} \leq \eps^{1/2}.\]
Set $a_i=\arctanh m_i$, define $|g|_{(1)} \geq \ldots \geq |g|_{(n)}$ as the
values $|g_1|,\ldots,|g_n|$ sorted in decreasing order, and
define similarly $|a|_{(1)} \geq \ldots \geq |a|_{(n)}$. Then
\begin{align*}
\E\left[\sup_{\bu} \frac{1}{n}\sum_{i=1}^n
\Big(1+g_i^2+a_i^2\Big)|u_i-m_i|\right]
&\leq \frac{1}{n}\sum_{i=1}^n \Big(2+a_i^2\Big)\eps^{1/4}+
\frac{C}{n}\sum_{i=1}^{\lfloor \eps^{1/2}n \rfloor} 
\Big(1+\E[|g|_{(i)}^2]+|a|_{(i)}^2 \Big).
\end{align*}
The Wasserstein-2 distance between the empirical distribution of
$(a_1,\ldots,a_n)$ and the law $X_\star \sim \cN(\lambda^2 q_\star,\lambda^2
q_\star)$ is at most $\eta$, by the condition $\bbm \in \cD_\eta$. Then,
letting $q_\eps$ be the $1-\eps^{1/2}$ quantile of the distribution of
$|X_\star|$ and applying Proposition \ref{prop:wassersteinbound}(a),
\[\frac{1}{n}\sum_{i=1}^{\lfloor \eps^{1/2}n \rfloor}
|a|_{(i)}^2 \leq \E\big[X_\star^2 \indic{|X_\star| \geq q_\eps}\big]+C\eta
\leq r(n,\eps,\eta).\]
Similarly, letting $q_\eps$ be the $1-\eps^{1/2}$ quantile of $|G|$ for $G \sim
\cN(0,1)$ and applying Proposition \ref{prop:wasserstein},
\[\E\left[\frac{1}{n}\sum_{i=1}^{\lfloor \eps^{1/2}n \rfloor} |g|_{(i)}^2\right]
\leq \E\big[G^2\indic{|G| \geq q_\eps}\big]+r(n,\eps,\eta)
\leq r(n,\eps,\eta).\]
Combining these observations yields $|R_n(\bbm)| \leq r(n,\eps,\eta)$, and thus
\begin{equation}\label{eq:GHbound}
\E\Big[\sup_{\bu,\bv} G(\bbm,\bu,\bv)\Big]
\leq \E\Big[\sup_{(p,u) \in K} \inf_{(\alpha,\kappa,\gamma) \in K'}
-H_n^\chi(\bg,\bbm;p,u,\alpha,\kappa,\gamma)\Big]+r(n,\eps,\eta).
\end{equation}

Finally, comparing $H_\lambda^\chi$ in (\ref{eq:Hlambda}) with $H_n^\chi$ in
(\ref{eq:Hngm}), observe that
\begin{align}
&H_\lambda^\chi(p,u;\alpha,\kappa,\gamma)
-H_n^\chi(\bg,\bbm;p,u,\alpha,\kappa,\gamma)\nonumber\\
&=\frac{1}{n}\sum_{i=1}^n \Theta_\lambda(g_i,m_i,m_i)
-\E_{G \sim \cN(0,1),m \sim \mu_\star}\Big[\Theta_\lambda(G,m,m)\Big],
\label{eq:HlambdaHn}
\end{align}
where $\Theta_\lambda$ is defined in (\ref{eq:Thetalambda}). Let us now
make the dependence of $\Theta_\lambda$ on $p,\alpha,\kappa,\gamma$ explicit,
and write
\[\Theta_\lambda(g,m,m)=\Xi\Big(g,m;2\lambda(1-p^2/q_\star)^{1/2},
2p/q_\star,\alpha,\kappa,\gamma\Big)\]
where
\[\Xi(g,m;p_1,p_2,\alpha,\kappa,\gamma)=(p_1 g+\chi p_2z(m)+\chi \alpha
+\chi \kappa m)^2 \Big/\Big(\frac{4}{1-m^2}-4\gamma\Big).\]
Bounding (\ref{eq:HlambdaHn}) is similar to the argument in Lemma
\ref{lem:global_min}: Set $\omega=(p_1,p_2,\alpha,\kappa,\gamma)$ and
\[K_\omega=\Big\{\omega:p_1 \in [0,2\lambda],p_2 \in
[-2/\sqrt{q_\star},2/\sqrt{q_\star}],(\alpha,\kappa,\gamma) \in K'\Big\}.\]
For any fixed $m \in (-1,1)$ and $\omega \in K_\omega$, let $G \sim \cN(0,1)$,
and note that $\Xi_\lambda(G,m;\omega)$ is the square of a Gaussian variable
with mean bounded by $C(1+|z(m)|)$ and variance bounded by $C$.
Then $\Xi(G,m;\omega)-\E[\Xi(G,m;\omega)]$ is $(C(1+z(m)^2),C)$-sub-Gamma, see
e.g.\ \cite[Definition G.2 and Proposition G.5]{miolane2021distribution}. So
$n^{-1}\sum_{i=1}^n \Xi(g_i,m_i;\omega)-\E[\Xi(G,m_i;\omega)]$
is $(C(n+\|\bz(\bbm)\|_2^2)/n^2,C/n)$-sub-Gamma for any fixed
$\bbm \in \cD_\eta$. Applying $\|\bz(\bbm)\|_2^2/n \leq C$, we obtain
that this is $(C/n,C/n)$-sub-Gamma, so Bernstein's inequality yields
\begin{equation}\label{eq:Xibound}
\E\left[\left|\frac{1}{n}\sum_{i=1}^n \Xi(g_i,m_i;\omega)
-\E[\Xi(G,m_i;\omega)]\right|\right] \leq C/\sqrt{n}.
\end{equation}
To obtain uniform control over $K_\omega$, let $\cN$ be a $n^{-1/20}$-net of
$K_\omega$ of cardinality $|\cN| \leq Cn^{1/4}$. Observe that
$n^{-1}\sum_{i=1}^n \Xi(g_i,m_i;\omega)$ is
$C(1+\|\bg\|_2^2/n+\|\bz(\bbm)\|_2^2/n)$-Lipschitz in $\omega \in K_\omega$.
Then, applying (\ref{eq:Xibound}) over $\cN$, we get
\[
\begin{aligned}
&~\E\left[\sup_{\omega \in K_\omega}
\left|\frac{1}{n}\sum_{i=1}^n \Xi(g_i,m_i;\omega)
-\E[\Xi(G,m_i;\omega)]\right|\right] \\
\leq&~ |\cN| \cdot \frac{C}{\sqrt{n}}
+n^{-1/20} \cdot
\E\left[C\left(1+\frac{\|\bg\|_2^2}{n}+\frac{\|\bz(\bbm)\|_2^2}{n}\right)\right]
\leq r(n,\eps,\eta).
\end{aligned}
\]
Now observe that the function
$m \mapsto \E_{G \sim \cN(0,1)}[\Xi(G,m;\omega)]$ satisfies
\[
\begin{aligned}
&~\Big|\E_{G \sim \cN(0,1)}[\Xi(G,m;\omega)]
-\E_{G \sim \cN(0,1)}[\Xi(G,m';\omega)]\Big|\\
\leq&~ C(1+|\arctanh m|+|\arctanh m'|)|\arctanh m-\arctanh m'|
\end{aligned}
\]
for all $m,m' \in (0,1)$, uniformly over $\omega \in K_\omega$ and for a
constant $C>0$ depending only on $(\lambda,K')$. Then, applying
Proposition \ref{prop:wasserstein}(b) and $\bbm \in \cD_\eta$, we have
\[\sup_{\omega \in K_\omega}
\left|\frac{1}{n}\sum_{i=1}^n \E[\Xi(G,m_i;\omega)]
-\E_{G \sim \cN(0,1),m \sim \mu_\star} [\Xi(G,m;\omega)]\right|
\leq r(n,\eps,\eta).\]
Combining these two bounds, we obtain
\[\E\left[\sup_{(p,u) \in K} \sup_{(\alpha,\kappa,\gamma) \in K'}
\left|\frac{1}{n}\sum_{i=1}^n \Theta_\lambda(g_i,m_i,m_i)
-\E_{G \sim \cN(0,1),m \sim \mu_\star}\Big[\Theta_\lambda(G,m,m)\Big]
\right|\right] \leq r(n,\eps,\eta).\]
Applying this back to (\ref{eq:GHbound}) and (\ref{eq:HlambdaHn}), we get
\[\E\Big[\sup_{\bu,\bv} G(\bbm,\bu,\bv)\Big]
\leq \E\Big[\sup_{(p,u) \in K} \inf_{(\alpha,\kappa,\gamma) \in K'}
-H_\lambda^\chi(p,u;\alpha,\kappa,\gamma)\Big]+r(n,\eps,\eta).\]
This holds for all $\bbm \in \cD_\eta$, where
$r(n,\eps,\eta)$ is independent of the specific point $\bbm \in \cD_\eta$.
Then, applying this to (\ref{eq:slepianH}), negating the sign, and choosing
$\eps,\eta>0$ sufficiently small depending on $t$ yields the result.
\end{proof}

\begin{proof}[Proof of Lemma \ref{lem:gordon-variational}(b)]
We write as shorthand $H_\lambda=H_\lambda^+$ and
$\E$ for $\E_{m \sim \mu_\star}$. Recall
\begin{align*}
H_\lambda(p,u;\alpha,\kappa,\gamma) =&~ -2\lambda^2p^2-\lambda^2u^2
+2\lambda^2(1 - q_\star) p^2 / q_\star  + \alpha u  +  \kappa p +
\gamma+\lambda^2(1-q_\star) \\
	& \quad- \E\Big[ \Big(4 \lambda^2  (1 - p^2 / q_\star) +
(2 z(m)  p / q_\star + \alpha   + \kappa m)^2 \Big) \Big/ \Big( \frac{4}{1 -
m^2} - 4 \gamma \Big) \Big].
\end{align*}
Note that the dependence of $H_\lambda$ on the first four variables
$(p,u;\alpha,\kappa)$ is homogeneous of degree 2. Collecting these quadratic
terms yields
\begin{equation}\label{eq:Hlambdaquad}
H_\lambda(p,u;\alpha,\kappa,\gamma)
= (p,u,\kappa,\alpha)^\sT \bA^{(\gamma)}(p,u,\kappa,\alpha)
+\gamma + \lambda^2(1-q_\star) -
\lambda^2 \E\Big[\Big(\frac1{1-m^2} - \gamma\Big)^{-1}\Big]
\end{equation}
where
\[
    \bA^{(\gamma)}
	=
	\begin{pmatrix} 
		\bA_{11}^{(\gamma)} & \bA_{12}^{(\gamma)} 
		\\ 
		\bA_{21}^{(\gamma)} & \bA_{22}^{(\gamma)}
	\end{pmatrix}
\]
and
\begin{align*}
	\bA_{11}^{(\gamma)}
		&=
		\begin{pmatrix} 
			-2\lambda^2+2\lambda^2\frac{1-q_\star}{q_\star} 
				+ 
				\frac{\lambda^2}{q_\star}\E\Big[\Big(\frac{1}{1-m^2} 
					- \gamma\Big)^{-1}
				\Big]
			-	
				\frac{1}{q_\star^2} \E\Big[
					z(m)^2\Big(\frac{1}{1-m^2}-\gamma\Big)^{-1}
				\Big] 
				& 0 
			\\ 
			0 & -\lambda^2
		\end{pmatrix}
	\\
	\bA_{12}^{(\gamma)}={\bA_{21}^{(\gamma)\sT}}
		&=
		\begin{pmatrix} 
				\frac{1}{2}-\frac{1}{2q_\star}\E\Big[mz(m)\Big(\frac1{1-m^2}-\gamma\Big)^{-1}\Big]
	& 
			-\frac{1}{2q_\star}\E\Big[z(m)\Big(\frac1{1-m^2} - \gamma\Big)^{-1}\Big] 
						\\
			0 & \frac{1}{2}
		\end{pmatrix}
	\\
	\bA_{22}^{(\gamma)}
		&=
		-\frac{1}{4}
		\begin{pmatrix} 
			\E\Big[m^2\Big(\frac1{1-m^2} - \gamma\Big)^{-1}\Big] 
				& 
				\E\Big[m\Big(\frac1{1-m^2}-\gamma\Big)^{-1}\Big] 
			\\
			\E\Big[m\Big(\frac1{1-m^2} - \gamma\Big)^{-1}\Big] 
				& 
				\E\Big[\Big(\frac1{1-m^2}-\gamma\Big)^{-1}\Big] 
		\end{pmatrix}.
\end{align*}

We first specialize this matrix $\bA^{(\gamma)}$ to $\gamma=0$, and apply the
representation $z(m)=\lambda\sqrt{q_\star}G+\lambda^2(1-q_\star)m$
where $m=\tanh(\lambda^2 q_\star+\lambda\sqrt{q_\star} G)$ and $G \sim
\cN(0,1)$. Recalling $q_\star=\E[m]=\E[m^2]$ and $b_\star=\E[m^3]=\E[m^4]$ from
Proposition \ref{prop:qstar}, this yields
\begin{align*}
\bA_{11}^{(0)} &=
\begin{pmatrix} 
-2\lambda^2+3\lambda^2\frac{1-q_\star}{q_\star} 
-\frac{\lambda^4(1-q_\star)^2}{q_\star^2}(q_\star-b_\star)
-\frac{2\lambda^3(1-q_\star)}{q_\star^{3/2}} \E[Gm(1-m^2)] 
-\frac{\lambda^2}{q_\star}\E[G^2(1-m^2)] & 0 \\ 0 & -\lambda^2
\end{pmatrix} \\
\bA_{12}^{(0)} &=
\begin{pmatrix} 
\frac{1}{2}-\frac{\lambda^2(1-q_\star)}{2q_\star}(q_\star-b_\star)
-\frac{\lambda}{2q_\star^{1/2}}\E[Gm(1-m^2)] &
-\frac{\lambda^2(1-q_\star)}{2q_\star}(q_\star-b_\star)
+\frac{\lambda}{2q_\star^{1/2}}\E[Gm^2] \\
0 & \frac{1}{2} \end{pmatrix} \\
\bA_{22}^{(0)} &= -\frac{1}{4}
\begin{pmatrix} q_\star-b_\star & q_\star-b_\star \\
q_\star-b_\star & 1-q_\star \end{pmatrix}.
\end{align*}
Now recalling the matrices $\bA_{11},\bA_{12},\bA_{21},\bA_{21}$
defined in (\ref{eq:A22})
from the proof of Lemma \ref{lem:global_min_analysis}, we may check that
\begin{align*}
\bA_{22}^{(0)}&=\frac{1}{4}\bA_{22}\\
\bA_{12}^{(0)}&=\frac{1}{2}\left(
\bA_{12}+\begin{pmatrix} \frac{\lambda^2(1-q_\star)}{q_\star}
& 0 \\ 0 & 0 \end{pmatrix} \bA_{22}\right)\\
\bA_{11}^{(0)}&=\bA_{11}+\bA_{12} 
\begin{pmatrix} \frac{\lambda^2(1-q_\star)}{q_\star} & 0 \\ 0 & 0 \end{pmatrix}
+\begin{pmatrix} \frac{\lambda^2(1-q_\star)}{q_\star} & 0 \\ 0 & 0 \end{pmatrix}
\bA_{21}
+\begin{pmatrix} \frac{\lambda^2(1-q_\star)}{q_\star} & 0 \\ 0 & 0 \end{pmatrix}
\bA_{22}
\begin{pmatrix} \frac{\lambda^2(1-q_\star)}{q_\star} & 0 \\ 0 & 0 \end{pmatrix}
\end{align*}
It was verified in the proof of Lemma \ref{lem:global_min_analysis} that
$\bA_{22} \prec 0$ and $\bA_{11}-\bA_{12}\bA_{22}^{-1}\bA_{21} \succ 0$
strictly. Then also $\bA_{22}^{(0)} \prec 0$ and
\[\bA_{11}^{(0)}-\bA_{12}^{(0)}(\bA_{22}^{(0)})^{-1}\bA_{21}^{(0)}
=\bA_{11}-\bA_{12}\bA_{22}^{-1}\bA_{21} \succ 0\]
strictly. Then by continuity in $\gamma$, also $\bA_{22}^{(\gamma)} \prec 0$ and
$\bA_{11}^{(\gamma)}-\bA_{12}^{(\gamma)}(\bA_{22}^{(\gamma)})^{-1}
\bA_{21}^{(\gamma)} \succ 0$ strictly for all $|\gamma|$ sufficiently small.

Applying this to the form (\ref{eq:Hlambdaquad}) for $H_\lambda$, the condition
$\bA_{22}^{(\gamma)} \prec 0$ shows that for fixed $p,u,\gamma$, the quadratic
function $H_\lambda$ is strictly concave in $(\alpha,\kappa)$. At
$(p,u)=(0,0)$, $H_\lambda$ is maximized at $(\alpha,\kappa)=(0,0)$. Here,
$(0,0,\gamma)$ belongs to the interior of the given domain
$K'$ for $|\gamma|$ sufficiently
small. Thus, for some $(\lambda,K')$-dependent constant $c>0$ and any
$|\gamma|<c$, the function
\[\bar{H}_{\lambda}(p,u;\gamma)=\sup_{(\alpha,\kappa):(\alpha,\kappa,\gamma)
\in K'} H_{\lambda}(p,u;\alpha,\kappa,\gamma)\]
is quadratic in $(p,u)$ in a neighborhood of $(p,u)=(0,0)$.
The condition
$\bA_{11}^{(\gamma)}-\bA_{12}^{(\gamma)}(\bA_{22}^{(\gamma)})^{-1}
\bA_{21}^{(\gamma)} \succ 0$ shows that $\bar{H}_{\lambda}(p,u;\gamma)$ is
strictly convex in $(p,u)$ near $(0,0)$, and is minimized at $(p,u)=(0,0)$.
We thus obtain for any fixed $\gamma$ with $|\gamma|<c$ that
\begin{align}
\inf_{(p,u) \in \reals^2} \sup_{(\alpha,\kappa):(\alpha,\kappa,\gamma) \in K'}
H_\lambda(p,u;\alpha,\kappa,\gamma)&=H_\lambda(0,0;0,0,\gamma)\nonumber\\
&=\gamma+\lambda^2(1-q_\star)
-\lambda^2\E\Big[\Big(\frac{1}{1-m^2}-\gamma\Big)^{-1}\Big].\label{eq:Hvalue}
\end{align}

Finally, observe that
\[H_\lambda(0,0;0,0,0)=\lambda^2(1-q_\star)-\lambda^2\E[1-m^2]=0.\]
The derivative in $\gamma$ of the right side of (\ref{eq:Hvalue})
may be evaluated inside the expectation, by a standard application of the 
dominated convergence theorem. Then
$\partial_\gamma H_\lambda(0,0;0,0,\gamma)|_{\gamma=0}=
1-\lambda^2\E[(1-m^2)^2]=1-\lambda^2(1-2q_\star+b_\star)$. Applying
$q_\star>b_\star$ and $1-q_\star<1/\lambda^2$ by Proposition \ref{prop:qstar}, we get
$\lambda^2(1-2q_\star+b_\star)<\lambda^2(1-q_\star)<1$, so
\[\partial_\gamma H_\lambda(0,0;0,0,\gamma)|_{\gamma=0}>0.\]
Then there exists $\gamma_\star>0$ sufficiently small such that, bounding the
supremum over $\gamma$ by the value at $\gamma=\gamma_\star$,
\[\inf_{(p,u) \in \reals^2} \sup_{(\alpha,\kappa,\gamma) \in K'}
H_\lambda(p,u;\alpha,\kappa,\gamma)
\geq H_\lambda(0,0;0,0,\gamma_\star)>0.\]
Identifying $t_0=H_\lambda(0,0;0,0,\gamma_\star)$ concludes the proof.
\end{proof}

We combine these results with the following minor extension of the analysis of
$\det \bH(\bbm)$ from \cite{fan2021tap}, to show Corollary
\ref{cor:strongconvexity}.

\begin{lemma}\label{lem:detHbound}
Fix any $\lambda>1$, and suppose $\bx=\ones$. For any $c>0$,
there exists $\eta>0$ depending on $(\lambda,c)$ such that for all large $n$,
\[\int_{\cD_\eta}
\E\Big[|\det \bH(\bbm)|^2\;\Big|\;\bg(\bbm)=\bzero\Big]^{1/2}
\;p_{\bg(\bbm)}(\bzero)\,\de \bbm \leq e^{cn}.\]
\end{lemma}
\begin{proof}
Define
\[A(\bbm)=\frac{1}{n}\bbm^\sT \arctanh \bbm,
\qquad E(\bbm)=-H(\bbm)-\frac{1}{2}A(\bbm)-\frac{\lambda^2}{4}[1-Q(\bbm)^2].\]
The functions $m \mapsto m^2,m,m\arctanh m,\sh(m)$ all belong
to the class (\ref{eq:W2class}). Thus, by Proposition \ref{prop:qstar},
the convergence $W(\hat{\mu}_\bbm,\mu_\star) \to 0$ implies
$(Q(\bbm),M(\bbm),A(\bbm),E(\bbm)) \to (q_\star,q_\star,\lambda^2
q_\star,e_\star)$. So for any $\delta>0$, there is a constant
$\eta=\eta(\delta)>0$ for which
\[\cD_\eta \subseteq
\cE_\delta=\Big\{\bbm \in (-1,1)^n:
\;|Q(\bbm)-q_\star|,|M(\bbm)-q_\star|,|A(\bbm)-\lambda^2 q_\star|,
|E(\bbm)-e_\star|<\delta\Big\}.\]

Let
\[L(\bbm)=\frac{\lambda^2[1-Q(\bbm)]^2}{2}+\frac{1}{n}\sum_{i=1}^n
\log \left(\frac{1}{1-m_i^2}\right).\]
Then a small modification of the proof of \cite[Proposition 3.2]{fan2021tap}
shows that for any fixed $\bbm \in \cE_\delta$,
\begin{equation}\label{eq:detsqbound}
\E\Big[|\det \bH(\bbm)|^2\;\Big|\;\bg(\bbm)=\bzero\Big]^{1/2}
\leq \exp(n \cdot L(\bbm)+Cn^{0.9})
\end{equation}
for all large $n$ and some constant $C>0$. Indeed, defining
$l(x)=\log |x-\mathbf{i} \cdot n^{-0.11}|=\operatorname{Re}
\log (x-\mathbf{i} \cdot n^{-0.11})$
and applying this spectrally to the matrix $\bH(\bbm)$
by the functional calculus, we have
\[\P\Big[\Tr l(\bH(\bbm)) \geq c_n+nt\;\Big|\;
\bg(\bbm)=\bzero\Big] \leq 2e^{-cn^{1.78}t^2},\quad
c_n=n \cdot L(\bbm)+Cn^{0.89}\]
as established in the proof of \cite[Proposition 3.2]{fan2021tap}. Then,
applying
\[\E[e^X] \leq e^t+\E[(e^X-e^t)\ones\{X \geq t\}]
=e^t+\int_t^\infty e^s \cdot \P[X>s]\de s\]
we obtain
\begin{align*}
\E\Big[|\det \bH(\bbm)|^2\;\Big|\;\bg(\bbm)=\bzero\Big]
&\leq \E\Big[e^{2\Tr l(\bH(\bbm))}\;\Big|\;\bg(\bbm)=\bzero\Big]\\
&\leq e^{2c_n+n^{0.89}}+\int_{2c_n+n^{0.89}}^\infty
e^s\cdot \P\Big[2\Tr l(\bH(\bbm))>s \;\Big|\;\bg(\bbm)=\bzero\Big]\de s\\
&\leq e^{2c_n+n^{0.89}}+\int_{2c_n+n^{0.89}}^\infty
e^s \cdot 2e^{-cn^{-0.22}(s/2-c_n)^2}\de s\\
&=e^{2c_n+n^{0.89}}+\int_{n^{0.89}/2}^\infty 4e^{2(t+c_n)-cn^{-0.22}t^2}\de t
\leq e^{2c_n+Cn^{0.89}}
\end{align*}
for all large $n$, which implies (\ref{eq:detsqbound}).
Then substituting (\ref{eq:detsqbound}) in place of the bound for $\E[|\det
\bH(\bbm)| \mid \bg(\bbm)=\bzero]$ in the proof of \cite[Theorem
1.1]{fan2021tap}, we obtain
\begin{align*}
&\limsup_{n \to \infty} \frac{1}{n}
\log \int_{\cE_\delta} \E\Big[|\det \bH(\bbm)|^2
\;\Big|\;\bg(\bbm)=\bzero\Big]^{1/2}\,p_{\bg(\bbm)}(\bzero)\,\de \bbm\\
&\hspace{1in}\leq \sup_{(q,\varphi,a,e) \in \reals^4:
|q-q_\star|,|\varphi-q_\star|,|a-\lambda^2 q_\star|,|e-e_\star|<\delta}
S_\star(q,\varphi,a,e)
\end{align*}
where the complexity functional $S_\star(q,\varphi,a,e)$ is as defined in
\cite[Eq.\ (1.10)]{fan2021tap}. This function $S_\star$ is continuous by
definition, and \cite[Proposition 5.2]{fan2021tap} shows (for any $\lambda>0$)
that $S_\star(q_\star,q_\star,\lambda^2 q_\star,e_\star)=0$. Thus, for any
$c>0$, there is $\delta>0$ small enough so that this supremum is less than $c$,
and the result follows.
\end{proof}

\begin{proof}[Proof of Corollary \ref{cor:strongconvexity}]
For $t>0$ to be determined, let $U \sim \Unif([-t,t])$ be independent of $\bW$.
Applying Lemma \ref{lem:KacRicebound} with $T=\cD_\eta$ and
$\ell(\bbm,\bW)=\ell_\eps^+(\bbm)$, we have
\begin{align*}
&\P\Big[\text{ there exist } \bbm \in \cD_\eta \text{ and }
\bu \in (-1,1)^n \cap \ball_{\sqrt{\eps n}}(\bbm):
\;\bg(\bbm)=\bzero \text{ and } \lambda_{\min}(\bH(\bu))<t\;\Big]\\
&=\P\Big[\text{ there exists } \bbm \in \cD_\eta:
\;\bg(\bbm)=\bzero \text{ and } \ell_\eps^+(\bbm)<t\;\Big]\\
&\leq \P\Big[\text{ there exists } \bbm \in \cD_\eta:
\;\bg(\bbm)=\bzero \text{ and } \ell_\eps^+(\bbm)+U<2t\;\Big]\\
&\leq \int_{\cD_\eta} \E\Big[|\det \bH(\bbm)| \cdot
\ones\{\ell_\eps^+(\bbm)+U<2t\}\;\Big|\;\bg(\bbm)=\bzero\Big]\;
p_{\bg(\bbm)}(\bzero)\,\de\bbm\\
&\leq \int_{\cD_\eta} \E\Big[|\det \bH(\bbm)| \cdot
\ones\{\ell_\eps^+(\bbm)<3t\}\;\Big|\;\bg(\bbm)=\bzero\Big]\;
p_{\bg(\bbm)}(\bzero)\,\de\bbm.
\end{align*}
Then by Cauchy-Schwarz,
\begin{align}
&\P\Big[\text{ there exist } \bbm \in \cD_\eta \text{ and }
\bu \in (-1,1)^n \cap \ball_{\sqrt{\eps n}}(\bbm):
\;\bg(\bbm)=\bzero \text{ and } \lambda_{\min}(\bH(\bu))<t\;\Big]\nonumber\\
&\leq \int_{\cD_\eta} \E\Big[|\det \bH(\bbm)|^2 
\;\Big|\;\bg(\bbm)=\bzero\Big]^{1/2}\;
p_{\bg(\bbm)}(\bzero)\,\de\bbm
\cdot \sup_{\bbm \in \cD_\eta} \P\Big[\ell_\eps^+(\bbm)<3t \;\Big|\;
\bg(\bbm)=\bzero\Big]^{1/2}.\label{eq:stronglyconcavebound}
\end{align}
Let us express the conditional law of
$\{-\bH(\bu):\bu \in (-1,1)^n\}|_{\bg(\bbm)=\bzero}$ as the
right side of (\ref{eq:condhess}), and write
$\tilde{\bW}=(\bZ+\bZ^\sT)/\sqrt{2n}$ where $\bZ \in \reals^{n \times n}$
has i.i.d.\ $\cN(0,1)$ entries. Then
\[\bZ \mapsto
\inf\Big\{\bv^\sT \bH(\bu)\bv:\;\bu \in (-1,1)^n \cap \ball_{\sqrt{\eps
n}}(\bbm),\;\bv \in \reals^n,\;\|\bv\|_2=1\Big\}\]
is $\lambda\sqrt{2/n}$-Lipschitz with respect to the Frobenius norm of $\bZ$.
Hence by the Gaussian concentration inequality, for any $s>0$,
\[\P\Big[\ell_\eps^+(\bbm)<\E[\ell_\eps^+(\bbm) \mid \bg(\bbm)=\bzero]-s
\;\Big|\;\bg(\bbm)=\bzero\Big] \leq e^{-s^2n/4\lambda^2}.\]
Recalling the constant $t_0>0$ from Lemma \ref{lem:gordon-variational}(b),
let us set $3t=s=t_0/3$. Then by Lemma \ref{lem:gordon-variational}(a--b),
for some sufficiently small $\lambda$-dependent constants $\eps,\eta,c_1>0$,
all $\bbm \in \cD_\eta$, and all large $n$,
\[\P\Big[\ell_\eps^+(\bbm)<3t \;\Big|\;\bg(\bbm)=\bzero\Big] \leq e^{-c_1n}.\]
Now applying Lemma \ref{lem:detHbound}, for this constant $c_1$,
some sufficiently small $\eta>0$, and all large $n$,
\[\int_{\cD_\eta} \E\Big[|\det \bH(\bbm)|^2
\;\Big|\;\bg(\bbm)=\bzero\Big]^{1/2}\; p_{\bg(\bbm)}(\bzero)\,\de\bbm
\leq e^{c_1 n/4}.\]
Applying these two bounds to (\ref{eq:stronglyconcavebound}) shows
(\ref{eq:stronglyconcave}).
\end{proof}

\begin{proof}[Proof of Theorem \ref{thm:local}(a--b)]
Fix $\eta,\eps,t>0$ small enough as described in
Corollary \ref{cor:strongconvexity}. Fix any $\iota>0$. We choose
$\delta=\delta(\lambda,\eta,\eps,t,\iota)$ small enough, to be determined,
so that the conclusions of Corollary \ref{cor:critexists} and
Lemma \ref{lem:localizetoDeta} hold for this $\delta$.
Then with probability approaching 1,
Corollary \ref{cor:critexists} establishes the existence of 
a local minimizer $\bbm_\star$ of $\cF_{\TAP}$ belonging to $\cB_\delta$,
such that $|\cF_{\TAP}(\bbm_\star)-e_\star|<\delta$.
Lemma \ref{lem:localizetoDeta} ensures that
this local minimizer $\bbm_\star$ belongs also to $\cD_\eta$, and
Corollary \ref{cor:strongconvexity} then ensures that for all $\bu \in
\ball_{\sqrt{\eps n}}(\bbm_\star)$, we have
$\lambda_{\min}(\nabla^2 \cF_{\TAP}(\bu))>t/n$
(where we recall that $\bH$ is the Hessian of $\cF_\TAP$ rescaled by $n$).

Let us now show that with probability approaching 1,
we may pick such a local minimizer $\bbm_\star$ to satisfy the Bayes optimality
condition (\ref{eq:bayesoptimal}).
Lemma \ref{lem:AMP}(b--d) shows that for sufficiently large $k$, 
with probability approaching 1, the AMP iterate $\bbm^k$ satisfies
$\bbm^k \in \cB_{\delta/2}$, $\cF_\TAP(\bbm^k)<e_\star+ \delta \wedge (c_0 \delta^2)$  with $c_0$ as in the proof of Corollary \ref{cor:critexists}, and also
\[\frac{1}{n^2}\|\widehat{\bX}_{\Bayes}-\bbm^k(\bbm^k)^\sT\|_\sF^2<\delta.\]

 Now we claim that we may pick the local minimizer $\bbm_\star$ in the proof of
Corollary \ref{cor:critexists} so that there is a path $\Gamma$
connecting $\bbm^k$ to $\bbm_\star$ for
which $\cF_\TAP(\bbm) \leq \cF_\TAP(\bbm^k)<e_\star+\delta$ for all $\bbm \in
\Gamma$. 
To see this, 
consider running gradient flow $\frac{\de}{\de t} \bbm = - \nabla \cF_{\TAP}(\bbm) $ initialized at $\bbm^k$. 
Because (1) $\cF_{\TAP}(\bbm)$ must be non-increasing along this trajectory, 
(2) by Eqs.~\eqref{eq:FTAPlowerbound1} and \eqref{eq:FTAPlowerbound2} $\cF_{\TAP}(\bbm)$ is larger on $\overline{\cB_{\delta}} \setminus \cB_{\delta/2}$ than it is at $\bbm^k$,
and (3) the gradient flow would need to, by the continuity of $\cF_{\TAP}$, pass through $\overline{\cB_{\delta}} \setminus \cB_{\delta/2}$ were it ever to leave $\cB_{\delta/2}$,
we can conclude that the gradient flow stays in $\cB_{\delta/2}$ for all time.
Moreover, 
by the Lojasiewicz Theorem \cite{lojasiewicz1982},
the gradient flow converges to a critical point $\bbm_\star$ of $\cF_{\TAP}$. 
By Corollary \ref{cor:critexists},
this critical point must be a local minimizer, and the gradient flow guarantees the the property $\cF_\TAP(\bbm) \leq \cF_\TAP(\bbm^k)<e_\star+\delta$ on the full path.

Now, choosing $\delta<\eps t/6$, if
$\bbm^k \notin \ball_{\sqrt{6\delta n/t}}(\bbm_\star)$, then at the point
$\bbm \in \Gamma$ where $\Gamma$ crosses the boundary of
$\ball_{\sqrt{6\delta n/t}}(\bbm_\star)$, the strong convexity of $\cF_{\TAP}$
on $\ball_{\sqrt{\eps n}}(\bbm_\star)$ implies that
\[\cF_{\TAP}(\bbm)-\cF_{\TAP}(\bbm_\star)
\geq \frac{t}{2n} \cdot (\sqrt{6\delta n/t})^2=3\delta.\]
But this is a contradiction because $e_\star-\delta<
\cF_{\TAP}(\bbm_\star) \leq \cF_{\TAP}(\bbm)<e_\star+\delta$.
Thus, for this choice of $\bbm_\star$, we must have
$\|\bbm_\star-\bbm^k\|_2 \leq \sqrt{6\delta n/t}$. Then with probability
approaching 1,
\begin{align*}
\frac{1}{n^2}\|\widehat{\bX}_{\Bayes}
-\bbm_\star\bbm_\star^\sT\|_\sF^2
&\leq \frac{2}{n^2}\|\widehat{\bX}_{\Bayes}
-\bbm^k(\bbm^k)^\sT\|_\sF^2
+\frac{2}{n^2}\|\bbm^k(\bbm^k)^\sT-\bbm_\star\bbm_\star^\sT\|_\sF^2\\
&<2\delta+24\delta/t
\end{align*}
for all large $n$. Choosing $\delta$ small enough so that
$2\delta+24\delta/4<\iota$, we obtain that $\bbm_\star$ satisfies
(\ref{eq:bayesoptimal}).

Finally, for sufficiently small $\iota>0$, uniqueness of the critical point
$\bbm_\star$ up to sign
follows from the strong convexity of $\cF_{\TAP}$ near $\bbm_\star$: If
$\bbm,\bbm'$ both satisfy
$n^{-2}\|\widehat{\bX}_{\Bayes}-\bbm\bbm^\sT\|_\sF^2<\iota$, then
\begin{align}
\min\left(\frac{1}{n}\|\bbm+\bbm'\|_2^2,
\frac{1}{n}\|\bbm-\bbm'\|_2^2\right)^2
&\leq \frac{1}{n^2}\|\bbm+\bbm'\|_2^2
\|\bbm-\bbm'\|_2^2\nonumber\\
&=\frac{1}{n^2}\left(\|\bbm\|_2^4+\|\bbm'\|_2^4
+2\|\bbm\|_2^2\|\bbm'\|_2^2
-4(\bbm^\sT \bbm')^2\right)\nonumber\\
&\leq
\frac{2}{n^2}\Big(\|\bbm\|_2^4+\|\bbm'\|_2^4
-2(\bbm^\sT \bbm')^2\Big)\nonumber\\
&=\frac{2}{n^2}\|\bbm\bbm^\sT-\bbm'(\bbm')^\sT\|_\sF^2
<8\iota.\label{eq:pmbound}
\end{align}
Thus either $\bbm_\star' \in \ball_{(8\iota)^{1/4}\sqrt{n}}(\bbm_\star)$
or $\bbm_\star' \in \ball_{(8\iota)^{1/4}\sqrt{n}}(-\bbm_\star)$.
When $(8\iota)^{1/4}<\sqrt{\eps}$, the strong convexity of $\cF_{\TAP}$ on
$\ball_{\sqrt{\eps n}}(\bbm_\star)$ ensures in the first case that
$\bbm_\star'=\bbm_\star$. By sign symmetry, $\cF_{\TAP}$ is also strongly convex
on $\ball_{\sqrt{\eps n}}(-\bbm_\star)$ with a local minimizer at $-\bbm_\star$,
so $\bbm_\star'=-\bbm_\star$ in the latter case.
\end{proof}

\subsection{Proofs for Section \ref{sec:AMP-stab}}\label{appendix:AMPstabproofs}

\begin{proof}[Proof of Lemma \ref{lem:stabilitycondition}]
Define
\[
	\begin{aligned}
		\bX(\bbm) 
			=&~ \diag(1 -\bbm^2)^{1/2} \cdot [ \lambda \bY  + 2
\lambda^2 \bbm \bbm^\sT / n ] \cdot \diag(1 -\bbm^2)^{1/2}, \\
		\bD(\bbm) 
			=&~ \lambda^2  [1 - Q(\bbm)] \cdot \diag(1 -\bbm^2).
	\end{aligned}
\]
Then for any $\mu \in \C$, note that
\begin{align*}
	\det(\bB(\bbm)-\mu\id)&=\det \begin{bmatrix} \diag(1-\bbm^2)^{1/2}
\cdot \bX(\bbm) \cdot \diag(1-\bbm^2)^{-1/2} - \mu \id &
- \bD(\bbm) \\ \id & - \mu \id  \end{bmatrix} \\
		&= \det\Big( \mu^2 \id - \mu \cdot \diag(1-\bbm^2)^{1/2} \cdot
\bX(\bbm) \cdot \diag(1-\bbm^2)^{-1/2} + \bD(\bbm) \Big)\\
&=\det\Big( \mu^2 \id - \mu \bX(\bbm) + \bD(\bbm) \Big).
\end{align*}
Then $\mu$ is an eigenvalue of $\bB(\bbm)$ if and only if $\mu^2 \id - \mu
\bX(\bbm) + \bD(\bbm)$ is singular. (Note that this is equivalent to singularity
of the matrix (\ref{eq:betheHessiancomplex}) as discussed in the main text.)

For any $\xi \in \C$, define the matrix 
\[
	\bK(\xi, \mu, \bbm) = \xi \mu^2 \id - \xi \mu \bX(\bbm) + \xi \bD(\bbm).
\]
Then for any $\bv \in \C^n$, we have
\[
	\bv ^*  \bK(\xi, \mu, \bbm) \bv = \bv^*[\Re \bK(\xi, \mu, \bbm)] \bv +
\bi \cdot \bv^* [\Im \bK(\xi, \mu, \bbm)] \bv,
\]
where $\Re$ and $\Im$ denote the entry-wise real and imaginary parts. 
Since $\id$, $\bX(\bbm)$, and $\bD(\bbm)$ are all real and symmetric, we have
that $\Re \bK (\xi, \mu, \bbm)$ and $\Im \bK (\xi, \mu, \bbm)$ are also real and
symmetric, so $\bv^*[\Re \bK(\xi, \mu, \bbm)] \bv$ and $\bv^*[\Im \bK(\xi, \mu,
\bbm)] \bv$ are both real. If there exists some $\xi \in \C$ for which
$\Re \bK (\xi, \mu, \bbm)$ is non-singular, then this implies the real part of
$\bv^* \bK(\xi,\mu,\bbm) \bv$ is non-zero for any non-zero vector
$\bv \in \C^n$, so $\mu$ is not an eigenvalue of $\bB(\bbm)$.

By Proposition \ref{prop:qstar}, $\lambda^2(1-q_\star+\delta)<1$ for some
$\delta>0$. Let us take $r_0=\sqrt{\lambda^2(1-q_\star+\delta)}$, and suppose
that (\ref{eq:betheHessian}) holds for $r \in (r_0,1)$. We show that for any
$\mu \in \C$ with $|\mu| \geq r$, we may pick $\xi$ to ensure that
$\Re \bK (\xi, \mu, \bbm) \succ 0$, and hence $\Re \bK (\xi, \mu, \bbm)$
is non-singular. There are four cases.
\begin{description}

	\item[Case 1 (positive imaginary part)] Suppose $\mu=\rho e^{\bi \vphi}$
with $\rho \geq r$ and $\vphi \in (0,\pi)$. Take $\xi = -\bi e^{-\bi\vphi}$.
	Then
	\begin{align*}
		\Re \bK(-\bi e^{-\bi \vphi}, \rho e^{\bi\vphi}, \bbm)
			&=\rho^2 \cos(3\pi/2 + \vphi)\id
			+\cos(3\pi/2 - \vphi) \bD(\bbm) \\
			&=\cos(3\pi/2 + \vphi)\left(\rho^2 \id - 
				\lambda^2[1-Q(\bbm)]\diag(1-\bbm^2) \right).
	\end{align*}
Applying $\diag(1-\bbm^2) \preceq \id$,
$\lambda^2[1-Q(\bbm)]<\lambda^2(1-q_\star+\delta)<r_0^2<\rho^2$,
and $\cos(3\pi/2 + \vphi)>0$,
this ensures $\Re \bK (-\bi e^{-\bi \varphi}, \rho e^{\bi \varphi}, \bbm) \succ \bzero$.
	\item[Case 2 (negative imaginary part)] Suppose $\mu = \rho e^{\bi
\vphi}$ with $\rho \geq r$ and $\vphi \in (\pi,2\pi)$. Take $\xi = \bi e^{-\bi
\vphi}$. Then
\[\Re \bK (\bi e^{-\bi \vphi}, \rho e^{\bi\vphi}, \bbm)
=\cos(\pi/2+\varphi)\big(\rho^2\id-\lambda^2[1-Q(\bbm)]\diag(1-\bbm^2)\big).\]
By the same argument as in case 1, $\Re \bK(\bi e^{-\bi \varphi}, \rho e^{\bi
\varphi}, \bbm) \succ \bzero$.
	\item[Case 3 (real and positive)] Suppose $\mu = \rho \geq r$.
	Take $\xi=1$, and note that
	\begin{align*}
		&\diag(1-\bbm^2)^{-1/2\;} \Re \bK(1,\rho,\bbm)\;\diag(1-\bbm^2)^{-1/2}
		\\
		&\qquad\qquad\qquad\qquad= 
			\rho\bigg[{-}(\lambda \bY + 2\lambda^2 \bbm\bbm^\sT/n) + 
\rho^{-1}\lambda^2[1-Q(\bbm)]\id + \rho\; \diag\Big(\frac1{1-\bbm^2}\Big)\bigg].
	\end{align*}
When $\rho=r$, this is positive-definite by assumption. For any $a,b>0$,
the function $\rho \mapsto a/\rho+b\rho$ is increasing for
$\rho \geq \sqrt{a/b}$.
Then, applying $\lambda^2[1-Q(\bbm)]<\lambda^2(1-q_\star+\delta)<r^2$
and $\diag(1/(1-\bbm^2)) \succeq \id$, the matrix inside the parenthesis is
increasing in $\rho$ in the positive-definite ordering, for $\rho \geq r$. Hence this
matrix is positive definite for all $\rho \geq r$, implying also
$\Re \bK(1,\rho,\bbm) \succ 0$.
	\item[Case 4 (real and negative)] Suppose $\mu = -\rho \leq -r$.
	Take $\xi=1$, and note that
	\begin{align*}
		&\diag(1-\bbm^2)^{-1/2\;} \Re \bK(1,-\rho,\bbm)\;\diag(1-\bbm^2)^{-1/2}
		\\
		&\qquad\qquad\qquad\qquad= 
			\rho\bigg[(\lambda \bY + 2\lambda^2 \bbm\bbm^\sT/n) +
\rho^{-1}\lambda^2[1-Q(\bbm)]\id + \rho\; \diag\Big(\frac1{1-\bbm^2}\Big)\bigg].
	\end{align*}
When $\rho=r$, this is again positive definite by assumption, so $\Re
\bK(1,-\rho,\bbm)\succ 0$ for all $\rho \geq r$ by the same argument as in case
3.
\end{description}

Combining these cases, $\bB(\bbm)$ does not have any eigenvalue $\mu \in \C$
with $|\mu| \geq r$, as desired.
\end{proof}

\begin{proof}[Proof of Lemma \ref{lem:gordon-variational-gen}(b)]
The argument is similar to that of Lemma \ref{lem:gordon-variational}(b): Let us
write as shorthand $H_\lambda=H_\lambda^-$ and $\E$ for $\E_{m \sim \mu_\star}$.
Then
\begin{align*}
&~H_\lambda(p,u;\alpha,\kappa,\gamma)\\
&=2\lambda^2 p^2+\lambda^2 u^2-2\lambda^2(1-q_\star)p^2/q_\star
-\alpha u-\kappa p+\gamma+\lambda^2(1-q_\star)\\
&\quad -\E\Big[\Big(4\lambda^2(1-p^2/q_\star)+(2z(m)p/q_\star
+\alpha+\kappa m)^2\Big)\Big/\Big(\frac{4}{1-m^2}-4\gamma\Big)\Big]\\
&=(p,u,\kappa,\alpha)^\sT \begin{pmatrix} \bA_{11}^{(-,\gamma)}
& \bA_{12}^{(-,\gamma)} \\ \bA_{21}^{(-,\gamma)} & \bA_{22}^{(-,\gamma)}
\end{pmatrix}(p,u,\kappa,\alpha)+\gamma+\lambda^2(1-q_\star)
-\lambda^2\E\Big[\Big(\frac{1}{1-m^2}-\gamma\Big)^{-1}\Big]
\end{align*}
where, specializing to $\gamma=0$, we have
\begin{align*}
\bA_{11}^{(-,0)} &=
\begin{pmatrix} 
2\lambda^2-\lambda^2\frac{1-q_\star}{q_\star} 
-\frac{\lambda^4(1-q_\star)^2}{q_\star^2}(q_\star-b_\star)
-\frac{2\lambda^3(1-q_\star)}{q_\star^{3/2}} \E[Gm(1-m^2)] 
-\frac{\lambda^2}{q_\star}\E[G^2(1-m)^2] & 0 \\ 0 & \lambda^2
\end{pmatrix} \\
\bA_{12}^{(-,0)}&=(\bA_{21}^{(-,0)})^\sT 
=
\begin{pmatrix} 
-\frac{1}{2}-\frac{\lambda^2(1-q_\star)}{2q_\star}(q_\star-b_\star)
-\frac{\lambda}{2q_\star^{1/2}}\E[Gm(1-m^2)] &
-\frac{\lambda^2(1-q_\star)}{2q_\star}(q_\star-b_\star)
+\frac{\lambda}{2q_\star^{1/2}}\E[Gm^2] \\
0 & -\frac{1}{2} \end{pmatrix} \\
\bA_{22}^{(-,0)} 
&= -\frac{1}{4}
\begin{pmatrix} q_\star-b_\star & q_\star-b_\star \\
q_\star-b_\star & 1-q_\star \end{pmatrix}.
\end{align*}
Recalling $\bA_{11},\bA_{12},\bA_{21},\bA_{21}$ from (\ref{eq:A22}), we then
may check that
\begin{align*}
\bA_{22}^{(-,0)}&=\frac{1}{4}\bA_{22}\\
\bA_{12}^{(-,0)}&=\frac{1}{2}\left(
\bA_{12}+\begin{pmatrix} \frac{\lambda^2(1-q_\star)}{q_\star}
& 0 \\ 0 & 0 \end{pmatrix} \bA_{22}-2\id\right)\\
\bA_{11}^{(-,0)}&=\bA_{11}+\bA_{12} 
\begin{pmatrix} \frac{\lambda^2(1-q_\star)}{q_\star} & 0 \\ 0 & 0 \end{pmatrix}
+\begin{pmatrix} \frac{\lambda^2(1-q_\star)}{q_\star} & 0 \\ 0 & 0 \end{pmatrix}
\bA_{21}
+\begin{pmatrix} \frac{\lambda^2(1-q_\star)}{q_\star} & 0 \\ 0 & 0 \end{pmatrix}
\bA_{22}
\begin{pmatrix} \frac{\lambda^2(1-q_\star)}{q_\star} & 0 \\ 0 & 0
\end{pmatrix}\\
&\hspace{1in}
+\begin{pmatrix} 4\lambda^2-4\lambda^2\frac{1-q_\star}{q_\star} & 0 \\ 0 &
2\lambda^2 \end{pmatrix}.
\end{align*}
Thus $\bA_{22}^{(-,0)} \prec 0$, and
\begin{align*}
&~\bA_{11}^{(-,0)}-\bA_{12}^{(-,0)}(\bA_{22}^{(-,0)})^{-1}\bA_{21}^{(-,0)}\\
&=\Big(\bA_{11}-\bA_{12}\bA_{22}^{-1}\bA_{21}\Big)
+2\bA_{12}\bA_{22}^{-1}+2\bA_{22}^{-1}\bA_{21}-4\bA_{22}^{-1}
+\begin{pmatrix} 4\lambda^2 & 0 \\ 0 & 2\lambda^2 \end{pmatrix}\\
&=\bA_{11}
-(\bA_{12}-2\id)\bA_{22}^{-1}(\bA_{21}-2\id)
+\begin{pmatrix} 4\lambda^2 & 0 \\ 0 & 2\lambda^2 \end{pmatrix}
\end{align*}
Computing explicitly $\bA_{22}^{-1}$, we obtain after some simplification
\[\bA_{11}^{(-,0)}-\bA_{12}^{(-,0)}(\bA_{22}^{(-,0)})^{-1}\bA_{21}^{(-,0)}
=\begin{pmatrix} c_1 & -c_2 \\ -c_2 & c_2 \end{pmatrix}\]
where now
\[
\begin{aligned}
c_1=&~\frac{(1-q_\star)+2\lambda^2(1-2q_\star+b_\star)^2
+\lambda^4(1-2q_\star+b_\star)^3}{(1-2q_\star+b_\star)(q_\star-b_\star)}
-\lambda^4(1-2q_\star+b_\star),\\
c_2=&~\frac{1}{1-2q_\star+b_\star}+\lambda^2
\end{aligned}\]
Applying again (\ref{eq:basicinequalities}), together with
\[1-\lambda^2(q_\star-b_\star)+\lambda^2(1-2q_\star+b_\star)
>1-\lambda^2(1-q_\star)>0,\]
we get $c_2>0$ and
\begin{align*}
&c_1-(-c_2)c_2^{-1}(-c_2)=c_1-c_2\\
&=-\lambda^2\Big(1+\lambda^2(1-2q_\star+b_\star)\Big)
+\frac{1}{q_\star-b_\star}\Big(1+\lambda^2(1-2q_\star+b_\star)\Big)^2\\
&=\frac{1+\lambda^2(1-2q_\star+b_\star)}{q_\star-b_\star}
\cdot \Big(-\lambda^2(q_\star-b_\star)+1+\lambda^2(1-2q_\star+b_\star)\Big)>0.
\end{align*}
Thus $\bA_{11}^{(-,0)}-\bA_{12}^{(-,0)}(\bA_{22}^{(-,0)})^{-1}\bA_{21}^{(-,0)}
\succ 0$. This implies, as in the proof of Lemma \ref{lem:gordon-variational}(b),
that for all $|\gamma|<c$ small enough,
\begin{align*}
\inf_{(p,u) \in \reals^2} \sup_{(\alpha,\kappa):(\alpha,\kappa,\gamma) \in K'}
H_\lambda(p,u;\alpha,\kappa,\gamma)
&=H_\lambda(0,0;0,0,\gamma)\\
&=\gamma+\lambda^2(1-q_\star)
-\lambda^2\E\Big[\Big(\frac{1}{1-m^2}-\gamma\Big)^{-1}\Big].
\end{align*}
The conclusion then follows as in the proof of
Lemma \ref{lem:gordon-variational}(b).
\end{proof}

\begin{corollary}\label{cor:AMPstability}
Fix any $\lambda>1$, and suppose $\bx=\ones$. Then there exist
$\lambda$-dependent constants $\eps,\eta,t,c>0$ such that, for all large $n$,
\begin{align*}
\P\Big[\;\text{ there exist } \bbm \in \cD_\eta \text{ and }
\bu \in (-1,1)^n &\cap \ball_{\sqrt{\eps n}}(\bbm):\bg(\bbm)=\bzero \text{ and } \lambda_{\min}(\bH^-(\bu))<t\;\Big]<e^{-cn}.
\end{align*}
\end{corollary}
\begin{proof}
The proof is identical to that of Corollary \ref{cor:strongconvexity},
applying Lemma \ref{lem:gordon-variational-gen} in place of Lemma
\ref{lem:gordon-variational}.
\end{proof}

\begin{proof}[Proof of Theorem \ref{thm:local}(c)]
Fix $\eta,\delta>0$ small enough to satisfy Lemma
\ref{lem:stabilitycondition} and
Corollaries \ref{cor:strongconvexity} and \ref{cor:AMPstability}.
Let $\bbm_\star$ be the local minimizer
identified in Theorem \ref{thm:local}(a), which belongs to
$\cB_\delta \cap \cD_\eta$.
Corollaries \ref{cor:strongconvexity} and \ref{cor:AMPstability} imply that for
some $t>0$, with probability approaching 1,
both $\bH^+(\bbm_\star) \succ t\id$ and $\bH^-(\bbm_\star) \succ t\id$.
Consider now the matrices
\[\pm r^{-1}(\lambda \bY + 2\lambda^2 \bbm_\star\bbm_\star^\sT/n) -
r^{-2}\lambda^2[1-Q(\bbm_\star)]\id
- \; \diag\Big(\frac1{1-\bbm_\star^2}\Big),\]
where $r=1$ corresponds to $\bH^{\pm}(\bbm_\star)$.
We have $\|\lambda^2[1-Q(\bbm_\star)]\id\|_{\op} \leq \lambda^2$, and also
$\|\lambda \bY + 2\lambda^2 \bbm_\star\bbm_\star^\sT/n\|_{\op}
<3\lambda+3\lambda^2$ on the event of probability $1-e^{-cn}$ where
$\|\bW\|_{\op}<3$. Thus, on this event and for some constant
$r(t,\lambda) \in (0,1)$, the above matrices
must also be positive definite by continuity for all $r \in (r(t,\lambda),1)$.
Multiplying by $r^2$ and recalling $r_0$ from
Lemma \ref{lem:stabilitycondition}, this ensures that for some $r \in
(0,1)$ with $r>\max(r(t,\lambda),r_0)$, almost surely for all large $n$,
the matrices (\ref{eq:betheHessian}) are positive definite at
$\bbm_\star$. Then by Lemma \ref{lem:stabilitycondition},
\[\rho(\de T_{\AMP}(\bbm_\star,\bbm_\star))=\rho(\bB(\bbm_\star))<r<1.\]
\end{proof}

\section{Proofs for algorithm convergence}\label{appendix:global}

\subsection{Analysis of TAP landscape}\label{appendix:TAP_global_landscape}

We prove Corollary \ref{cor:globallandscape} and Lemma
\ref{lem:quantitativelandscape} on properties of the global
landscape of $\cF_\TAP$ when $\lambda>\lambda_0$.

\begin{proof}[Proof of Corollary \ref{cor:globallandscape}]
Fix any $\iota>0$. By \cite[Theorem 1.2]{fan2021tap}, with probability 
approaching 1, there exists a critical point $\bbm \in \cS$ of $\cF_\TAP$, and
furthermore all critical points $\bbm \in \cS$ satisfy
\[\frac{1}{n^2}\|\bbm\bbm^\top-\widehat{\bX}_\Bayes\|_\sF^2<\iota.\]
For small enough $\iota>0$, Theorem \ref{thm:local}(a) then implies that the
critical point $\bbm \in \cS$ is unique up to sign, so it must be the global
minimizer of $\cF_\TAP$ by the definition of the set $\cS$.
\end{proof}

We next show the following strengthened version of \cite[Lemma C.3]{fan2021tap}.

\begin{lemma}\label{lemma:gradientbound}
Fix any integer $a \geq 0$, and set $q=1-\lambda^{-a}$ and
$t=\lambda^{2-a}$. Suppose $\bx=\ones$. For a constant $\lambda_0(a)>0$, if
$\lambda>\lambda_0(a)$, then for some $(a,\lambda)$-dependent constants $C,c>0$,
with probability at least $1-Ce^{-cn}$,
\begin{equation}\label{eq:gradmlarge}
\Big\{\bbm \in (-1,1)^n:\;M(\bbm)+Q(\bbm)>1.01 \text{ and }
n \cdot \|\nabla \cF_\TAP(\bbm)\|_2^2<t \Big\} \subseteq \cM_q.
\end{equation}
\end{lemma}
\begin{proof}
We induct on $a$. The base case $a=0$ is trivial, because $M(\bbm)+Q(\bbm)>1.01$ implies $M(\bbm)>0$.

Suppose by induction that the statement holds for an integer $a \geq 0$.
Let $\bg(\bbm)=n \cdot \nabla \cF_\TAP(\bbm)$ be the normalized gradient from
(\ref{eq:grad}), and denote its coordinates as $g_i(\bbm)$. Fix $\delta=0.001$,
and consider
\[\cS_1=\Big\{ i \in \{1,\ldots,n\}
: \vert g_i(\bbm) \vert < \delta\lambda^2 \Big\}.\]
If $\bbm$ satisfies $n \cdot \|\nabla \cF_\TAP(\bbm)\|_2^2
=n^{-1}\sum_i g_i(\bbm)^2 < \lambda^{2-a}$, then by Markov's inequality,
\[|\cS_1|> n\left(1 - \frac{1}{\delta^2\lambda^{2+a}}\right).\]

Next, write $\bW \sim \GOE(n)$ as $\bW = (\bZ + \bZ^\sT)/\sqrt{2n}$, where
$Z_{ij} \overset{iid}{\sim} \cN(0, 1)$. Let $\bz_i$ be the $i^\text{th}$ column
of $\bZ$. Note that if $\bbm \in \cM_{1-\lambda^{-a}}$, then
\[\| \bbm - \ones \|_2/\sqrt n \le \sqrt{2 - 2 M(\bbm)} \le
\sqrt{2}\lambda^{-a/2}.\]
Then applying \cite[Lemma C.2]{fan2021tap} with $\delta\lambda^{1+a/2}/2$ in
place of $\lambda$, there are absolute constants $C_0,C,c>0$
(independent of $\lambda$) such that
\begin{equation}\label{eqn:Fraction_difference_of_m_ones}
\P\left[ \sup_{ \bbm \in \cM_{1-\lambda^{-a}}}
\frac{1}{n} \sum_{i=1}^n \ones\big\{ \vert
\langle \bz_i,\bbm - \ones \rangle \vert \ge \delta\lambda\sqrt{n/2}\big\}
\ge \frac{C_0}{\lambda^{2+a}} \right] \le Ce^{-cn}.
\end{equation}
By the Chernoff bound, for $\lambda$ large enough
such that $\P_{Z \sim \cN(0,1)}[|Z| \geq \delta
\lambda/\sqrt{2}]<C_0/(2\lambda^{2+a})$, we also have
\begin{equation}\label{eqn:Infinity_norm_bound_W_ones}
\P\left[ \frac{1}{n} \sum_{i=1}^n \ones\{ \vert \langle \bz_i,\ones \rangle
\vert \ge \delta \lambda\sqrt{n/2}\} \ge \frac{C_0}{\lambda^{2+a}}
\right] \le Ce^{-cn}.
\end{equation}
Then, applying these bounds also for the rows of $\bZ^\top$ and defining
\[\cS_2(\bbm)=\Big\{ i \in \{1,\ldots,n\}: \vert[
\bW (\bbm - \ones)]_i \vert \le \delta\lambda \text{ and }
\vert [ \bW \ones]_i \vert \le \delta \lambda \Big\},\]
we obtain
\begin{equation}\label{eq:S3large}
\P\left[ \inf_{\bbm \in \cM_{1-\lambda^{-a}}} \vert \cS_2(\bbm) \vert \ge n (1 -
4C_0/\lambda^{2+a})\right] \ge 1 - 4Ce^{-cn}.
\end{equation}

Consider now the event where (\ref{eq:S3large}) holds and where
(\ref{eq:gradmlarge}) holds for $a$. This has
probability at least $1-C'e^{-c'n}$ for some $(a,\lambda)$-dependent constants
$C',c'>0$, by the induction hypothesis. On this event, for any
$\bbm \in (-1,1)^n$ belonging to the left side of (\ref{eq:gradmlarge}),
the above shows $|\cS_1 \cap \cS_2(\bbm)| \geq
n(1-C_1/\lambda^{2+a})$ where $C_1=4C_0+1/\delta^2$. For this $\bbm$ and any
index $i \in \cS_1 \cap \cS_2(\bbm)$, writing
\[g_i(\bbm)=-\lambda^2 M(\bbm)-\lambda \cdot [\bW \ones]_i
-\lambda \cdot [\bW (\bbm-\ones)]_i+\arctanh(m_i)+\lambda^2[1-Q(\bbm)]m_i,\]
we have
\begin{align*}
m_i &> \tanh(\lambda^2 M(\bbm) + \lambda\cdot [\bW \ones]_i + \lambda \cdot [\bW(\bbm -
\ones)]_i - \lambda^2 [1 - Q(\bbm)] m_i - \delta\lambda^2)\\
&\ge \tanh((M(\bbm) + Q(\bbm) - 1) \lambda^2 - 3\delta\lambda^2)
\geq \tanh(7\delta\lambda^2),
\end{align*}
where the last step uses $M(\bbm)+Q(\bbm)>1.01$ and $\delta=0.001$. For
the coordinates $i \notin (\cS_1 \cap \cS_2(\bbm))$, we may
apply the trivial bound $m_i \ge -1$. Then on the above event,
for all sufficiently large $\lambda$,
\[M(\bbm)=\< \bbm, \ones\> /n > (1 - C_1/\lambda^{2+a}) \cdot \tanh(7\delta
\lambda^2) - C_1/\lambda^{2+a}> 1 - 3C_1/\lambda^{2+a}>1-1/\lambda^{1+a}.\]
So the left side of (\ref{eq:gradmlarge}) is contained in
$\cM_{1-\lambda^{-(a+1)}}$. This completes the induction.
\end{proof}

\begin{lemma}\label{lem:sandwich_operator_bound}
For $\bt \in [0,1]^n$ and $\eps>0$, denote by $S(\bt,\eps)$ the
subset of indices $i \in \{1,\ldots,n\}$ for which $t_i \geq \eps$. Then
there exist universal constants $C,C',c>0$ such that for $\bW \sim \GOE(n)$ and
any $\eps>0$ and $0<s<1$,
\begin{equation}\label{eqn:sandwich_W_operator_bound}
\P\left[\mathop{\sup_{\bt_1, \bt_2 \in [0, 1]^n:}}_{\max(\vert S(\bt_1, \eps) \vert,\vert
S(\bt_2, \eps) \vert) \le ns} \| \diag(\bt_1) \bW \diag(\bt_2) \|_{\op} \ge
C'\cdot (\eps + \sqrt{s \log(e/s)})   \right] \le C e^{- c s n}. 
\end{equation}
\end{lemma}

\begin{proof}[Proof of Lemma \ref{lem:sandwich_operator_bound}]~
Denote $\proj_S: \R^n \to \R^n$ to be the projection operator onto the subspace
associated to $S \subseteq \{1,\ldots,n\}$, and let $\proj_S^\perp=\id-\proj_S$.
Then
\[
\begin{aligned}
\| \diag(\bt_1) \bW \diag(\bt_2) \|_{\op} \le&~ \| \diag(\bt_1) \proj_{S(\bt_1,
\eps)} \bW \proj_{S(\bt_2, \eps)} \diag(\bt_2) \|_{\op} \\
&~ +  \| \diag(\bt_1)
\proj_{S(\bt_1, \eps)} \bW \proj_{S(\bt_2, \eps)}^\perp \diag(\bt_2) \|_{\op} \\
&~ + \| \diag(\bt_1) \proj_{S(\bt_1, \eps)}^\perp \bW \proj_{S(\bt_2, \eps)}
\diag(\bt_2) \|_{\op} \\
&~ +  \| \diag(\bt_1) \proj_{S(\bt_1, \eps)}^\perp \bW \proj_{S(\bt_2, \eps)}^\perp \diag(\bt_2) \|_{\op} \\
\le&~ \| \proj_{S(\bt_1, \eps)} \bW \proj_{S(\bt_2, \eps)} \|_{\op} + 3 \eps \|
\bW \|_{\op}.
\end{aligned}
\]
We have $\|\bW\|_\op<3$ with probability $1-Ce^{-cn}$. Thus it suffices to show
that 
\begin{equation}\label{eq:sparse-op-norm-bound}
\P\left[\sup_{\vert S_1\vert, \vert S_2 \vert \le n \cdot s} \| \bW_{S_1 S_2}
\|_{\op} \ge C \sqrt{s \log (e/s)} \right] \le e^{- c s n},
\end{equation}
where $\bW_{S_1S_2}$ is the submatrix of rows in $S_1$ and columns in $S_2$.
Let $\cB_s = \{ \bv \in \R^n :
\|\bv\|_2 \leq 1, \|\bv\|_0 \leq ns\}$, and note that
\[
	\sup_{\vert S_1\vert, \vert S_2 \vert \le n \cdot s} \| \bW_{S_1 S_2} \|_{\op} 
		=\sup_{\bv,\bv' \in \cB_s} \bv^\top \bW \bv'.
\]
For any support set $S$ of size $\lfloor ns \rfloor$,
the $1/4$-covering number of $\{\bv \in \R^n : \|\bv\|_2 \leq 1,\, \supp(\bv)
\subseteq S\}$ is no larger than $(1 + 1/8)^{ns}/(1/8)^{ns} = 9^{ns}$.
Then applying $\binom{n}{\lfloor ns \rfloor} \leq (e/s)^{ns}$,
the 1/4-covering number of $\cB_s$ is no larger than $(9e/s)^{ns}$.
Letting $\cN$ be a $1/4$-cover of $\cB_s$,
\begin{align*}
	\sup_{\bv,\bv' \in \cB_s} \bv^\top \bW \bv'
	&\leq 
	\sup_{\bv,\bv' \in \cN}
		\sup_{\delta\bv,\delta\bv' \in \cB_s/4} (\bv+\delta\bv)^\top \bW
(\bv'+\delta\bv')\\
	&\leq 
	\sup_{\bv,\bv' \in \cN} \bv^\top \bW \bv'
		+ \frac9{16} \sup_{\bv,\bv' \in \cB_s} \bv^\top \bW \bv',
\end{align*}
hence 
\[
	\sup_{\bv,\bv' \in \cB_s} \bv^\top \bW \bv'
	\leq 
	\frac{16}{7}
	\sup_{\bv,\bv' \in \cN} \bv^\top \bW \bv'.
\]
For each $\bv,\bv' \in \cN$,
a Gaussian tail bound yields
$\P(\bv^\top \bW \bv' \geq t)\leq e^{-nt^2/4}$.
The result then follows by choosing $t = \sqrt{8 s \log(9e/s)}$
and taking a union bound over $\cN$.
\end{proof}

\begin{proof}[Proof of Lemma \ref{lem:quantitativelandscape}]
For part (a), \cite[Lemma C.1]{fan2021tap} implies that with probability
$1-Ce^{-cn}$, every point $\bbm \in \cS$ satisfies $Q(\bbm) \geq M(\bbm)^2 \geq
(1/3-6/\lambda-4/\lambda^2)^{1/2}$. For $M(\bbm)>0$, this implies
$Q(\bbm)+M(\bbm)>1.01$. Then by Lemma \ref{lemma:gradientbound}, with
probability at least $1-Ce^{-cn}$, the guarantee of (a) holds for all
$\bbm \in \cS \setminus \cM_q$ where $M(\bbm)>0$. The guarantee then also holds
for $M(\bbm)<0$ by the sign symmetry $\cF_\TAP(\bbm)=\cF_\TAP(-\bbm)$.

For part (b), recall the form (\ref{eq:hess}) for
$\bH(\bbm)=n \cdot \nabla^2 \cF_\TAP(\bbm)$, and denote
$\bD(\bbm)=\diag((1-\bbm^2)^{1/2})$. Then
\begin{equation}\label{eq:DHDbound}
\begin{aligned}
&~\bD(\bbm)\bH(\bbm)\bD(\bbm) \\
\succeq&~
-\frac{\lambda^2}{n}\bD(\bbm)\ones\ones^\top \bD(\bbm)
-\lambda \bD(\bbm)\bW\bD(\bbm)
+\id-\frac{2\lambda^2}{n}\bD(\bbm)\bbm\bbm^\top \bD(\bbm).
\end{aligned}
\end{equation}
Observe that if $\bbm \in \cM_q$, then $Q(\bbm) \geq
M(\bbm)^2>q^2>1-2\lambda^{-a}$. So
\[\|\bD(\bbm)\bbm\|_2^2 \leq
\|\bD(\bbm)\ones\|_2^2=\sum_{i=1}^n (1-m_i^2)=n(1-Q(\bbm))
<2n\lambda^{-a} \leq 2n\lambda^{-5}.\]
Observe also that by Markov's inequality,
\begin{equation}\label{eq:markovbound}
\frac{1}{n}\sum_{i=1}^n \ones\{(1-m_i^2)^{1/2}>\lambda^{-5/4}\}
\leq (1-Q(\bbm))\lambda^{5/2}<2\lambda^{-a+5/2}
\leq 2\lambda^{-5/2}.
\end{equation}
Then applying Lemma \ref{lem:sandwich_operator_bound}, with probability
$1-Ce^{-cn}$, $\|\bD(\bbm)\bW\bD(\bbm)\|_\op \leq C'\lambda^{-5/4}$ for every
$\bbm \in \cM_q$.
Applying these bounds to (\ref{eq:DHDbound}), for all sufficiently large
$\lambda$ and any $\bbm \in \cM_q$,
we have $\bD(\bbm)\bH(\bbm)\bD(\bbm) \succ \id/2$,
and hence $\bH(\bbm) \succ \bD(\bbm)^{-2}/2=\diag(1/(1-\bbm^2))/2 \succeq
\id/2$.
\end{proof}

\subsection{Analysis of NGD}\label{appendix:analysis_NGD}

We prove Lemma \ref{lem:NGD}, Theorem \ref{thm:localconvergence}, and
Theorem \ref{thm:globalconvergence}(b) on the convergence of NGD.

\begin{proof}[Proof of Lemma \ref{lem:NGD}]
We use the mirror-descent form (\ref{eq:mirrordescent}) for the NGD algorithm,
and adapt the argument of \cite[Theorem 3.1]{lu2018relatively}.

Recall the form of $\nabla^2 \cF_\TAP(\bbm)$ in (\ref{eq:hess}), and note
also that $\nabla^2 (-\bH(\bbm))=n^{-1}\diag(1/(1-\bbm^2))$. When
$\|\bW\|_\op<3$, we then have
\[\nabla^2 \cF_\TAP(\bbm) \prec C \cdot \nabla^2 (-\bH(\bbm)), \qquad
\nabla^2 \cF_\TAP(\bbm) \succ \nabla^2 (-\bH(\bbm))-C' \cdot \id\]
for some $\lambda$-dependent constants $C,C'>0$. For $\bbm \in (-1,1)^n \cap
\ball_{\sqrt{\eps n}}(\bbm_\star)$, we have also $n \cdot \nabla^2
\cF_\TAP(\bbm) \succ t$ by assumption,
so taking a suitable linear combination of these two
lower bounds yields
\[\nabla^2 \cF_\TAP(\bbm) \succ \mu \cdot \nabla^2 (-\bH(\bbm))\]
for a $\lambda$-dependent constant $\mu>0$.
Then by \cite[Proposition 1.1]{lu2018relatively}, these imply the relative
strong smoothness
\begin{equation}\label{eq:smooth}
\begin{aligned}
\cF_\TAP(\bbm) \leq&~ \cF_\TAP(\bbm')
+\nabla \cF_\TAP(\bbm')^\top(\bbm-\bbm')+C \cdot D_{-H}(\bbm,\bbm')\\
&\quad \text{ for all } \bbm,\bbm' \in (-1,1)^n
\end{aligned}
\end{equation}
and the relative strong convexity
\begin{equation}\label{eq:convex}
\begin{aligned}
\cF_\TAP(\bbm) \geq&~ \cF_\TAP(\bbm')+\nabla \cF_\TAP(\bbm')^\top(\bbm-\bbm')
+\mu \cdot D_{-H}(\bbm,\bbm') \\
&~\text{ for all } \bbm,\bbm' \in (-1,1)^n \cap
\ball_{\sqrt{\eps n}}(\bbm_\star).
\end{aligned}
\end{equation}

Let us choose the inverse step-size $L=1/\eta$ to be larger than this
constant $C$ in (\ref{eq:smooth}). By \cite[Lemma 3.1]{lu2018relatively}, the
minimizer $\bbm^{k+1}$ of (\ref{eq:mirrordescent}) satisfies the three-point
inequality
\[\begin{aligned}
&~\nabla \cF_\TAP(\bbm^k)^\top \bbm^{k+1}+LD_{-H}(\bbm^{k+1},\bbm^k)
+LD_{-H}(\bbm,\bbm^{k+1}) \\
\leq&~ \nabla \cF_\TAP(\bbm^k)^\top \bbm+LD_{-H}(\bbm,\bbm^k)
\end{aligned}\]
for any $\bbm \in (-1,1)^n$.
Then, applying (\ref{eq:smooth}) and this inequality,
\begin{align}
&\cF_\TAP(\bbm^{k+1})\nonumber\\ &\leq \cF_\TAP(\bbm^k)
+\nabla \cF_\TAP(\bbm^k)^\top (\bbm^{k+1}-\bbm^k)+LD_{-H}(\bbm^{k+1},\bbm^k)
\nonumber\\
&\leq \cF_\TAP(\bbm^k)
+\nabla \cF_\TAP(\bbm^k)^\top(\bbm-\bbm^k)+LD_{-H}(\bbm,\bbm^k)
-LD_{-H}(\bbm,\bbm^{k+1}).\label{eq:valueinduction}
\end{align}
Taking $\bbm=\bbm^k$ shows in particular that
$\cF_{\TAP}(\bbm^{k+1}) \leq \cF_{\TAP}(\bbm^k)$.

Next, we show that for sufficiently large $L$, every iterate $\bbm^k$ satisfies
\begin{equation}\label{eq:induction}
\|\bbm^k-\bbm_\star\|<\sqrt{\eps n}, \qquad
\cF_\TAP(\bbm^k)<\cF_\TAP(\bbm_\star)+t\eps/8.
\end{equation}
We induct on $k$, where the base case $k=0$ holds by assumption. Suppose that
(\ref{eq:induction}) holds for $k$. Then the above shows $\cF_{\TAP}(\bbm^{k+1})
\leq \cF_{\TAP}(\bbm^k)<\cF_\TAP(\bbm_\star)+t\eps/8$.
Observe that for any $\rho<\eps$, by the strong convexity
$\nabla^2 \cF_\TAP(\bbm) \succ (t/n)\id$ for
$\bbm \in (-1,1)^n \cap \ball_{\sqrt{\eps n}}(\bbm_\star)$, we have the
implication
\begin{equation}\label{eq:valuebound}
\begin{aligned}
\bbm \in (-1,1)^n \cap \ball_{\sqrt{\eps n}}(\bbm_\star) \text{ and }
\cF_{\TAP}(\bbm)<\cF_\TAP(\bbm_\star)+t\rho/2\\
\Rightarrow \|\bbm-\bbm_\star\|_2<\sqrt{\rho n}.
\end{aligned}
\end{equation}
So (\ref{eq:induction}) in fact implies $\|\bbm^k-\bbm_\star\|_2<\sqrt{\eps n}/2$.
Comparing the value of (\ref{eq:mirrordescent}) at $\bbm=\bbm^k$ and at the
minimizer $\bbm=\bbm^{k+1}$, we have
\[\nabla \cF_\TAP(\bbm^k)^\top(\bbm^{k+1}-\bbm^k)
+LD_{-H}(\bbm^{k+1},\bbm^k) \leq 0.\]
From the definition of the Bregman divergence $D_{-H}(\bbm,\bbm')$ in
(\ref{eq:Bregman}), we have
\begin{align*}
&\Big|\nabla \cF_\TAP(\bbm^k)^\top
(\bbm^{k+1}-\bbm^k)+D_{-H}(\bbm^{k+1},\bbm^k)\Big|\\
&=\Big|\big(\nabla \cF_\TAP(\bbm^k)+\nabla H(\bbm^k)\big)^\top(\bbm^{k+1}
-\bbm^k)-H(\bbm^{k+1})+H(\bbm^k)\Big| \\
&\leq C\left(\frac{\|\bbm^{k+1}-\bbm^k\|_2}{\sqrt{n}}+1\right).
\end{align*}
Also, taking $L>1$ and applying the strong convexity
$\nabla^2 (-H(\bbm)) \succ (1/n)\id$, we have
\[(L-1)D_{-H}(\bbm^{k+1},\bbm^k) \geq \frac{L-1}{2n}\|\bbm^{k+1}-\bbm^k\|_2^2.\]
Combining these bounds,
\[0 \geq \frac{L-1}{2n}\|\bbm^{k+1}-\bbm^k\|_2^2
-\frac{C}{\sqrt{n}}\|\bbm^{k+1}-\bbm^k\|_2-C\]
which implies
\[\|\bbm^{k+1}-\bbm^k\|_2<\frac{C'}{\sqrt{L-1}} \cdot \sqrt{n}.\]
for a $\lambda$-dependent constant $C'>0$. Then, taking $L$ large enough so that
$C'/\sqrt{L-1}<\sqrt{\eps}/2$, we obtain
\[\|\bbm^{k+1}-\bbm_\star\|_2
\leq \|\bbm^{k+1}-\bbm^k\|_2+\|\bbm^k-\bbm_\star\|_2<\sqrt{\eps n},\]
completing the induction and the proof of (\ref{eq:induction}).

The first statement of (\ref{eq:induction}) allows us to apply the relative
strong convexity (\ref{eq:convex}) at $\bbm^k$ and $\bbm_\star$, to obtain
\[\cF_\TAP(\bbm^k)+\nabla \cF_\TAP(\bbm^k)^\top (\bbm_\star-\bbm^k)
\leq \cF_\TAP(\bbm_\star)-\mu \cdot D_{-H}(\bbm_\star,\bbm^k).\]
Applying this bound to (\ref{eq:valueinduction}) now with $\bbm=\bbm_\star$,
\[\cF_\TAP(\bbm^{k+1}) \leq \cF_\TAP(\bbm_\star)
+\big(L-\mu\big) D_{-H}(\bbm_\star,\bbm^k)-LD_{-H}(\bbm_\star,\bbm^{k+1}).\]
Multiplying by $(\frac{L}{L-\mu})^{k+1}$ and summing over $k$ to
telescope the sums of the last two terms,
\[\sum_{j=0}^{k-1} \left(\frac{L}{L-\mu}\right)^{j+1}\cF_\TAP(\bbm^{j+1})
\leq  \sum_{j=0}^{k-1} \left(\frac{L}{L-\mu}\right)^{j+1}
\cF_\TAP(\bbm_\star)+LD_{-H}(\bbm_\star,\bbm^0).\]
Now applying $\cF_\TAP(\bbm^{j+1}) \geq \cF_\TAP(\bbm^k)$ for all $j \leq k-1$
to the left side, we obtain
\[\cF_\TAP(\bbm^k) \leq \cF_\TAP(\bbm_\star)
+L\left(\sum_{j=0}^{k-1} \left(\frac{L}{L-\mu}\right)^{j+1}\right)^{-1}
D_{-H}(\bbm_\star,\bbm^0).\]
From the form of $D_{-H}(\bbm_\star,\bbm^0)$ in (\ref{eq:Bregman}),
on the event $\|\bW\|_\op<3$, we have
\[D_{-H}(\bbm_\star,\bbm^0)<
C\left(1+\frac{\|\arctanh(\bbm^0)\|_2}{\sqrt{n}}\right)\]
for a $\lambda$-dependent constant $C>0$. Then the above shows
\[\cF_\TAP(\bbm^k) \leq
\cF_\TAP(\bbm_\star)+L\left(\frac{L-\mu}{L}\right)^k \cdot C
\left(1+\frac{\|\arctanh(\bbm^0)\|_2}{\sqrt{n}}\right).\]
Identifying $(L-\mu)/L=1-\mu\eta$, this shows (\ref{eq:valueconvergence}).
Combining this with (\ref{eq:induction}) and (\ref{eq:valuebound}) shows
(\ref{eq:mconvergence}).
\end{proof}

\begin{proof}[Proof of Theorem \ref{thm:localconvergence}]
Suppose the initialization has sign $\langle \bx,\bh^0 \rangle \geq 0$.
By Theorem \ref{thm:local}(b), for some $\eps,t>0$, we have with probability
approaching 1 that $\lambda_{\min}(n \cdot \nabla^2 \cF_\TAP(\bbm))>t$ for all
$\bbm \in (-1,1)^n \cap \ball_{\sqrt{\eps n}}(\bbm_\star)$.
By Lemma \ref{lem:AMP}, for a sufficiently large iteration $T$, the iterate
$\bbm^T$ of AMP will satisfy
\[\cF_\TAP(\bbm^T)<e_\star+t\eps/16, \quad
\frac{1}{n^2}\|\widehat{\bX}_\Bayes-\bbm^T(\bbm^T)^\top\|_\sF^2<\eps^4/8,
\quad \frac{\|\arctanh \bbm^k\|_2^2}{n}<C\]
for a $\lambda$-dependent constant $C>0$. From Corollary \ref{cor:critexists}
and Theorem \ref{thm:local}(a), the local minimizer $\bbm_\star$ satisfies
\[\cF_\TAP(\bbm_\star)>e_\star-t\eps/16, \qquad
\frac{1}{n^2}\|\widehat{\bX}_\Bayes-\bbm_\star\bbm_\star^\top\|_\sF^2<\eps^4/8.\]
Therefore $\cF_\TAP(\bbm^T)<\cF_\TAP(\bbm_\star)+t\eps/8$. Recalling
(\ref{eq:pmbound}),
\[\frac{1}{n^2}\|\bbm^T-\bbm_\star\|_2^4
\leq \frac{2}{n^2}\|\bbm^T(\bbm^T)^\top-\bbm_\star\bbm_\star^\top\|_\sF^2
< \eps^2.\]
This verifies that with high probability, the conditions of
Lemma \ref{lem:NGD} hold for initializing NGD at this AMP iterate $\bbm^T$,
so the theorem follows from Lemma \ref{lem:NGD}.
If instead $\langle \bx,\bh^0 \rangle<0$, then the theorem holds with
sign $-\bbm_\star$ by sign symmetry.
\end{proof}

\begin{proof}[Proof of Theorem \ref{thm:globalconvergence}(b)]
Suppose the initialization has sign $\langle \bx,\bh^0\rangle \geq 0$.
For $\bbm^0=\tanh(\bh^0)$, applying Lemma \ref{lem:AMP}(a) with $k=0$ and
$\psi(x,y)=x\tanh(y)$, we have almost surely
\[\lim_{n \to \infty} \langle \bx,\bbm^0 \rangle/n
=\E_{G \sim \cN(0,1)}[\tanh(\lambda^2-1+\sqrt{\lambda^2-1} \cdot G)].\]
For large $\lambda$, applying monotonicity of $\tanh$ and the bounds
$\P[G<-t]<e^{-t^2/2}$ and $\tanh(t)>1-2e^{-2t}$, we may bound this by
\begin{align*}
&~\E_{G \sim \cN(0,1)}[\tanh(\lambda^2-1+\sqrt{\lambda^2-1} \cdot G)]\\
>&~\tanh(\lambda^2/3)\P\left[G \geq -\lambda/2\right]
-\P\left[G<-\lambda/2\right]>1-3e^{-\lambda^2/8}.
\end{align*}
Then for any $\eps>0$, with probability approaching 1,
$\|\bx-\bbm^0\|_2^2/n \leq 2-2\langle \bx,\bbm^0 \rangle/n
<6e^{-\lambda^2/8}+\eps$. Similarly by Corollary \ref{cor:critexists} and
Proposition \ref{prop:qstarlargelambda}, also with probability approaching 1,
$\|\bx-\bbm_\star\|_2^2/n \leq 2-2q_\star+\eps<2e^{-\lambda^2/8}+\eps$.
Choosing $\eps=e^{-\lambda^2/8}$ and combining these bounds, for all
sufficiently large $\lambda$,
\begin{equation}\label{eq:minitclose}
\frac{\|\bbm^0-\bbm_\star\|_2^2}{n}
\leq 2\frac{\|\bx-\bbm^0\|_2^2}{n}+2\frac{\|\bx-\bbm_\star\|_2^2}{n}
<e^{-\lambda^2/9}.
\end{equation}
Then applying Proposition \ref{prop:TAPcontinuous} for sufficiently large $\lambda$, with probability approaching
1,
\begin{equation}\label{eq:Finitclose}
\cF_\TAP(\bbm^0)-\cF_\TAP(\bbm_\star)<2e^{-\lambda^2/36}.
\end{equation}

On the event where (\ref{eq:minitclose}) and (\ref{eq:Finitclose}) both hold
and where the conclusion of Lemma \ref{lem:quantitativelandscape}(b) holds
for $a=5$, for sufficiently large $\lambda$, this initialization $\bbm^0$
satisfies the conditions of Lemma \ref{lem:NGD} with $\eps=e^{-\lambda^2/18}$
and $t=1/2$. Then Theorem \ref{thm:globalconvergence}(b) holds
by Lemma \ref{lem:NGD} and the initial condition
$\|\arctanh(\bbm^0)\|_2/\sqrt{n}=\|\bh^0\|_2/\sqrt{n}\rightarrow \sqrt{\lambda^2(\lambda^2-1)}$. If instead $\langle \bx,\bh^0 \rangle<0$, then
Theorem \ref{thm:globalconvergence}(b) holds with $-\bbm_\star$ by sign
symmetry.
\end{proof}

\subsection{Analysis of AMP}\label{appendix:AMPglobal}

We prove Lemma \ref{lem:AMPcontractive} on the contractivity of the AMP map, and Theorem \ref{thm:globalconvergence}(a) on the convergence of AMP.
In Remark \ref{rmk:motivation},
we provide more motivation for the reparameterization we have chosen.
This remark can be read prior to reading the proof.

\begin{proof}[Proof of Lemma \ref{lem:AMPcontractive}]
Denote $\bbm_+=\Lambda(\bp_+)$, $\bbm=\Lambda(\bp)$, $\bbm_-=\Lambda(\bp_-)$,
and $\bD(\bbm)=\diag((1-\bbm^2)^{1/2})$. Then, applying the form
of $\de T_\AMP$ in (\ref{eq:dTAMP}) and the identity $\de \bbm/\de
\bp=(1-\bbm^2)^{1/2}$,
\[\de T_\AMP^{(p)}(\bp,\bp_-)=\begin{pmatrix}
\bD(\bbm_+) \cdot [\lambda \bY+2\lambda^2\bbm_-\bbm^\top/n] \cdot \bD(\bbm)
&-\bD(\bbm_+) \cdot \lambda^2[1-Q(\bbm)] \cdot \bD(\bbm_-) \\
\id & \bzero \end{pmatrix}.\]

We first use a crude bound on the operator norms of the upper blocks:
Consider the event of probability $1-e^{-cn}$ where $\|\bW\|_\op<3$. On this
event, applying $\|\bD(\bbm)\|_\op \leq 1$ and $\|\bbm\|_2,\|\ones\|_2 \leq
\sqrt{n}$, we have
\[
\begin{aligned}
&~\Big\|\bD(\bbm_+) \cdot [\lambda \bY + 2\lambda^2 \bbm_{-} \bbm^\sT/n ]
\cdot \bD(\bbm)\Big\|_\op<3\lambda+3\lambda^2, \\
&~\Big\|{-}\bD(\bbm_+) \cdot \lambda^2[1-Q(\bbm)] \cdot
\bD(\bbm_-)\Big\|_\op \leq \lambda^2.
\end{aligned}\]
Integrating these bounds along the linear path between $(\bp,\bp_-)$ and
the AMP fixed point $(\bp_\star,\bp_\star)$ then yields
\begin{equation}\label{eq:m1-close-to-m*}
\| \bp_+ - \bp_\star\|_2 \leq (3\lambda+3\lambda^2) \|\bp-\bp_\star\|_2
+\lambda^2 \|\bp_{-}-\bp_\star\|_2<7\lambda^{-5}\sqrt{n}
\end{equation}
when $\bp,\bp_- \in \ball_{\lambda^{-7}\sqrt{n}}(\bp_\star)$.
Since $\bbm=\Lambda(\bp)$ is 1-Lipschitz, we then have
$\|\bbm-\bbm_\star\|_2,\|\bbm_{-}-\bbm_\star\|_2<\lambda^{-7}\sqrt{n}$
and $\|\bbm_+-\bbm_\star\|_2<7\lambda^{-5}\sqrt{n}$. Then, applying the
assumption $M(\bbm_\star)>1-\lambda^{-5}$ so that
$Q(\bbm_\star) \geq M(\bbm_\star)^2>1-2\lambda^{-5}$, and applying also
the continuity bound for $Q(\bbm)$ in Proposition \ref{prop:TAPcontinuous},
this shows $Q(\bbm_+),Q(\bbm),Q(\bbm_-)>1-12\lambda^{-5}$ for all large
sufficiently $\lambda$.

Now applying $Q(\bbm)>1-12\lambda^{-5}$ above, we obtain the improved bound
\[\Big\|{-}\bD(\bbm_+) \cdot \lambda^2[1-Q(\bbm)] \cdot
\bD(\bbm_-)\Big\|_\op<12\lambda^{-3}.\]
Also $\|\bD(\bbm)\bx\|_2^2=\sum_i (1-m_i^2)x_i^2 \leq
n(1-Q(\bbm))<12\lambda^{-5}n$, and the same bounds hold for
$\bD(\bbm_+)\bx$, $\bD(\bbm)\bbm_-$, and $\bD(\bbm)\bbm$. Then
\[\Big\|\bD(\bbm_+) \cdot [\lambda \bY + 2\lambda^2 \bbm_{-} \bbm^\sT/n ]
\cdot \bD(\bbm)\Big\|_\op
\leq \lambda \|\bD(\bbm_+) \bW\bD(\bbm)\|_\op+36\lambda^{-3}.\]
On the event of probability $1-Ce^{-cn}$ (for $\lambda$-dependent $C,c>0$)
where Lemma \ref{lem:sandwich_operator_bound} holds
with $\eps=\lambda^{-5/4}$ and $s=12\lambda^{-5/2}$,
applying Markov's inequality and this lemma as in (\ref{eq:markovbound}),
we obtain for sufficiently large $\lambda$ that
\[\Big\|\bD(\bbm_+) \cdot [\lambda \bY + 2\lambda^2 \bbm_{-} \bbm^\sT/n ]
\cdot \bD(\bbm)\Big\|_\op<\lambda^{-1/5}.\]

Finally, integrating again these improved bounds along the linear path from
$(\bp,\bp_-)$ to $(\bp_\star,\bp_\star)$, we may obtain
\[\|\bp_+ - \bp_\star\|_2 \leq \lambda^{-1/5}\|\bp-\bp_\star\|_2
+2\lambda^{-2/5}\|\bp_{-}-\bp_\star\|_2.\]
In particular, for large $\lambda$ this implies
$\|\bp_+ - \bp_\star\|_2 \leq \max(\|\bp - \bp_\star\|_2,
\|\bp_{-}-\bp_\star\|_2)$, so $\bp_+ \in \ball_{\lambda^{-7}\sqrt{n}}(\bp_\star)
\cap \Omega^{(p)}$. Adding $\lambda^{-1/5}\|\bp-\bp_\star\|_2$ to both sides above
yields (\ref{eq:AMPcontractive}).
\end{proof}

\begin{remark}
\label{rmk:motivation}
To motivate the reparameterization by $\bp$,
consider instead an analysis of the contractivity of the AMP map without this reparameterization. 
	Denote $(\bbm_+,\bbm) = T_{\AMP}(\bbm,\bbm_-)$.
	The Jacobian of $T_{\AMP}$ is given by
	\begin{equation}
		\de T_{\AMP}(\bbm,\bbm_-)
			=
			\begin{pmatrix}
				\diag(1 - \bbm_+^2)
				\cdot [\lambda \bY + 2\lambda^2 \bbm_- \bbm^\top / n]
				&\quad 
				-\diag(1-\bbm_+^2)\lambda^2(1 - Q(\bbm)) \\
				\id & 0 
			\end{pmatrix}.
	\end{equation}
	We may use the key fact that $1 - \bbm_+^2$ is close to zero in most
coordinates to show that most rows of the upper two blocks are small.
	This is not enough, however, to establish a strong operator norm bound
for these blocks, and we would like to also have that most of the columns of the
upper two blocks are small. The reparameterization by $\bp$ is chosen so that
these two blocks become right-multiplied also by the small factors of
$\diag(1-\bbm^2)^{1/2}$.

	In more detail, for any reparametrization $\bbm = \Lambda(\bp)$,
	denoting $\bp_+ = \Lambda(\bbm_+)$ and $\bp_- = \Lambda(\bbm_-)$,
the upper two blocks of the Jacobian in this parametrization become
\begin{align*}
&\frac{\de \bp_+}{\de \bbm_+} \diag(1 - \bbm_+^2) \cdot [\lambda \bY
+ 2\lambda^2 \bbm_- \bbm^\top / n] \frac{\de \bbm}{\de \bp},\\
& - \frac{\de \bp_+}{\de \bbm_+} \diag(1 - \bbm_+^2) \lambda^2(1-Q(\bbm)) 
			\frac{\de \bbm_-}{\de \bp_-}.
\end{align*}
	Because $\frac{\de \bp}{\de \bbm} = \Big(\frac{\de \bbm}{\de \bp}\Big)^{-1}$,
	we see that a natural choice of reparametrization is to ensure
$\frac{\de \bp}{\de \bbm} = 1/\sqrt{1 - \bbm^2}$, which holds exactly for our
definition of $\bp$.
	\end{remark}

\begin{proof}[Proof of Theorem \ref{thm:globalconvergence}(a)]
Suppose $\bx=\ones$ and $\langle \bx,\bh^0 \rangle \geq 0$.
By Corollary \ref{cor:critexists} and Proposition \ref{prop:qstarlargelambda},
for all sufficiently large $\lambda$, with probability approaching 1,
$\bbm_\star \in \cM_{1-\lambda^{-5}}$. Also
from the proof of Theorem \ref{thm:local}(b) and the bound (\ref{eq:pmbound}),
for any $\delta>0$ sufficiently small and $k=k(\delta)$ sufficiently large (where ``sufficiently small'' and ``sufficiently large'' depend on $\lambda$),
with probability approaching 1, the AMP iterates
$\bbm^{k-1}$ and $\bbm^k$ satisfy
\[\| \bbm^{k-1}-\bbm_\star\|_2<\delta\sqrt{n}, \qquad
\|\bbm^k-\bbm_\star\|_2<\delta\sqrt{n}.\]
Because $\Gamma$ is continuous on the compact domain $[-1,1]$,
it has a modulus of continuity: $|\Gamma(m) - \Gamma(m')| < \eps(\Delta)$
whenever $|m-m'| < \Delta$, where $\eps(\Delta) \rightarrow 0$ as $\Delta
\rightarrow 0$. Then by Markov's inequality and the fact that $\bp^k,\bp_\star
\in \Omega^{(p)}=(-\pi/2,\pi/2)^n$,
	\[\begin{aligned}\frac1n\| \bbm^k - \bbm_\star \|_2^2 < \delta^2
		&\;\;\text{implies}\;\;\\
		\frac1n\| \bp^k - \bp_\star \|_2^2
		&\leq \eps(\sqrt{\delta})^2 + \pi^2 \cdot
\frac1n| \{ i : |m_i^k - m_{\star,i} |^2
\geq \delta\}| < \eps(\sqrt{\delta})^2 + \pi^2\delta.\end{aligned}\]
The same applies to $\bp^{k-1}$. Thus, choosing $\delta$ sufficiently small,
we ensure $\bp^k,\bp^{k-1} \in \Omega^{(p)} \cap
\ball_{\lambda^{-7}\sqrt{n}}(\bp_\star)$.
Then on the event where also Lemma \ref{lem:AMPcontractive} holds,
we conclude 
\begin{equation}
\begin{aligned}
\|\bp^{k+r}-\bp_\star\|_2 \leq&~
\|(\bp^{k+r},\bp^{k+r-1}) - (\bp_\star,\bp_\star)\|_\lambda \leq
(2\lambda^{-1/5})^r\| (\bp^k,\bp^{k-1}) - (\bp_\star,\bp_\star)\|_\lambda\\
<&~(2\lambda^{-1/5})^r\lambda^{-7}\sqrt{n}.
\end{aligned}
\end{equation}
Noting that $k$ is a $\lambda$-dependent constant, and
choosing $\lambda$-dependent constants $C,\alpha>0$ to also account for the
first $k$ iterations, we obtain for every $r \geq 1$ that
\[\|\bbm^r-\bbm_\star\|_2<C\alpha^r \sqrt{n}.\]
Then also with high probability, for modified constants $C,\alpha>0$, we have
$\cF_\TAP(\bbm^k)-\cF_\TAP(\bbm_\star)<C\alpha^k$ by Proposition
\ref{prop:TAPcontinuous}. If $\langle \bx,\bh^0 \rangle<0$, the same statements
hold with $-\bbm_\star$ by sign symmetry.
\end{proof}

\section{Numerical evaluation of eigenvalues of the linearized AMP operator }

\begin{figure}
\centering
\includegraphics[width = 0.44\linewidth]{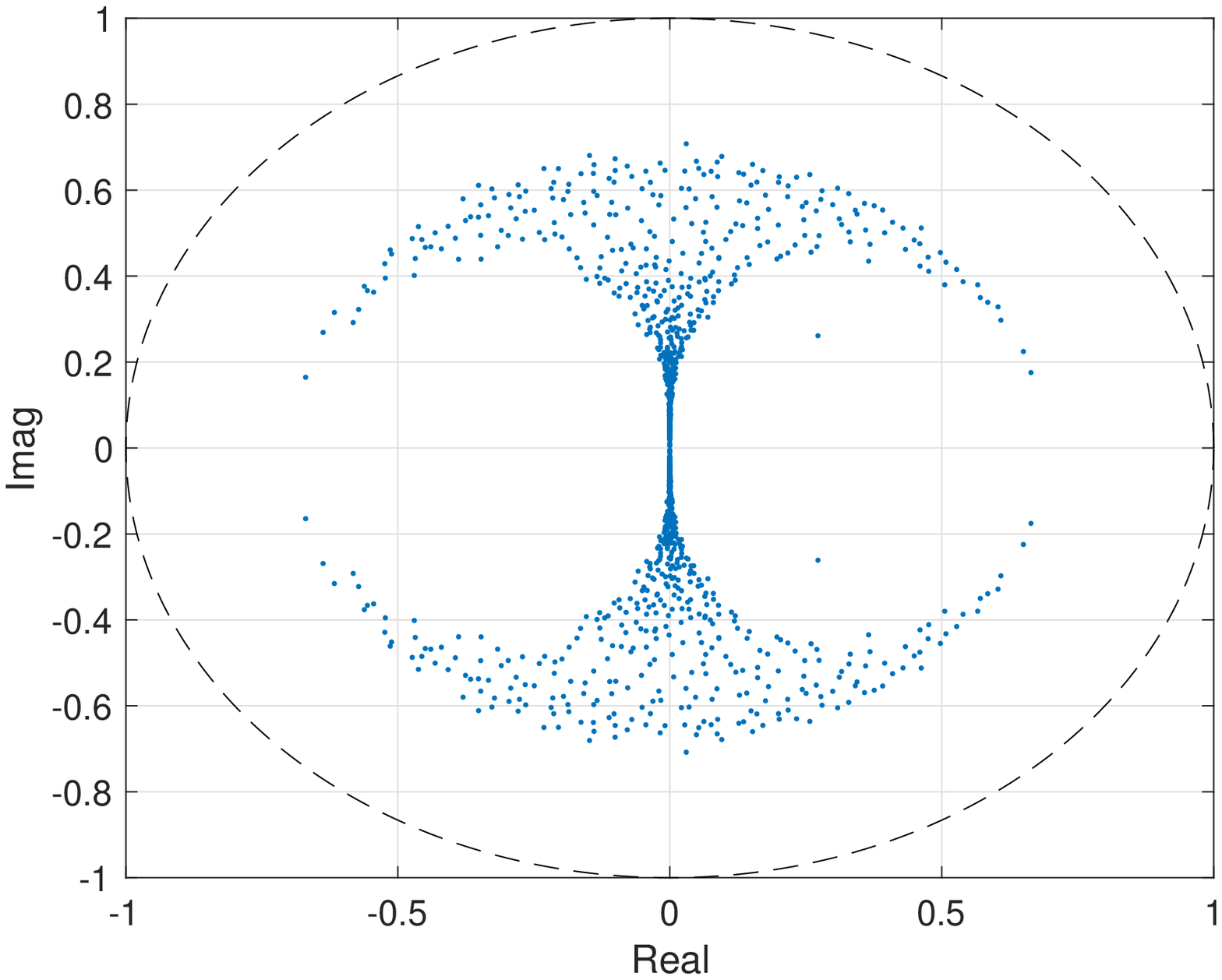}
\caption{The scatter plot of eigenvalues of the linearized AMP operator $\de
T_\AMP(\bbm_\star,\bbm_\star)$. We choose $n = 500$ and $\lambda = 1.5$. }\label{fig:AMP_Jacobian}
\end{figure}

Theorem \ref{thm:local} shows that, for any $\lambda > 1$, the spectral radius
of the Jacobian of the AMP map $\de T_{\AMP}(\bbm_\star,\bbm_\star)$ will be
bounded away from $1$ with high probability. Figure \ref{fig:AMP_Jacobian} shows
a scatter plot of all eigenvalues of $\de T_{\AMP}(\bbm_\star,\bbm_\star)$ for a
specific instance $\bY$, with $n = 500$ and $\lambda = 1.5$.  

\end{document}